\documentclass{article}
\usepackage{verbatim} 
\usepackage[latin1]{inputenc}
\usepackage[german,english]{babel}
\usepackage{mathbbol,mathrsfs,amsfonts,amssymb,amsthm,stmaryrd}
\usepackage[intlimits]{amsmath}

\usepackage[all]{xypic}
\usepackage{mathpazo}
\usepackage[mathcal]{euscript}

\usepackage{graphics,supertabular}

\newcommand{\C}{\mathbb{C}}
\newcommand{\R}{\mathbb{R}}
\newcommand{\N}{\mathbb{N}}
\newcommand{\Z}{\mathbb{Z}}
\newcommand{\T}{\mathbb{T}}
\newcommand{\A}{\mathcal{A}}
\newcommand{\U}{{\rm U}}
\newcommand{\UA}{{\rm U}^{\rm ab}}
\newcommand{\PU}{{\rm PU}}
\newcommand{\PP}{{\rm P}}
\newcommand{\HH}{\mathcal{H}}
\newcommand{\K}{\mathcal{K}}
\newcommand{\LL}{\mathcal{L}}
\newcommand{\B}{\mathcal{B}}
\newcommand{\E}{\mathcal{E}}
\newcommand{\M}{\mathcal{M}}
\newcommand{\G}{\mathcal{G}}
\newcommand{\nou}{ {u\!\!\!\!\! \bigcirc}}
\newtheoremstyle{ansgarthmstyle}{}{}{\it}{}{\bf}{}{0.4cm}{}
\theoremstyle{ansgarthmstyle}
\newtheorem{exa}{Example}[section]

\newtheorem{defi}{Definition}[section]
\newtheorem{prop}{Proposition}[section]
\newtheorem{rem}{Remark}[section]
\newtheorem{lem}{Lemma}[section]
\newtheorem{thm}{Theorem}[section]
\newtheorem{cor}{Corollary}[section]
\newenvironment{pf}{\begin{proof}[\textrm{{\bf Proof\ :\!}}]} {\end{proof}}

\begin{document}
\selectlanguage{german}

\thispagestyle{empty}
\begin{center}
\begin{tabular}{c}
\huge {\sc Die lokale Struktur}  \\\\
\huge {\sc von  }\\\\
\huge  {\sc T-Dualitätstripeln}\\\\
\\\\\\\\\\\\\
\large Dissertation zur Erlangung des Doktorgrades\\\\\\
\large der Mathematisch-Naturwissenschaftlichen Fakultäten  \\\\\\
\large der Georg-August-Universität   Göttingen\\\\\\
\large  angefertigt von \\\\\\

\LARGE {\sc Ansgar Schneider}\\
\\\\
\large aus\\\\\\
\LARGE {\sc Foerde}\\\\\\
\large geboren in \\\\\\
\LARGE {\sc Siegen}
\\\\\\\\
\large Göttingen, im Oktober 2007
\end{tabular}
\end{center}

\newpage
\thispagestyle{empty}
\ 
\vspace{15cm}
\ 
\\
D7
\\\\
Referrent: Prof. Dr. Ulrich Bunke	
\\\\
Koreferent: Prof. Dr. Thomas Schick
\\\\
Tag der mündlichen Prüfung: 5.11. 2007

\newpage

\thispagestyle{empty}

\ 
\vspace{5cm}

\begin{center}

\large{
Meinen Eltern und Großeltern 

--  Jetzt habt ihr den Salat!}

\end{center}

\newpage
\thispagestyle{empty}
\ 

\newpage
\section*{Vorrede und Danksagung}
 
Die vorliegende Arbeit ist das wesentliche Ergebnis meiner
Zeit am mathematischen  Institut der Universität Göttingen. Dort war ich vom Sommer 2005 bis zum Herbst 2007 
Promotionsstudent, unterstützt  durch den  dortigen 
Graduiertenkolleg Gruppen und Geometrie.
\thispagestyle{empty}

Wie der Titel sagt, beschäftigt sich die Arbeit mit der lokalen Struktur von T-Du\-ali\-täts\-tri\-peln.
T-Dualitätstripel (das T steht für Torus) sind mathematische Objekte, die sich als zweckmäßig 
erwiesen haben, T-Dualität mittels topologischer Methoden  zu beschreiben \cite{BS}.
T-Dualität selbst ist eine Dualität von Stringtheorien \cite{Po} 
auf zwei verschiedenen Raumzeitmanigfaltigkeiten.
Im einfachsten Falle sind diese 
durch eine Transformation 
Radius $\mapsto$ 1/Radius entlang einer 
kompaktifizierten Raumdimension (einem Torus)
miteinander identifiziert.
Das Verhalten der zugrundeliegenden 
Felder und die Geometrie und Topologie der Raumzeitmanigfaltigkeiten
unter dem Übergang von einer Raumzeit zu der anderen ist das, 
was durch T-Dualität beschrieben wird.
Die Literatur zu diesem Thema ist reichhaltig, und man möge etwa
\cite{BEM,BHM1,MR,BS} und der darin zitierten Literatur folgen, um 
einen Überblick darüber zu erhalten.
Neben der schon erwähnten topologischen Beschreibungsweise,
ist  in  \cite{MR}  ein $C^*$-algebraischer Zugang beschrieben, 
dessen zentrale Objekte  gewisse $C^*$-dynamische Systeme sind, deren 
Dualitätstheorie wiederum eine andere mathematische Beschreibung von 
T-Du\-ali\-tät liefert.

Das Ziel dieser Arbeit ist es, einen expliziten Zusammenhang zwischen dem $C^*$-algebraischen Zugang und den Resultaten über topologische T-Dualität herzustellen. Dabei wird sich zeigen, daß wir dieses Ziel erreichen können, indem wir die lokale Struktur der zugrundeliegenden Objekte analysieren und mit Hilfe
der gewonnenen lokalen Daten zeigen, daß in beiden Fällen 
die jeweiligen, geeignet gewählten  Äquivalenzklassen 
der topologischen und $C^*$-algebraischen Objekte übereinstimmen.
\\

Ich will mich an dieser Stelle recht herzlich bei all denjenigen bedanken, die zum Gelingen und Entstehen dieser Arbeit beigetragen haben.
An erster Stelle gilt mein Dank natürlich meinem Doktorvater Herrn Ulrich Bunke, der sich (erstaunlicherweise) dazu bereit erklärte, mich als Doktorand zu betreuen, und mich immer
wieder auf vielfältige Weise gefordert und gefördert hat. 
Ohne ihn wäre diese Arbeit tat\-säch\-lich nicht möglich gewesen, und 
da die mir aufgetragene Fragestellung im Grenzgebiet verschiedener, mathematischer Disziplinen liegt, konnte ich Dinge lernen, die mir unter anderen Umständen sicher verwehrt geblieben wären.  
Meinem Doktoronkel Herrn Thomas Schick und meinem Doktorbruder Herrn Moritz Wiethaup sei hier ausdrücklich für die zahlreichen, fruchtbaren Gespräche gedankt, ohne die ich einige Dinge nicht hätte so einfach oder schnell oder überhaupt hätte bewerkstelligen können. 
Dem Graduiertenkolleg Gruppen und Geometrie danke ich für das entgegengebrachte Vertrauen und die großzügige Unterstützung in den letzten zwei Jahren.

Nachdem mir von prominenter Stelle zugetragen wurde, daß es wohl angebracht ist,
auch denjenigen zu danken, die indirekt an dieser Arbeit beteiligt sind, will ich
manchen von meinen Lehrern und Hochschullehrern danken, von denen ich
(natürlich in alphabetischer Folge) besonders die Herren
Helmut Becker,
Detlev Buchholz,
Rüdiger Heidersdorf,
Bernhard Meyer,
Karl-Henning Rehren
und  Wolfgang Watzlawek 
erwähnen will, die alle ihren besonderen Anteil an meinem Bildungsweg hatten und haben.
Herrn Friedrich von Schiller  danke ich für die Ode an die Freude,
Herrn Ralf Lindemann für alles und  Herrn Jens Latsch dafür, daß er mir im ersten Semester meine Oberonprogramme zugeschickt hat.
Zum Schluß will ich meinen lieben Eltern und meiner lieben Großmutter dafür danken, daß diese Arbeit, nicht nur in der offensichtlichen Weise, ohne sie nicht hätte entstehen können.
\vspace{1.5cm}

\begin{flushright}
Göttingen, im Oktober 2007

Ansgar Schneider
\end{flushright}
\thispagestyle{empty}

\selectlanguage{english}

\newpage
\tableofcontents

\newpage
\section{Introduction and Summary}
String theory \cite{Po} is based on the physical idea to describe nature not only by point particles but with the concept of higher dimensional objects.
Although it is still unclear how close string theory is to physics
a lot of interesting new phenomena have been observed 
which have led to many fruitfully new ideas and have made it a 
rich theory. In particular, a lot of new mathematical ideas  have 
arose from the desire to understand the discovered 
structures. One of those is the  concept  of T-duality.

T-duality is a duality of string theories (type IIA and IIB) on
different underlying space-time manifolds $E$ 
and\footnote{The accent circonflexe $\hat\ $ is going to be the most
overloaded symbol  in this work.}  $\widehat E$ which are (in the simplest case) related by a transformation of type: radius $\mapsto$ 1/radius along a compactified  space-time dimension. 
A duality between two theories on different space-time  manifolds gives a prescription how the fields and their correlation functions 
transform under the change of the underlying manifolds.  
In the present case one may take as an example the charges of the D-branes which 
take values in the twisted K-theory of $E$ \cite{BM}, where the twist is given by a background field on $E$ (a 3-form $H$ called H-flux).
Then T-duality must give the answer how the background fields transform and 
should lead to an isomorphism of the twisted  K-theories of the underlying  manifolds.

We cannot give a summary of the whole subject, but we can
try to point out some of its  mathematical issues. 
In literature there are  different approaches of a mathematical understanding of T-duality. One is based on the theory of $C^*$-dynamical systems 
which serve a notion of T-duality using crossed product $C^*$-algebras \cite{BHM2,MR},
another is by  geometric and topological means \cite{BEM,BHM1,BS1,BS},
a third using methods from algebraic geometry \cite{BSST}.
We focus our attention to the first and second approach and 
continue to describe some features of the geometric-topological side in more detail. 
\\

Let us think of the manifolds $E$ and $\widehat E$ as  principal circle bundles which have isomorphic quotients $E/S^1\cong\widehat E/S^1=:B$.
In \cite{BEM} it is described in terms of differential geometry how the data
of the curvature $F$ of $E$ and of the H-flux $H$ on $E$ are related to the corresponding dual data 
$\widehat F$ and $\widehat H$ of $\widehat E$. The result is that integration of  $H$ along the fibres of $E$ yields the dual curvature and vice versa.
In the case of $S^1$-bundles $E$ and $\widehat E$ we can identify the classes of the curvatures with the realifications of the first Chern classes $c_1,\hat c_1\in H^2(B,\Z)$ of the respective bundles, and there also exist integer cohomology classes $h\in H^3(E,\Z), \hat h\in H^3(\widehat E,\Z)$  whose realifications are $h_\R=[H]$ and $\hat h_\R=[\widehat H]$.
In this sense, we  forgot geometry and now may only consider these topological data. 
This is the point of view which was adopted in \cite{BS1},
wherein among other things the results of \cite{BEM} are restated on a purely topological level.
The higher dimensional case, where the circle $S^1$ is replaced 
by the $n$-dimensional torus $\T^n=(S^1)^{\times n}$, i.e. $E$ is thought of 
a principal $\T^n$-bundle, is described 
in \cite{BHM1} in terms of differential geometry. Its topological structure is 
described in \cite{BS} which we want to discuss in more detail.
They introduce so-called  T-duality triples and define that 
a pair $(E,h)$ is dual to a pair $(\widehat E,\hat h)$ if there is a T-duality triple connecting them.
The notion of T-duality triples which we are going to call topological triples 
(Definition \ref{DefiTopTriples}) is central for this work, so
let us clarify what it means that a T-duality triple connects the 
pairs $(E,h)$ and $(\widehat E,\hat h):$

A T-duality triple is a commutative diagram
\begin{equation}\label{DiagTheFirstTDualityDiag}
\xymatrix{
&P\times_B\widehat E\ar[dl]\ar[dr]&&E\times_B\widehat P\ar[dl]\ar[dr]\ar[ll]_\cong^\kappa&\\
P\ar[dr]&&E\times_B\widehat{E}\ar[dr]\ar[dl]&&\widehat{P}\ar[dl]\\
&E\ar[dr]&&\widehat E\ar[dl]&\\
&&B,&&
}
\end{equation}
wherein $P\to E$ and $\widehat P\to \widehat E$ are principal bundles with
structure group $\PU(\HH)$, the projective unitary group of some infinite dimensional, separable Hilbert space $\HH$, such that both of these bundles
are trivialisable when restricted to the fibres of $E\to B$ or $\widehat E\to B$
respectively.
Moreover, the top-isomorphism $\kappa$ satisfies the following local condition:
Due to the triviality condition on the bundles $P,\widehat P$,  we can trivialise (\ref{DiagTheFirstTDualityDiag}) over each $u\in B$ such that
the isomorphism $\kappa$ induces a map $\kappa(u):\T^n\times\T^n\to \PU(\HH)$ which implements the isomorphism.
Now, $\PU(\HH)$ is an Eilenberg-McLane space of type $K(\Z,2)$ so
$\kappa(u)$ defines a class $[\kappa(u)]\in H^2(\T^n\times \T^n,\Z)$.
We force this class to satisfy $[\kappa(u)]\in \pi + {\rm im}({\rm pr}_1^*)+{\rm im}({\rm pr}_2^*)$, where ${\rm pr}_{1,2}:\T^n\times \T^n\to\T^n$ are the projections
and $\pi$ is the class of the tautological line bundle over $\T^n\times \T^n$ which
is $\pi=y_1\cup \hat y_1+\dots +y_n\cup\hat y_n$, for the generators
$y_1,\dots,y_n,\hat y_1,\dots,\hat y_n$ of $H^1(\T^n\times \T^n,\Z)$.

\begin{rem}\label{RemThefirstRemOnTriples}
In fact, this is not the definition of \cite{BS}. Firstly, they uses the more general notion of twists instead of the bundles $P$ and $\widehat P$, but the category 
of $\PU(\HH)$-principal bundles with homotopy classes of bundle isomorphisms 
as morphisms is a model of twists.
Secondly, they require the (a priori) more restrictive condition
that the class of the bundles $[P]\in H^3(E,\Z)$ (analogously for 
$[\widehat P]\in H^3(\widehat E,\Z)$)
lies in the second step $F^2H^2(E,\Z)$ of the filtration
$\{0\}\subset F^3H^3(E,\Z)\subset  F^2H^3(E,\Z)\subset  F^1H^3(E,\Z)\subset H^3(E,\Z)$ associated to the Leray-Serre spectral sequence. The requirement 
of triviality over the fibres of $E$ we stated above precisely means that the class lies in the first step $F^1H^3(E,\Z)$ of the filtration.
However, these two conditions on the classes are equivalent as we 
show in Lemma \ref{LemTopTriplesAreTDualityTriples}.
\end{rem}

A T-duality triple which we denote by 
$(\kappa, (P,E),(\widehat P,\widehat E))$ 
connects the two pairs $(E,h)$ and $(\widehat E, \hat h)$ if
we have an equality of the classes $[P]=h$ and 
$[\widehat P]=\hat h$.
\\

To a T-duality triple we can associate 
two $C^*$-algebras, namely the $C^*$-algebra of 
sections $\Gamma(E,F)$ and $\Gamma(\widehat E,\widehat F)$
of the associated bundles $F:=P\times_{\PU(\HH)}\K(\HH)$ and 
$\widehat F:=\widehat P\times_{\PU(\HH)}\K(\HH)$. 
It is the very 
aim of this work to understand how these two $C^*$-algebras are 
related to each other.
This issue turns our focus on the $C^*$-algebraic approach 
to T-duality \cite{MR, BHM2} which is based on the understanding of abelian $C^*$-dynamical systems.
\\
\\
We shortly summarise some $C^*$-algebraic background.

The duality theory of abelian $C^*$-dynamical systems
which  has been investigated for a quite long time \cite{Pe}
 is the foundation to understand T-duality by  $C^*$-algebraic means.
The dual of an abelian $C^*$dynamical system $(A,G,\alpha)$, 
i.e. $A$ a $C^*$-algebra with strongly continuous action 
$\alpha: G\to {\rm Aut}(A)$ of a locally compact, abelian
group $G$, 
is the
crossed product $C^*$-algebra $G\times_\alpha A$ 
equipped with the natural action $\hat \alpha$ of the dual group $\widehat G$,
i.e. $(G\times_\alpha A,\widehat G,\hat \alpha)$ becomes again 
a $C^*$-dynamical system (see \cite{Pe} or section \ref{SecTheSecOnCrossedProducts}).
A central result is the  Takai duality theorem (Theorem \ref{ThmTakaiduality})
which states in particular that the  bi-dual $C^*$-algebra
is stably isomorphic to the original one, i.e. 
they are Morita equivalent.
Thus, it is completely trivial to  understand the structure 
of the bi-dual and the difficult task is to understand the dual
$G\times_\alpha A$.

In the 80s and 90s big progress has been made to understand
the dual in case $A$ is a continuous trace algebra which we  assume 
from now on.
The basic structure theorem of Dixmier and Douady (see e.g. \cite{Di})
says that any separable, stable continuous trace algebra $A$
is isomorphic to $\Gamma_0(E,F)$ the $C^*$-algebra of sections vanishing at
infinity, where $E:={\rm spec}(A)$ is the 
spectrum of $A$ and $F\to E$ is a locally trivial bundle with each fibre
isomorphic to the compacts $\K(\HH)$. Their isomorphism 
classes are classified by $\check H^2(E,\underline{\U(1)})\cong H^3(E,\Z)$ 
(cp. section \ref{SecTheUnitaryAndProjectiveUnitaryGroup}),
and the class in $H^3(E,\Z)$ which determines the isomorphism type
of $A\cong \Gamma_0(E,F)$ is called the Dixmier-Douady invariant of $A$.

A first result \cite{Pe2,RW} for an understanding the crossed product
$G\times_\alpha A$
was that if $G$ is compactly generated 
and the induced action of $G$  on the spectrum $E$ of $A$ is trivial,
then  the crossed product $G\times_\alpha A$ is  isomorphic
to the balanced tensor product
$C_0(\widehat E)\otimes_{C_0(B)}A=: p^* A$, wherein
$p:\widehat E\to B$
is  a $\widehat G$-principal bundle consisting of the spaces
$\widehat E:={\rm spec}(G\times_\alpha A)$ and  $B:=E$.

The more general situation
wherein the action of $G$ does not fix the spectrum $E$ of
$A$ but has constant isotropy group $N$ for each 
$\pi\in E$ is concerned in \cite{RR}.  One of 
the statements therein is the following. 
Assume that $E$ with the induced action of $G/N$ 
is a principal fibre bundle $E\to B:= E/(G/N)$ and that
the restricted action $\alpha|_N$ of $N$ on $A$ is locally
unitary, then there is a pull back diagram of 
principal fibre bundles 
$$
\xymatrix{
&\hspace{1cm}E\times_B\widehat E\cong\hspace{-1cm}&{\rm spec}(N\times_{\alpha|_N} A)\ar[rd]^{\hat p}\ar[ld]_p&\\
\hspace{1.1cm} E:=\hspace{-1.1cm}&{\rm spec}(A)\ar[dr] 
&\hspace{1.8cm}\widehat E:=\hspace{-1.8cm}& {\rm spec}(G\times_\alpha A)\ar[dl]\\
&&B,&
}
$$
wherein the down-right arrows have fibre $G/N$
and the down-left arrows have fibre 
$\widehat N\cong \widehat G/N^\perp$.
($N^\perp$ is the annihilator of $N$ which is the 
set of characters of $G$ whose restriction to $N$ is 
identically 1.)
Moreover, $p^* A$ is isomorphic to $N\times_{\alpha|_N} A$
and Morita equivalent to $\hat p^*(G\times_\alpha A)$.
Thus, we  have the following  schematic situation of 
$C^*$-algebras over their spectra
\begin{equation}\label{DiagSimilarToTDualDiag}
\xymatrix{
&p^* A\ar@{.>}[dr]&&N\times_{\alpha|_N} A\ar@{..>}[dl]\ar[ll]_\cong
&\\
A\ar@{..>}[dr]&&E\times_B\widehat E\ar[rd]^{\hat p}\ar[ld]_p&
&G\times_\alpha A\ar@{.>}[dl]\\
&E\ar[dr] &&\widehat E\ar[dl]&\\
&&B&&
}
\end{equation}
\vspace{-5.0cm}
\begin{eqnarray*}
\hspace{3.5cm}
&\wr\ _{^{^{\rm  Morita}}}&\\
& \hat p^*(G\times_\alpha A)&
\end{eqnarray*}
\vspace{3.0cm}

\noindent
which obviously is similar to 
diagram (\ref{DiagTheFirstTDualityDiag}).
The question is whether or not it is possible that 
both $A$ and $G\times_\alpha A$ are separable, stable
continuous trace algebras. This question has been answered  in \cite[Thm. 6]{ER}.
In particular, this is true if $\alpha|_N$ is point-wise unitary and the 
action of $G/N$ on $E\times_B E\cong {\rm spec}(N\times_{\alpha|_N} A)$ is 
proper, e.g. $G/N$ is compact.

This finishes our summary.
\\

The approach to T-duality of \cite{MR} 
considers the following set-up.
Let $E$ be a locally compact space 
(with certain finiteness assumptions) 
with an action of the torus $\T^n$ such that
$E\to B:=E/\T^n$ becomes a principal torus bundle,
and let $h\in H^3(E,\Z)$. They
concern these data as a stable continuous trace algebra
$\Gamma_0(E,F)$ that has Dixmier-Douady invariant
$h$. 
Under which circumstances is it possible to lift 
the $\T^n$-action from the spectrum $E$ 
 to an action of $\R^n=:G$ on
$\Gamma_0(E,F)$? If so, is it further 
possible to obtain an action whose restriction 
to $N:=\Z^n$ is point-wise unitary for all $\pi\in E$?
These questions are answered in
\cite[Thm 3.1]{MR}. The general answer is 
no, but if we make further restrictions to the
class $h$ one achieves a positive answer.
Namely, a lift $\alpha$  to an $\R^n$-action exists 
if and only if $h\in F^1H^3(E,\Z)$ the first 
step of the filtration 
$$
0\subset F^3H^3(E,\Z)\subset F^2H^3(E,\Z)\subset
F^1H^3(E,\Z)\subset H^3(E,\Z)
$$ 
associated to the Leray-Serre spectral sequence,
and there exists a point-wise  unitary action
if even $h\in F^2H^3(E,\Z)$.
Consequently, in the latter  case there exists a
T-dual in the sense that there exists a 
dual space-time $\widehat E$ in diagram (\ref{DiagSimilarToTDualDiag})
which is a $\widehat{\Z^n}\cong \T^n$-principal fibre
bundle over $B$ and is the spectrum of the
stable continuous trace algebra 
$\R^n\times_\alpha\Gamma(E,F)\cong\Gamma_0(\widehat E,\widehat F)$.
\\

In this work we are going  to show that the two approaches 
to T-duality are essentially equivalent, i.e. there is no difference 
between (equivalence classes of) T-duality triples and
(equivalence classes of)  abelian $C^*$-dynamical systems 
which are obtained as described above. 
To fulfil this  task we must develop a technique 
which enables us to compare such different objects.
The methods of \cite{BS} and \cite{MR} are not 
applicable for such a manoeuvre as they are
too less explicit: The explicit description of the local structure of these 
two different kinds of objects  can be used to describe  
a transformation as desired.

Our method is general enough that we do not
have to restrict ourselves to the case of the
groups $\R^n$ and $\T^n\cong\R^n/\Z^n$. -
We develop a theory for all second countable, locally compact, abelian 
groups $G$ with lattice $N\subset G$, i.e. $N$ is a discrete, cocompact  subgroup. 
\\
\\
We give an overview  of this work.

A basic technique we use throughout the whole of 
this work is to lift all local, projective unitary families of functions 
(e.g. the transition functions of a $\PU(\HH)$-principal bundle
$P\to E$) to unitary Borel- or $L^\infty$-functions in direction
of the fibres $G/N$ of $E\to B$ and to think of them
as  new unitary (multiplication) operators in $\U(L^2(G/N)\otimes\HH)$.
We call this procedure Borel lifting technique.
The technical condition we must assume is that 
the base $B$ is a paracompact Hausdorff space which
is locally contractible.
We call spaces with these properties base spaces, and 
the whole theory we develop is a theory over base spaces.

In sections \ref{subsecTheCategoryOfPairs} and 
\ref{SubsecThelocalStrOfPairsAndTheisClassifying} we introduce 
the notion of  pairs. A pair is a $G/N$-principal fibre bundle 
$E\to B$ over a base space $B$ and a $\PU(\HH)$-principal fibre bundle 
$P\to E$ which is trivialisable over the fibres of $E$.
We explain their local structure and give a quick 
result on their classification.

In section \ref{SecPairsAndTwistedCech} we introduce 
a twisted version of $\check{\rm C}$ech cohomology on the
base $B$, where the twist is given by the bundle $E\to B$.
The Borel lifting technique mentioned above defines 
a map from (equivalence classes of) pairs into the
second twisted $\check{\rm C}$ech cohomology
(similar to the ordinary definition of the second  
$\check{\rm C}$ech class of a $\PU(\HH)$-bundle).

In sections \ref{SubsecDynamicalTriples} we extend the same procedure 
to dynamical triples $(\rho,E,P)$ which are pairs $(P,E)$ 
equipped  with a  lift $\rho$ of the $G/N$
action on $E$ to a $G$-action on $P$.
Due to  this  action,  group cohomological 
expressions arise in the description of their local structure   and lead finally  
to a map from (equivalence classes of) dynamical triples to
the second  cohomology  of a double  complex
which has  one  (twisted) 
$\check{\rm C}$ech cohomological direction
and a second  group cohomological direction.
A two cocycle (or two cohomology class) of this complex has three entries:
a  pure $\check{\rm C}$ech part due to the 
transition functions of the pair, a group cohomological 
part and a third mixed term. 
In section \ref{SecTheDualityTheoryOfDynTrip}
we focus our attention to those triples for which the
group cohomological entry vanishes.

Section \ref{SecDualPairsAndTriples}
just contains the definition of dual pairs and dual
dynamical triples which are nothing more than
pairs and dynamical triples, but as underlying groups we
take  the dual group $\widehat G$ of $G$ with dual lattice 
$N^\perp$ the annihilator of $N$.

Then in section \ref{SecTopologicalTriples} we 
introduce topological triples.
Our definition is a straight forward  
generalisation to arbitrary locally compact,
abelian groups $G$ with lattice $N$ from the notion of a 
T-duality triple ($G=\R^n, N=\Z^n$).

In section \ref{SecTheDualityTheoryOfDynTrip} we state our
first main result. We first
single out a specific subclass of dynamical triples which we
call dualisable. They are those dynamical triples for which
the group cohomological entry of the associated 
two cohomology class (of the double complex) vanishes.
We construct an explicit map 
$[(\rho,E,P)]\mapsto[(\hat \rho,\widehat E,\widehat P)]$ from the set
of equivalence classes of dualisable dynamical triples
to the set of equivalence classes of dualisable dual
dynamical triples, and show that it is a bijection
whose inverse is given by the dual map 
(defined by replacing everything by its dual counterpart).
In this sense dualisable dynamical triples and dual dualisable 
triples are in duality.

Section \ref{SecTheRelationToTopologicalTDuality} contains
two important statements. The first is that we have  a
map $\tau(B)$ from the set of equivalence classes of dualisable dynamical triples to
the set of equivalence classes of topological triples 
(everything understood over a base space $B$).
This map is defined by the duality theorem
of section \ref{SecTheDualityTheoryOfDynTrip}, i.e.
the topological triple we define consists of the 
pairs of two dynamical triples in duality.
Then we try to define a map $\delta(B)$ in the opposite 
direction which generally fails 
as an obstruction occurs.
However, on  the subset of those topological triples which
have a vanishing obstruction we then find a construction 
of a whole family of dualisable dynamical triples
which is
associated to a topological triple.
This construction, when restricted to the image of the first
map, can be turned into  an honest map, i.e. we have a preferred choice of an
element of the family, and this map is inverse to $\tau(B)$.

Section \ref{SecTheCaseOfRAndZ} is devoted to the special case of
the group $G=\R^n$ with lattice $N=\Z^n$. 
In this situation the construction of the  map
$\delta(B)$  simplifies drastically,
the group wherein the obstruction lives vanishes.
As a result, we can associate to each topological triple  
a dynamical triple which is unique (up to equivalence) 
because the family of dynamical triples degenerates to a family of
one single element only.
As it turns out the two maps  $\tau(B)$ and $\delta(B)$
are  bijections and inverse to each other.

Moreover, the four maps mentioned above are natural
in the base, so they define natural transformations of
functors. 
Thus, the main result of section \ref{SecTheCaseOfRAndZ}
can be restated as follows.
In the case of $G=\R^n$ and $N=\Z^n$, we have a completely 
explicit construction of  equivalences of functors
$$
{\rm Top}\cong{\rm Dyn}^\dag\cong \widehat{\rm Dyn}^\dag.
$$
Therein, ${\rm Dyn}^\dag$ ($\widehat{\rm Dyn}^\dag$) is the functor sending 
a base space to the set of equivalence classes of
dualisable (dual) dynamical triples over it. ${\rm Top}$
is the functor which sends a base space to the set
of equivalence classes of topological (T-duality) triples over
it.

In section \ref{SecTheStructureOfTheAssociatedCStarDynSys}
we point out that the theory developed so far is 
connected to the theory of $C^*$-dynamical systems precisely
as one expects. Namely,
the  $C^*$-dynamical systems 
$(G\times_{\alpha^\rho} \Gamma(E,F),\widehat G,\widehat{\alpha^\rho})$ 
and $(\Gamma(\widehat  E,\widehat F),\widehat G,\alpha^{\hat\rho})$
are isomorphic, wherein
$(\hat \rho, \widehat E,\widehat P)$ is 
the dual of $(\rho, E, P)$ and
$(\rho,E,P)\mapsto (\Gamma(E,F),G,\alpha^\rho)$
is the functor which sends a dualisable (dual) dynamical triple 
to its corresponding $C^*$-dynamical system, i.e.
$F$ is the 
associated $\K(\HH)$-bundle 
to $P$ and 
$\alpha^\rho$ is the by $\rho$ induced action on the
$C^*$-algebra of sections $\Gamma(E,F)$.

In an appendix we put some technical lemmata and 
notation we are going to use.

\newpage
\section{Pairs and Triples}

\subsection{The Category of Pairs}
\label{subsecTheCategoryOfPairs}

\begin{defi}\label{DefiOfBaseSpace}
A {\bf base space} $B$ is a topological space which is
Hausdorff, paracompact and locally contractible.
\end{defi}
The category of base spaces consists
of bases spaces as objects and  continuous maps between them
as morphisms.
A typical class of base spaces are CW-complexes \cite[Thm. 1.3.2, Thm. 1.3.5]{FP}.

By $G$ we will always denote a second countable,
Hausdorff, locally compact abelian group and 
by $N$  a discrete, cocompact subgroup, i.e.
the quotient $G/N$ is compact.

Let $\HH$ be an infinite dimensional, separable Hilbert space.
Let $E\to B$ be a $G/N$-principal fibre bundle and
$P\to E$ be a $\PU(\HH)$-principal fibre bundle.

\begin{defi}\label{DefiPairOverB}
We call the data\footnote{differing from \cite{BS}} 
$(P,E)$ a {\bf pair} over $B$ with underlying Hilbert
space $\HH$ if 
\begin{itemize}
 \item[$i)$] $B$ is a base space,
\item[$ii)$]
the restriction of the bundle $P\to E$
to the fibres of $E\to B$ is trivialisable.
\end{itemize}
\end{defi}

\begin{rem}
We do not require local compactness for $B$,
because we want to develop
a theory that includes non locally compact spaces 
such as  classifying spaces of groups (s. next section).
Therefore $E$ need not to be locally compact 
and  thus need not equal  the spectrum of
any  (continuos trace) $C^*$-algebra such
as $\Gamma(E,F)$ the $C^*$-algebra of bounded sections
(or $\Gamma_0(E,F)$ the $C^*$-algebra of sections
vanishing at infinity [A section vanishing at infinity
vanishes already identically on the set of points which
don't have a compact neighbourhood.])
 of the associated 
$C^*$-bundle $F:=P\times_{\PU(\HH)}\K(\HH)$. 
\end{rem}

A {\bf morphism} $(\varphi,\vartheta,\theta)$ over $B$ from a   pair
$P{\to} E\to B$ with underlying Hilbert space $\HH$ 
to a pair  $P'{\to} E'\to B$ with underlying Hilbert space $\HH'$  is
a commutative diagram of bundle isomorphisms
\begin{equation}\label{DiagramMorphismOfPairs}
\xymatrix{
P\ar[r]^\vartheta        \ar[d]  & \varphi^*P'\ar[d]&\\
E\ar[r]^{\theta} \ar[d]  & E'\ar[d]& \\
{B}\ar[r]^{=}             &B,      &
}
\end{equation}
wherein $\varphi\in\PU(\HH,\HH'):=\U(\HH,\HH')/\U(1)$ is the
class of a unitary isomorphism $\HH\to \HH'$.
$\varphi^*P'$ is the $\PU(\HH)$-bundle with
total space  $\varphi^*P'=P'$, but with 
$\PU(\HH)$-action that is induced by 
$\varphi^*:\PU(\HH)\to \PU(\HH')$, i.e. 
$x'\cdot U:=x'\cdot (\varphi^*)^{-1} U, x'\in P', U\in \PU(\HH).$
Pairs over a base space $B$
and their morphisms form a category;
composition of morphisms $(\varphi,\vartheta,\theta)$ and $(\varphi',\vartheta',\theta')$ 
is just component-wise composition 
$(\varphi'\circ \varphi,\vartheta'\circ \vartheta,\theta'\circ \theta).$
This category  is even a groupoid, i.e. every  morphism is 
an isomorphism.

This notion of morphism
is well-behaved under stabilisation
in the following sense.
Let $(P,E)$ be a pair over $B$ with underlying Hilbert 
space $\HH$.
Let $\HH_0$ be any separable Hilbert space, not neccessarily 
 infinite dimensional. Then
$\PU(\HH)$ is isomorphic to the subgroup $\Eins\otimes\PU(\HH)$
of $\PU(\HH_0\otimes\HH)$ and thus 
$\PU(\HH)$ acts  on $\PU(\HH_0\otimes,\HH)$ by 
left multiplication: $(U,V)\mapsto (\Eins_{\HH_0}\otimes U) V$.
Then the associated (stabilised) bundle
\begin{eqnarray}
P_{\HH_0} :=\PU(\HH_0)\otimes P:= P\times_{\PU(\HH)}\PU(\HH_0\otimes \HH)
\label{EqStablilisedPUBundleP}
\end{eqnarray}
is a $\PU(\HH_0\otimes\HH)$-principal bundle and
$(P_{\HH_0},E)$ is a pair over $B$ with underlying Hilbert space
$\HH_0\otimes\HH$. 
We call  
two pairs 
$(P,E)$ and $(P',E')$ with underlying Hilbert 
spaces $\HH$ and $\HH'$  
{\bf stably isomorphic} if there exists a Hilbert
space $\HH_0$ such that
 the pairs $(P_{\HH_0},E)$ and $(P'_{\HH_0},E')$
are isomorphic.

\begin{prop}
Two pairs over $B$ are stably isomorphic if and only if
they are isomorphic.
\end{prop}
\begin{pf}
It is clear that isomorphic pairs are stably isomorphic.
To prove the converse it suffices to show that
$P$ and $\varphi^*P_{\HH_0}$ are isomorphic over $E$, for
an isomorphism $\varphi:\HH\cong\HH_0\otimes\HH$.
To do so, we show that $P$ and $\varphi^*P_{\HH_0}$ 
define the same $\check{\rm C}$ech
class, hence the classification theorem of $\PU(\HH)$-bundles 
(Theorem \ref{ThmClassifacationOfPUBundlesAndAutos})
implies that the two bundles are isomorphic.

In fact, if $\zeta^1_{ji}:V_{ji}\to \PU(\HH)$
are transition functions for $P$ (here $\{V_i\}$ is a covering of $E$) , then
$\zeta^2_{ji}:=\varphi^*(\Eins\otimes\zeta_{ji})$
are transition functions for $\varphi^*P_{\HH_0}$.
If we refine the covering such that the transition functions $\zeta^1_{ji}$ 
lift to unitary-valued functions $\overline\zeta^1_{ji}$ (Lemma \ref{LemLiftForLocTrivFB}), 
then these lifts define also lifts $\overline\zeta^2_{ji}$ for the other family of transition  
functions. Thus, on threefold intersections $V_{kji}$ the cocycle 
identities of the two families 
are perturbed by the same $\check{\rm C}$ech 2-cocycle 
$$
c_{kji}:=\overline\zeta^n_{ji}{\overline\zeta^n_{ki}}^{-1}\overline\zeta^n_{kj}:V_{kji}\to\U(1)\cdot\Eins,
\quad n=1,2.
$$
Hence, their classes agree.
\end{pf}

\begin{rem}
In case of an  additional structure  such as a group action of $G$ on $P$  
stable isomorphism and  isomorphism are different notions when we force 
them to preserve the extra structure.
This will be important in section \ref{SubsecDynamicalTriples}. 
\end{rem}

Let us denote by ${\rm Par}$ the set valued contravariant functor
that sends a base space $B$ to the set of stable isomorphism classes 
of pairs over $B$, i.e.
$$
{\rm Par}(B):= \{\textrm{ pairs over } B\}\big/_{\rm stable\ isomorphism}\ ,
$$
and if $f:B'\to B$ is a continuous map between base spaces, 
then pullback defines a map
$f^*:{\rm Par}(B)\to {\rm Par}(B')$.

There is a  subcategory of pairs over $B$ consisting of pairs 
with a fixed $G/N$-bundle $E\to B$ and morphisms of the form
$(\varphi,\vartheta,{\rm id}_E)$, and we call pairs 
$(P,E)$ and $(P',E)$ {\bf stably isomorphic over} $E$
if there is a  isomorphism of this special form between
$(P_{\HH_1},E)$ and $(P'_{\HH_1},E)$, for 
a Hilbert space $\HH_1$.
We define
$$
{\rm Par}(E,B):= \{\textrm{ pairs over } B \textrm{ with fixed } E\}{\big/_{
\textrm{stable isomorphism over } E}}\ ,
$$
and  for a bundle morphism
$$
\xymatrix{
E'\ar[d]\ar[r]^\theta & E\ar[d]\\
B'\ar[r] &B
}
$$ 
we define by pullback a map
$\theta^* :{\rm Par}(E,B)\to{\rm Par}(E',B')$,
so ${\rm Par}(\ .\ ,..)$ becomes a contravariant functor
from the category of $G/N$-principal fibre bundles over
base spaces to sets.
The bundle automorphisms ${\rm Aut}_B(E)$
of $E$ over ${\rm id}_B$
act on ${\rm Par}(E,B)$ by pullback, 
and we have a decomposition 
$
{\rm Par}(B)\cong \coprod_{ [E]} ({\rm Par}(E,B)/{\rm Aut}_B(E)),
$
wherein the 
disjoint union runs over all isomorphism classes of
$G/N$-bundles $E\to B$.

\begin{rem}
For each fixed $E\to B$ the set ${\rm Par}(E,B)$ has a 
natural group structure.
If $[(P,E)]$ and $[(P',E)]$ are two classes of pairs,
then we let
$[(P,E)]+[(P',E)]:=
[(P\otimes P',E)]$,
wherein $P\otimes P'$
is the $\PU(\HH\otimes\HH')$-bundle
which is associated to the 
$\PU(\HH)\times\PU(\HH')$-bundle
$P\times_E P'$,
$$
P\otimes P':= 
(P\times_E P')\times_{\PU(\HH)\times\PU(\HH')} \PU(\HH\otimes\HH').
$$
The unit element is given by the class
of a trivial bundle and the inverse of
$[(P,E)]$  
is  given by the class $[( P^\#,E)]$ 
of the complex conjugate bundle $P^\#$
which is as space  the bundle $P$ but has the
action $(x,U)\mapsto x\cdot U^\#$.
$U^\#$ is here the complex conjugate (not the adjoint) 
of $U\in\PU(\HH)$  (which may be defined by identifying $\HH=l^2(\N)$ and 
taking the complex conjugate  matrix of $u=(u_{ij})_{i,j\in \N}$, for
$U={\rm Ad}(u)$ ).

In this way we just mimic the group structure 
of $\check H^2(E,\underline{\U(1)})$, i.e.
the classification map\footnote{See Theorem \ref{ThmClassifacationOfPUBundlesAndAutos}.} ${\rm Par}(E,B)\to \check H^2(E,\underline{\U(1)})$
is turned into a group homomorphism.
\end{rem}

An {\bf automorphism} of a pair is
a morphism from a pair onto itself.
The group of automorphisms of a pair
is denoted by ${\rm Aut}(P,E)$.
It becomes a topological group when equipped 
with the initial topology of the forgetful 
map
${\rm Aut}(P,E)\to \PU(\HH,\HH')\times{\rm Map}(P,P)$ which 
sends a morphism $(\varphi,\vartheta,\theta)$ to $(\varphi,\vartheta)$, wherein 
$\U(\HH,\HH')$
has the strong topology, i.e.
the topology of point-wise convergence.
${\rm Map}(P,P)$ has the compact open topology.
\label{PageofAut}

Since $G/N$ is a commutative group
mappings  of the form
$\theta_z:E\ni e\mapsto e\cdot z\in E$,
$z\in G/N$, are bundle morphisms.
They give rise to a subgroup
${\rm Aut_1}(P,E)$ which
consist of all  morphisms $({\rm id_\HH},\vartheta,\theta_z)$.
Let ${\rm Aut_0}(P,E)$ denote 
the subgroup consisting of morphisms 
$({\rm id_\HH},\vartheta,{\rm id}_E)$, then we find a
short exact sequence of topological groups
$$
\Eins\to {\rm Aut_0}(P,E)\to{\rm Aut_1}(P,E)\to G/N\to 0.
$$

\subsection{The Local Structure of Pairs and their Classifying Space}
\label{SubsecThelocalStrOfPairsAndTheisClassifying}
By $\A$ we denote  the group of automorphisms of the trivial pair over a point.
In particular, $a=(\varphi,\vartheta,\theta)\in \A$ makes
$$
\xymatrix{
G/N\times\PU(\HH)\ar[r]^{\vartheta}
\ar[d] & G/N\times\PU(\HH)\ar[d]\\
G/N\ar[r]^{\theta=\_+z} \ar[d]         &G/N\ar[d] \\
{*}\ar[r]                    &{*}
}
$$
commute, for some $z\in G/N.$ It is immediate that  there is an
isomorphism
\begin{eqnarray*}
\A&\cong& G/N \ltimes {\rm Map}(G/N,\PU(\HH))\rtimes\PU(\HH)
\end{eqnarray*}
of topological groups, wherein
$ G/N \ltimes {\rm Map}(G/N,\PU(\HH))\rtimes\PU(\HH)$ 
is the semi-direct product with multiplication  
$(y,\eta,u)\cdot(z,\zeta,v):= (y+z,(z\cdot \eta\cdot v)\  \zeta,uv)$.
The action on the continuous functions ${\rm Map}(G/N,\PU(\HH))$ is 
$$
(z\cdot \eta\cdot v)(x):=v^{-1} \eta(x+z) v.
$$
Let $(P,E)$ be any pair over $B$ with
underlying Hilbert space $\HH$.
$B$ is a base space, and so we can choose 
a covering $\{U_i | i\in I\} $ of $B$ \
of  open sets such that for each $U_i$ there is a 
commutative diagram
\begin{equation}\label{DiagTheLocalStructureOfAPair}
\xymatrix{
U_i\times G/N\times \PU(\HH)
\ar[r]_{h_i}^\cong\ar[d]&q^{-1}(p^{-1}(U_i))\ar[r]  \ar[d]  & P\ar[d]^q&\\
U_i\times G/N\ar[r]_{k_i}^\cong \ar[d]&p^{-1}(U_i)\ar[r] \ar[d]       & E\ar[d]^p& \\
\ U_i\ \ar[r]^{=} 			&\ U_i\ \ar[r]^{\subset}             &B,       &
}
\end{equation}
with bundle isomorphisms $k_i,h_i$.
We refer
to such a covering as an {\bf atlas} $U_\bullet=\{(U_i,k_i,h_i)\ |\ i\in I\}$ consisting
of the {\bf charts} $(U_i,k_i,h_i)$.
The transition from one chart to another is 
described by a set of transition functions. For 
a pair this consists of two families of continuous  functions 
$g_{ij}:U_{ij}\to G/N$ and $\zeta_{ij}:U_{ij}\to {\rm Map}( G/N, \PU(\HH))$
which appear in\label{PageOfTransiFunc}
\begin{eqnarray*}
h_j^{-1}\circ h_i:U_{ji}\times G/N\times \PU(\HH)&\to& U_{ij}\times G/N\times \PU(\HH).\\
(u,z,U)&\mapsto& (u,g_{ji}(u)+z,\zeta_{ji}(u)(z)U)
\end{eqnarray*}
It follows that on threefold
intersections $U_{ijk}$
the relations $g_{ki}(u)=g_{kj}(u)+g_{ji}(u)$ and
\begin{eqnarray}\label{EqTransitForPuBundle}
\zeta_{ki}(u)(z)=\zeta_{kj}(u)(g_{ji}(u)+z)\ \zeta_{ji}(u)(z)\in \PU(\HH)
\end{eqnarray}
are valid; equivalently, the family of functions
$$
a_{ij}:=g_{ij}\times\zeta_{ij}:U_{ij}\to G/N\ltimes {\rm Map}(G/N,\PU(\HH))=:\A_1\subset\A
$$
satisfies the  $\check{\rm C}$ech 1-cocycle condition
$a_{ij}(u)a_{jk}(u)=a_{ik}(u)$.\label{PageOfA1EA1BA1}
Now, let $\E\A_1\to\B\A_1$ be the universal $\A_1$-principal
fibre bundle.
We call the associated pair
$$
\xymatrix{
&P_{\rm univ}\ar[d]&
\!\!\!\!\!\!\!\!\!\!\!\!\!\!:=\E\A_1\times_{\A_1}(G/N\times\PU(\HH))\\
&E_{\rm univ}\ar[d]&
\!\!\!\!\!\!\!\!\!\!\!\!\!\!\!\!\!\!\!\!\!\!\!
\!\!\!\!\!\!\!\!\!\!\!\!\!\!\!\!\!\!\!\!:=\E\A_1\times_{\A_1} G/N \\
&B_{\rm univ}&
\!\!\!\!\!\!\!\!\!\!\!\!\!\!\!\!\!\!\!\!\!\!\!
\!\!\!\!\!\!\!\!\!\!\!\!\!\!\!\!\!\!\!\!\!\!\!
\!\!\!\!\!\!\!\!\!\!\!\!\!\!\!\!\!\!\!\!:=\B\A_1
}
$$
the {\bf universal pair}. Indeed we
can choose a CW-model for $\B\A_1$ such
that the universal pair is a pair in the sense of 
Definition \ref{DefiPairOverB}.
Its name is due to  the following  universal
property.

\begin{prop}\label{PropTheClassifyingSpaceOfPairs}
The space $B_{\rm univ}$ classifies pairs over 
(pointed) CW-complexes, i.e. if $B$ is a (pointed)
CW-complex, then
$$
[B,B_{\rm univ}]\cong {\rm Par}(B),
$$
wherein the left hand side is the set of (pointed) 
homotopy classes of maps $B\to B_{\rm univ}$.
\end{prop}
\begin{pf}
We already observed that 
the transition functions of a pair define 
a $\check{\rm C}$ech class
$[a_{..}]=[g_{..}\times\zeta_{..}]\in\check H^1(B,\underline{\A_1}).$
This class is independent of the chosen atlas.
Now, let $(\varphi,\vartheta,\theta):(P,E)\to(P',E')$ be an isomorphism of pairs,
then the transition functions  $a'_{..}$ for $(P',E')$ can be turned
into a $\check{\rm C}$ech $\A_1$-cocycle induced by
$\varphi^*:\PU(\HH)\cong\PU(\HH')$, and this class  depends 
on the isomorphism class of the pair only. 
Conversely, the associated pairs of isomorphic $\A_1$-principal bundles
are isomorphic, and each isomorphism class arises. 

Thus, pairs and $\A_1$-principal bundles over $B$ have the same isomorphism
classes, but for a (pointed) CW-complex $B$ the latter isomorphism class
is  given by homotopy classes of maps to $\B\A_1$.
\end{pf}

\subsection{Pairs and Twisted $\check {\rm C}$ech Cohomology}
\label{SecPairsAndTwistedCech}

Let $\M$ and $\G$ be  abelian (pre-)sheaves on a space $B$ and
assume that $\M$ is a right $\G$ module. 
We are going to twist the $\check{\rm C}$ech coboundary operator of $\M$
by a $\check{\rm C}$ech $\G$ 1-cocycle.
Fix an open covering $U_{\bullet}=\{U_{i}|i\in I\}$ of $B$ and
let $g\in \check Z^1(U_{\bullet},\G)$.
Then for $\varphi\in \check C^{n-1}(U_{\bullet},\M),n=1,2,\dots,$
we define
$$
\delta_{g}\varphi :=\delta\varphi +g^\star\varphi,
$$
wherein $\delta$ is the ordinary $\check{\rm C}$ech coboundary operator\footnote{See \ref{SecCechandGroupCoho}} and
$$
(g^\star\varphi)_{i_0\dots i_{n}}:=
(-1)^{n-1}\big(
\varphi_{i_0\dots i_{n-1}}-\varphi_{i_0\dots i_{n-1}}\cdot g_{i_{n-1}i_n}\big)|_{U_{i_0 \dots i_n}}.
$$
The choice of the sign $(-1)^{n-1}$ is such that the last term of $\delta \varphi$ and the
first term of $g^\star\varphi$ cancel.
We obtain a sequence 
$$
\xymatrix{
0\ar[r]&\check C^{0}(U_{\bullet},\M)\ar[r]^{\delta_{g}}&
\check C^{1}(U_{\bullet},\M)\ar[r]^{\delta_{g}}&
\check C^{2}(U_{\bullet},\M)\ar[r]^{\quad\delta_{g}}&\cdots
}
$$

\begin{lem}
$(\check C^{\bullet}(U_{\bullet},\M),\delta_{g})$ is a cochain complex.
\end{lem}
\begin{pf}
We have to show that the square of $\delta_{g}$ vanishes.
Let $\varphi\in \check C^{n-1}(U_{\bullet},\M)$, then
\begin{eqnarray*}
\delta_{g}\delta_{g}\varphi&=&
\delta\delta\varphi+\delta(g^\star\varphi)+g^\star\delta\varphi
+g^\star(g^\star\varphi)\\
&=&
\delta(g^\star\varphi)+g^\star\delta\varphi
+g^\star(g^\star\varphi)\\
\end{eqnarray*}
and therefore
\begin{eqnarray*}
(\delta_{g}\delta_{g}\varphi)_{i_0 \dots i_{n+1}}&=&
\sum_{k=0}^{n+1}(-1)^k(g^\star\varphi)_{i_0 \dots \hat i_k \dots i_{n+1}}\\
&&+(-1)^n\big(\sum_{k=0}^n(-1)^k\varphi_{i_0 \dots \hat i_k\dots i_{n}}
-\sum_{k=0}^n(-1)^k\varphi_{i_0\dots\hat i_k\dots i_{n}}\cdot g_{i_n i_{n+1}}\big)\\
&&+(-1)^n (g^\star\varphi)_{i_0\dots i_n}-(-1)^n
(g^\star\varphi)_{i_0\dots i_n}\cdot g_{i_n i_{n+1}}\\
&=&\sum_{k=0}^{n-1}(-1)^k
(-1)^{n-1}(\varphi_{i_0\dots \hat i_k\dots i_n}-
\varphi_{i_0\dots\hat i_k\dots i_n}\cdot g_{i_n i_{n+1}})\\
&&+(-1)^n (-1)^{n-1}(\varphi_{i_0\dots i_{n-1}}-\varphi_{i_0\dots i_{n-1}}\cdot g_{i_{n-1} i_{n+1}})\\
&&+(-1)^{n+1} (-1)^{n-1}(\varphi_{i_0\dots i_{n-1}}-\varphi_{i_0\dots i_{n-1}}\cdot g_{i_{n-1} i_{n}})\\
&&+(-1)^n\big(\sum_{k=0}^n(-1)^k\varphi_{i_0 \dots \hat i_k\dots i_{n}}
-\sum_{k=0}^n(-1)^k\varphi_{i_0\dots\hat i_k\dots i_{n}}\cdot g_{i_n i_{n+1}}\big)\\
&&+(-1)^n(-1)^{n-1}\big(\varphi_{i_0\dots i_{n-1}}-
\varphi_{i_0\dots i_{n-1}}\cdot g_{i_{n-1}i_n}\big)\\
&&-(-1)^{n}(-1)^{n-1}\big(\varphi_{i_0\dots i_{n-1}}-
\varphi_{i_0\dots i_{n-1}}\cdot g_{i_{n-1}i_n}\cdot g_{i_n i_{n+1}}\big)\\
&=&0.
\end{eqnarray*}
The last equality holds due to the cocycle relation $g_{ij}+g_{jk}=g_{ik}$.
\end{pf}

By the last lemma we have well-defined  cohomology groups 
$\check H^n(U_\bullet,\M,g)$ for  $n=0,1,2,\dots$ and any open cover
$U_\bullet=\{U_i | i\in I\}$ of $B$.
And as in the untwisted case 
we  define the twisted cohomology groups $\check H^n(B,\M,g)$
of $B$ by passing  to the limit 
$$
\check H^n(B,\M,g):={\rm lim}_{V_\bullet}\check H^n(V_\bullet,\M,g|),
$$
wherein the limit runs over all refinements $V_\bullet$ of $U_\bullet$.
To be precise,
consider a refinement  $V_\bullet=\{V_k | k\in K\}$ 
of $U_\bullet=\{U_i | i\in I\}$ with refinement
map $\iota:K\to I$ , i.e.
$V_k\subset U_{\iota(k)}$ 
then we define 
$(\iota^*g)_{kl}:= g_{\iota(k)\iota(l)}|_{V_{kl}}$ 
and similarly $(\iota^*\varphi)_{k_0\dots k_n}:= 
\varphi_{\iota(k_0)\dots\iota(k_n)}|_{V_{k_0\dots k_n}}$ 
to obtain  a  cochain map 
$\iota^*:(\check C^n(U_\bullet,\M),\delta_g)\to (\check C^n(V_\bullet,\M),\delta_{\iota^*g})$.
This construction defines a functor from the category of 
coverings with refinement maps as morphisms 
to 
the category of cochain complexes (and after
taking homology to the category of graded abelian groups).
The category of coverings is filtered in the following sense:
\label{PageForLimit}

$(i)$ Any two coverings have a common refinement, 
i.e. for any two objects $(U_\bullet,I)$ and $(V_\bullet,K)$ 
there is a third object $(W_\bullet,L)$ with morphisms
$(U_\bullet,I)\to (W_\bullet,L)$ and $(V_\bullet,K)\to(W_\bullet,L)$.

$(ii)$ Any two refinement maps become equal finally, i.e.
for any to morphisms 
$\iota:(U_\bullet,I)\to (V_\bullet,K)$ and
$\kappa:(U_\bullet,I)\to (V_\bullet,K)$ there exists an object
$(W_\bullet,L)$ and morphisms
$\iota':(V_\bullet,K)\to (W_\bullet,L)$ and
$\kappa':(V_\bullet,K)\to (W_\bullet,L)$ such that
$\kappa'\circ\kappa=\iota'\circ\iota$.

Due to $(i)$ and $(ii)$ the  limit 
${\rm lim}_{V_\bullet}\check H^n(V_\bullet,\M,g|)$
is independent of the choice of the first covering 
and independent of  the 
refinement maps.
\\

So far,  in our construction  we referred explicitly to 
a choice of a cocycle $ g\in \check Z^1(U_\bullet,\G)$,
but up to   isomorphism $\check H^n(B,\M,g)$
depends only on the class of $g$.
In fact, let $V_\bullet=\{V_k | k\in K\}$ and
 $U_\bullet=\{U_i | i\in I\}$ be open coverings of $B$, and let 
 $g'\in \check Z^1(V_\bullet,\G)$  and $g\in \check Z^1(U_\bullet,\G)$ 
represent the same element in $\check H^1(B,\G)$. 
 If $W_\bullet=\{W_m | m\in M\}$ is a common refinement
with  refinement maps $\iota:M\to I$ and $\kappa: M\to K$
then there are  $r_{m}\in \G(W_m)$ such that
$g'_{\kappa(m)\kappa(n)}|_{W_{mn}} = r_{m}|_{W_{mn}}\ +
g_{\iota(m)\iota(n)}|_{W_{mn}}\ -r_{n}|_{W_{mn}}.$
They give rise to the following diagram  of cochain complexes
$$
\xymatrix{
(\check C^{\bullet}(U_{\bullet},\M),\delta_{g})\ar[d]^{\iota^*}&
(\check C^{\bullet}(V_{\bullet},\M),\delta_{g'})\ar[d]^{\kappa^*}\\
(\check C^{\bullet}(W_{\bullet},\M),\delta_{\iota^*g})&
(\check C^{\bullet}(W_{\bullet},\M),\delta_{\kappa^*g'})\ar[l]_\cong^{r^\#}&
}
$$
wherein $r^\#$ is defined by 
$(r^\#\varphi)_{k_0\dots k_{n}}:=\varphi_{k_0\dots k_{n}}\cdot
r_{k_n}|_{V_{k_0\dots k_n}}$, for $\varphi\in\check C^{n}(V_{\bullet},\M)$.
One easily derives $\delta_{\iota^*g}(r^\#\varphi)=r^\#\delta_{\kappa^*g}\varphi$, for
$\varphi\in \check C^{\bullet}(V_{\bullet},\M)$, i.e. $r^\#$ is a cochain map and even
an isomorphism. Thus the corresponding cohomology groups 
are isomorphic and this isomorphism passes to the  limit.
Therefore the twisted $\check{\rm C}$ech groups $\check H^n(B,\M,[g])$
are well defined for the class $[g]\in \check H^1(B,\G)$ up to the
considered isomorphism.
\\\\
We now consider  the relation between pairs and twisted
$\check{\rm C}$ech cohomology.
One should note at this point that, just  as in the untwisted  
case, the first $\check{\rm C}$ech ``group`` $\check H^1(B,\M,g)$
is a well-defined set even in case the sheaf $\M$ is just a 
sheaf of groups and not necessarily  abelian. In that case
the additive (commutative) relation of being cohomologous 
$\zeta_{ji}\sim\zeta_{ji}+(\delta_{g}\eta)_{ji}$ is replaced by 
the multiplicative  relation 
$\zeta_{ji}\sim \eta_j|_{U_{ij}}\cdot g_{ji}\ \zeta_{ji}\ \eta_i^{-1}|_{U_{ij}}$, for
$\zeta\in \check Z^1(U_\bullet,\M,g),
\eta\in \check C^0(U_\bullet,\M)$.
The next proposition is then just a reformulation of 
what we already observed in Proposition \ref{PropTheClassifyingSpaceOfPairs}.

\begin{prop}
Let $B$ be a base space.
Let $\M$ be the sheaf of topological groups on $B$
defined by $\M(U):=C(U,{\rm Map}(G/N,\PU(\HH)))$, and let
$\G:=\underline{G/N}$, i.e. $\underline{G/N}(U):=C(U,G/N)$, 
for $U\subset B$. $\underline{G/N}$ acts on $\M$
in the obvious way, i.e. by translation in the arguments. 
Then the first twisted  cohomology
classifies pairs over $B$, i.e. we have a bijection
$$
\qquad \check H^1(B,\M,g)\cong
\textrm{Par}(E,B),
$$
if the class $[g]\in \check H^1(B,\underline{G/N})$ is the class for  the bundle $E$.
\end{prop}
\begin{pf}
Let $g_{..},\zeta_{..}$ be the transition
functions of a pair over $B$ for an atlas
$U_\bullet$. So $g\in \check Z^1(U_\bullet, \underline{G/N})$,
and the crucial point is to observe that
equation (\ref{EqTransitForPuBundle}) is 
equivalent to $\delta_{g}\zeta=\Eins$
and therefore each pair defines an element in
$\check H^1(B,\M,g)$. In fact, this is well defined,
because if $g'_{..},\zeta'_{..}$ are transition functions for
another atlas then (after choosing a common refinement)
the two classes match under the isomorphism $r^\#$, i.e.
$r^\#\zeta'$ is cohomologous to $ \zeta$.
Similarly, an isomorphism of pairs leads to cohomologous
cocycles. 
Conversely, any $\zeta\in \check Z^1(U_\bullet,\M,g)$
defines an associated pair, and if $\zeta'\in \check Z^1(V_\bullet,\M,g)$
defines  the same class as $\zeta_{..}$ then the two associated pairs are isomorphic. 
Since each class arises in such a way the assertion is proven.
\end{pf}

Let $g_{..},\zeta_{..}$ be the transition
functions of a pair over $B$.
Since $G/N$ is compact and $B$ paracompact we
can apply  Lemma \ref{LemMapIsBaseOfPFB}
and Lemma \ref{LemLiftForLocTrivFB} for 
the family of transition functions 
$\zeta_{ij}:U_{ij}
\to {\rm Map}(G/N,\PU(\HH)$.
I.e. we can find a refined atlas
$\{ V_{k} | k \in K:= I\times B \}$, $V_{k}\subset U_i$ if
$k=(i,x)$,
such that on its twofold intersections
$V_{kl}$
the restricted transition functions 
lift to continuous functions 
$\overline\zeta_{kl}:V_{kl}\to {\rm Bor}(G/N,\U(\HH))$. 
These lifts are unique up to 
continuous functions  $V_{kl}\to{\rm Bor}(G/N,\U(1))$.
Let  us denote by $g_{lk}$  the restriction
$g_{ji}|_{V_{lk}}$ in case $l=(j,y),k=(i,x)\in I\times B$.
On threefold intersections 
the function $V_{klm}\ni u\mapsto \overline\zeta_{kl}(u)(g_{lm}(u)+\_)$
won't be continuos as a function to 
${\rm Bor}(G/N,\U(\HH))$ in general, but it will as 
a function to\footnote{$L^\infty(G/N,\U(\HH))$ has the weak topology. See equation (\ref{EqTheDefiOfLinfty}) in section \ref{SecTheUnitaryAndProjectiveUnitaryGroup}.} 
$L^\infty(G/N,\U(\HH))$.
So 
equation (\ref{EqTransitForPuBundle})
implies 
that there are continuous $\psi_{mlk}:V_{mlk}\to L^\infty(G/N,\U(1))$
such that \label{PageOfPsi}
\begin{eqnarray}
\overline\zeta_{ml}(u)(g_{lk}(u)+z)\
\overline\zeta_{lk}(u)(z) 
&=&
\overline\zeta_{mk}(u)(z)\
\psi_{mlk}(u)(z)\cdot\Eins\nonumber\\
\Leftrightarrow \delta_{g}\overline\zeta&=&\psi\ ,
\label{EqUnitCocycleConditonForTransFDuex}
\end{eqnarray}
and it follows that $\delta_g\psi=1$.
The functions $\psi_{klm}$ therefore
define a twisted $\check{\rm C}$ech 2-cocycle
$\psi_{...}\in\check Z^2(V_\bullet,\underline{L^\infty(G/N,\U(1))},g)$.

\begin{prop}\label{PropWouldLikeSecondTwistedCechClassifiesPairs}
The construction of $\psi_{...}$ defines 
a homomorphism of groups
$$
\textrm{Par}(E,B)
\to\check H^2(B,\underline{L^\infty(G/N,\U(1))},g)
$$
if $[g]\in\check H^1(B,\underline{G/N})$ is the class 
of the bundle $E\to B$.
\end{prop}
\begin{pf}
We must check that the class $[\psi_{...}]$ is independent of 
all choices.
In fact, if $\overline\zeta_{ji}$ and $\overline\zeta_{ji}'$ are 
different choices of lifts of $\zeta_{ji}$, they differ
by a scalar function $s_{ji}=\overline\zeta_{ji}^{-1} \overline\zeta_{ji}'$
and $\psi \delta_g s$ is the cocycle obtained from $\overline\zeta_{ji}'$,
so the class $[\psi]$ is not effected.
It is also easy that the class does not 
change under the choice of the atlas
or by considering a pair stably isomorphic to the first one.
\end{pf}

We do not achieve the statement that the above homomorphism
is injective or surjective, so we are far from classifying pairs 
by this map.

\subsection{Dynamical Triples and their Local Structure}\label{SubsecDynamicalTriples}
Let $P\to E\to B$ be a pair. 
The quotient map $G \ni g\mapsto gN\in G/N$  
induces a $G$ action on $E$.

\begin{defi}
A {\bf decker} is just a continuous 
action $\rho: P\times G\to P$ that
lifts the induced $G$-action on $E$
such that 
$\rho(\ . \ ,g):P\to P$ is a bundle automorphism,
for all $g\in G$.
\end{defi}

The existence of deckers can be a very 
restrictive condition  on the bundle $P\to E$. 
(See e.g. Prop. \ref{PropDeckerForG/NExIff} below.)
In fact, they  need not exist and need not to be unique
in general, but the play a central rôle in what follows, 
therefore we introduced an extra name.

In  context of $C^*$-dynamical systems, 
i.e. $C^*$-algebras with (strongly continuous) group actions,
concretely, in context of the equivariant Brauer group 
several notions 
of equivalence of actions occur \cite{CKRW}. 
In particular, the notions of  isomorphic actions, 
stably isomorphic actions and exterior equivalent actions
are combined to the notion of stably outer conjugate actions.
We  slightly  modify  these notion for our purposes.
However, we postpone the definition until we made 
ourselves familiar with the local structure of 
dynamical triples.

\begin{defi}
A {\bf dynamical  triple} $(\rho,P,E)$ over
$B$ is a pair $(P,E)$ over $B$ together  with
a decker $\rho: P\times G\to P$. 
\end{defi}

Let $(\rho,P,E)$ be a dynamical tripel over $B$.
 
\begin{prop}\label{PropTheLocalStrOfDeckers}
\begin{itemize}
\item[$i)$]
If we define $\rho^\tau(g):=\rho(\ .\ ,g):P\to P$, we obtain 
a diagram of  topolgical groups 
\begin{equation}\label{EqLift}
\xymatrix{
0\ar[r]&N\ar[r]\ar[d]^{\rho^\tau|_N}&G\ar[d]^{\rho^\tau}\ar[r]&G/N\ar[r]\ar[d]^=&0\\
\Eins\ar[r]& {\rm Aut_0}(P,E)\ar[r]&{\rm Aut_1}(P,E)\ar[r]& G/N\ar[r]& 0.
}
\end{equation}
Conversely, if $B$ is locally compact then a commutative 
diagram (\ref{EqLift}) defines a decker.
\item[$ii)$]
Locally, i.e. after choosing  charts 
$U_i$,
a decker defines a family  of   continuous 
cocycles
$\mu_i:U_i\to Z^1_{\rm cont}(G,{\rm Map}(G/N,\PU(\HH)))$
such that on twofold intersections $U_{ij}\ni u$ the transition functions
of the pair and the cocycles are related by
\begin{eqnarray}\label{EqLocStucOfDeckers}
\mu_i(u)(g,z)=\zeta_{ji}(u)(z+gN)^{-1}\ \mu_j(u)(g,g_{ji}(u)+z)\ \zeta_{ji}(u)(z).
\end{eqnarray}
Conversely, any  family of cocycles $\{\mu_i\}_{i\in I}$ that fulfils 
eq. (\ref{EqLocStucOfDeckers}) determines a unique decker.
\end{itemize}
\end{prop}

\begin{pf}
\begin{itemize}
\item[$i)$] The origin of the diagram is obvious.
For the converse,
it is sufficient to prove the result locally 
because the action of $G$ on $F$ preserves 
charts. Explicitly, over a chart $(U_i,k_i,h_i)$ the action
of $g\in G$ is
\begin{eqnarray*}
h_i^*(\phi(g)):U_i\times G/N\times\PU(\HH)&\to&U_i\times G/N\times\PU(\HH).\\
(u,z,U)&\mapsto&(u,gN+z,\mu''_i(g)(u,z)U)
\end{eqnarray*}
If $B$ is locally compact the
exponential law  (Lemma \ref{LemExpLawForMap})
ensures that all functions
$\mu'_i:(u,g,z)\mapsto \mu''_i(g)(u,z)$ are jointly
continuous.
\item[$ii)$] Locally,
Lemma \ref{LemExpLawForMap} ensures
that the transposed functions 
$\mu_i:U_i\to{\rm Map}(G\times G/N,\PU(\HH))$
made out of $\mu'_i$ are well defined. The cocycle condition
$$
\mu_i(u)(g+h,z)=\mu_i(u)(g,z+hN)\ \mu_i(u)(h,z)
$$
and the validity of (\ref{EqLocStucOfDeckers}) are immediate 
as well as the converse statement.
\end{itemize}
\end{pf}

It should be mentioned at this point
that equation (\ref{EqLocStucOfDeckers}) 
(and its unitary version we consider later)
is  quite powerful  as turns out. 
A first application is given in the next 
proposition.
It is a complete answer to the existence of deckers 
in the case of $N=0$, i.e. $G=G/N$.
The result is well-known, but we state a proof
using (\ref{EqLocStucOfDeckers}) 
for the convenience of the reader.

\begin{prop}\label{PropDeckerForG/NExIff}
Assume  $N=0$ and let $P\stackrel{q}{\to} E\stackrel{p}{\to}B$ be a pair.
Then a decker exists if and only if $P\cong p^*P'$ for 
a $\PU(\HH)$-bundle  $P'\to B$.
\end{prop}

\begin{pf}
If $P\cong p^*P'$ for some 
$\PU(\HH)$-bundle $P'\to B$ we obtain a decker 
by acting on the first entry of the fibered product 
$p^*P'=E\times_B P'.$

Conversely, if a decker is given
then we define a family 
$\zeta'_{ji}:U_{ji}\to\PU(\HH)$ by
$\zeta'_{ji}(u):=\mu_j(u)(g_{ji}(u),0)^{-1}\ \zeta_{ji}(u)(0)$,
for $u\in U_{ji}\subset B$.
This is well defined since $G=G/N.$
\\
Claim 1 : $\{\zeta'_{ji}\}_{j,i\in I}$ are transition functions
for a $\PU(\HH)$-bundle $P'\to B.$\\
Proof : Let $u\in U_i\cap U_j\cap U_k.$ Then
\begin{eqnarray*}
&&\zeta'_{kj}(u)\ \zeta'_{ji}(u)\\
&=&
\mu_k(u)(g_{kj}(u),0)^{-1}\ \zeta_{kj}(u)(0)\ \mu_j(u)(g_{ji}(u),0)^{-1}\ \zeta_{ji}(u)(0)\\
&\stackrel{(\ref{EqLocStucOfDeckers})}{=}&
\underbrace{\mu_k(u)(g_{kj}(u),0)^{-1}\  \mu_k(u)(g_{ji}(u),g_{kj}(u)+0)^{-1}}\
\underbrace{ \zeta_{kj}(u)(0+g_{ji}(u))\ \zeta_{ji}(u)(0)}\\
&&\qquad\qquad\qquad\qquad
\stackrel{\rm cocy. cond.}{=}\mu_k(u)(g_{ji}(u)+g_{kj}(u),0)^{-1}\qquad
\stackrel{(\ref{EqTransitForPuBundle})}{=}\zeta_{ki}(u)(0)\\
&=&\zeta'_{ki}(u).
\end{eqnarray*} 
\\
Claim 2 : $P\cong p^*P'$.\\
Proof : The bundle $q':p^*P'\to E$ has
transition functions
$$
p^*\zeta'_{ji} : p^{-1}(U_{ij})\cong U_{ij}\times G/N \ni (u,z)\mapsto \zeta'_{ji}(u).
$$
We define an isomorphism $f:P\to p^*P'$ locally
$f_i:= f|_{q^{-1}(p^{-1}(U_i))}$ by
\begin{eqnarray*}
P\supset q^{-1}(p^{-1}(U_i)  &\stackrel{f_i}
{\longrightarrow}& q'^{-1}(p^{-1}(U_i))\subset p^*P'\\
\downarrow\cong\qquad&&\qquad\downarrow\cong\\
 U_i\times G/N\times \U(\HH)
&\longrightarrow& U_i\times G/N\times \U(\HH)\\
(u,z,U)&\longmapsto&(u,z,\mu_i(u)(z,0)^{-1}U).
\end{eqnarray*}
This is in fact a well defined global isomorphism
since 
form eq. (\ref{EqLocStucOfDeckers})
it follows for $G=G/N$ and $u\in U_i\cap U_j$
\begin{eqnarray*}
\mu_i(u)(z,0)^{-1}&=&\zeta_{ji}(u)(0)^{-1}\ \mu_j(u)(z,g_{ji}(u))^{-1}\ \zeta_{ji}(u)(z)\\
&=&\zeta_{ji}(u)(0)^{-1}\ \mu_j(g_{ji}(u))(u,0)\ \mu_j(u)(z+g_{ji}(u),0)^{-1}\ \zeta_{ji}(u)(z)\\
&=&\zeta'_{ji}(u)^{-1}\ \mu_j(u)(z+g_{ji}(u),0)^{-1}\ \zeta_{ji}(u)(z).
\end{eqnarray*}
Thus the local definition of $f$ is independent of the chosen
chart.
\end{pf}

We now 
introduce the notions of equivalence we mentioned earlier.

\begin{defi}\label{DefiOfExteriorEquivalence}
Two deckers $\rho,\rho':G\times P\to P$ on a pair $(P,E)$ 
are called {\bf exterior equivalent}
if the  continuous  function $c:P\times G\to P$ which is defined  
by  $\rho(\ .\ ,-g)\circ \rho'(\ .\ ,g)=c(\ .\ ,g):P\to P$, for $g\in G$,
is locally of the following form:
There exists an atlas $U_\bullet$ such that for each chart $(U_i,k_i,h_i)$ 
there is a continuous map  $c_i:U_i\times G\times G/N\to \U(\HH)$ 
which satisfies  the three conditions
\footnote{Here and in what follows 
we will always use the notation 
$\zeta(c), \zeta\in\PU(\HH),c\in\U(\HH)$, for the 
action of $\PU(\HH)$ on $\U(\HH)$ by conjugation.}
\begin{enumerate}
\item[\rm{(E0)}] $(h_i^{-1}\circ c(\_,g)\circ h_i) (u,z,U)= (u,z,{\rm Ad}(c_i(u,g,z)) U)
\in U_i\times G/N\times \PU(\HH)$,

\item[\rm{(E1)}] $c_j(u,g,g_{ji}(u)+z)= \zeta_{ji}(u)(z)\big(c_i(u,g,z)\big)\in \U(\HH)$ and
\item[{\rm (E2)}]
$
c_i(u,h+g,z)= 
 \mu_i(u)(g,z)^{-1}\big(c_i(u,h,z+gN)\big)\ c_i(u,g,z)\in \U(\HH),
$
\end{enumerate}
for the transition functions $g_{ji},\zeta_{ji}$  of the pair
and  the cocycles $\mu_i$ of the decker $\rho$ as above.
\end{defi}

It is  clear that,
if one has given a family of continuous unitary functions  
$\{c_i\}$ which satisfies  {\rm (E0), (E1)} and {\rm (E2)} for
the cocycles $\{\mu_i\}$ of a decker $\rho$,  
the family of cocycles
$$
\mu'_i(u)(g,z):= \mu_i(u)(g,z)\ {\rm Ad}(c_i(u,g,z))
$$ defines 
an exterior equivalent decker $\rho'$ to $\rho$.

Let $c^\tau:G\to{\rm Aut}_0(P,E)$ 
be defined by \label{PageOfCTau}
$c^\tau(g):=c(\ .\ ,g)$.
It satisfies the cocycle condition 
$$
c^\tau(g+h) = c^\tau(g)\cdot\rho^\tau(h)\ c^\tau(h),
$$
where $\cdot$ is the right action of 
${\rm Aut}_1(P,E)$ on ${\rm Aut}_0(P,E)$ given
by conjugation. 
This right action lifts to the sections 
$\Gamma(E,P\times_{\PU(\HH)}\U(\HH))$ 
of the to $P$ associated $\U(\HH)$-bundle
in the following diagram, i.e. 
the vertical map is ${\rm Aut}_1(P,E)$-equivariant,
$$
\xymatrix{
&&&\Gamma(E,P\times_{\PU(\HH)}\U(\HH))\ar[d]^{\circlearrowleft {\rm Aut}_1(P,E)}\\
G\ar[rr]^{c^\tau}\ar@{-->}@/^0.6cm/[rrru]^{\overline{c^\tau}}&&{\rm Aut}_0(P,E)\ar[r]^{\!\!\!\!\!\!\!\!\!\!\!\!\!\!\!\!\!\!\!\!\!\cong} 
&\Gamma(E,P\times_{\PU(\HH)}\PU(\HH)),
}
$$
Therein the associated bundles are both obtained by the
the conjugate action of $\PU(\HH)$ on the respective groups.
Now, conditions (E0) - (E2) imply that
we can lift 
$c^\tau$ to a unitary cocycle $\overline{c^\tau}:G\to\Gamma(E,P\times_{\PU(\HH)}\U(\HH))$,
i.e.
$$
\overline{c^\tau}(g+h) = \overline{c^\tau}(g)\cdot\rho^\tau(h)\ \overline{c^\tau}(h)
$$
Similar to Proposition \ref{PropTheLocalStrOfDeckers}
this global statement 
is an equivalent formulation of
exterior equivalence if the base $B$ is locally compact.
\begin{prop}
Let $B$ be locally compact, and let
$\rho$ and $ {\rho'}$ be deckers on a pair 
$(P,E)$.
These two deckers are exterior equivalent
if and only if $c^\tau$ as defined above lifts to a unitary cocycle $\overline{c^\tau}$.
\end{prop}
We do not give a detailed proof of this fact, it is 
again just an application of the exponential law
for locally compact spaces.
\\

The next statement gives an important example 
of exterior equivalent deckers.

\begin{exa}\label{ExaOfExteriorEquivalence}
Let $\rho$ be a decker on an arbitrary  
pair $(P,E)$, and let 
$v:P\to P$ be a bundle automorphism.
Then the conjugate decker $\rho^\nu$
is exterior equivalent to $\rho$ if
the class  $[v]\in\check H^1(E,\underline{\U(1)})$
of the automorphism vanishes.
\end{exa}
\begin{pf}
Let us denote by $\mu_i$  the cocycles of the decker $\rho$
which satisfy (\ref{EqLocStucOfDeckers}) on a chosen
atlas $\{U_i\}_{i\in I}$.
Because the class of $v$ vanishes, 
we can assume without restriction that it is locally
implemented by  unitary functions $v_i:U_i\to {\rm Map}(G/N,\U(\HH))$ 
such that
$\zeta_{ji}(u)(z)(v_i(u)(z))=v_j(u)(g_{ji}(u)+z)$.

Locally, $\rho(\_,-g)\circ\rho^\nu(\_,g)$ is given by
\begin{eqnarray*}
&&\mu_i(u)(-g,z+gN){\rm Ad}(v_i(u)(z+gN))\mu_i(u)(g,z){\rm Ad}(v_i(u)(z)^{-1})\\
&=&\mu_i(u)(g,z)^{-1}{\rm Ad}(v_i(u)(z+gN))\mu_i(u)(g,z){\rm Ad}(v_i(u)(z)^{-1})\\
&=&{\rm Ad}(c_i(u,g,z)),
\end{eqnarray*}
for 
$c_i(u,g,z):=\mu_i(u)(g,z)^{-1}\big(v_i(u)(z+gN)\big)\ v_i(u)(z)^{-1}\in\U(\HH)$.
So condition (E0) is satisfied.
We check that the conditions (E1) and (E2) also holds. 
In fact,
\begin{eqnarray*}
&&\zeta_{ji}(u)(z)\big(c_i(u,g,z)\big)\\
&=&\zeta_{ji}(u)(z)\Big(
\mu_i(u)(g+h,z)^{-1})\big(v_i(u)(z+gN+hN)\big)\ v_i(u)(z)^{-1}\Big)\\
&\stackrel{(\ref{EqLocStucOfDeckers})}{=}&
\mu_j(u)(g,g_{ji}(u)+z)^{-1}\Big(
\zeta_{ji}(u)(z+gN)\big(
v_i(u)(z+gN)\big)\Big) \zeta_{ji}(u)(z)\big(v_i(u)(z)\big)^{-1}\\
&\stackrel{[\nu]=1}{=}&c_j(u)(g,g_{ji}(u)+z)
\end{eqnarray*}
which proves (E1), and 
\begin{eqnarray*}
&&c_i(u,g+h,z)\\
&=&\mu_i(u)(g+h,z)^{-1}\big(v_i(u)(z+gN+hN)\big)\ v_i(u)(z)^{-1}\\
&=&
\mu_i(u)(g,z)^{-1}\Big(\mu_i(u)(h,z+gN)^{-1}\big(v_i(u)(z+gN+hN)\big)\Big)\ v_i(u)(z)^{-1}\\
&=&
\mu_i(u)(g,z)^{-1}\Big(\mu_i(u)(h,z+gN)^{-1}\big(v_i(u)(z+gN+hN)\big)\ v_i(u)(z+gN)\Big)\\
&&
\mu_i(u)(g,z)^{-1}\big(v_i(u)(z+gN)\big)\ v_i(u)(z)^{-1}
\end{eqnarray*}
which proves (E2).
\end{pf}

Between two dynamical triples we introduce 
a notion of equivalence based on
exterior equivalence.

Two dynamical triples $(\rho,P,E)$ and $(\rho',P',E')$ are 
{\bf  isomorphic } if there is a morphism $(\varphi,\vartheta,\theta)$ of the underlying 
pairs such that $\rho=\vartheta^*\rho'$. 
The triples  are {\bf  outer conjugate } if there is a morphism 
$(\varphi,\vartheta,\theta)$ of the underlying 
pairs such that $\rho$ and $\vartheta^*\rho'$ are exterior equivalent on $(P,E)$.
Furthermore, we call the triples {\bf stably isomorphic} (respectively {\bf stably outer conjugate})
if the triples $(\Eins\otimes\rho,\PU(\HH)\otimes P,E)$
and $(\Eins\otimes\rho',\PU(\HH)\otimes P',E')$
are isomorphic (resp.  outer conjugate).

We can arrange these notions in a diagram of implications. 
$$
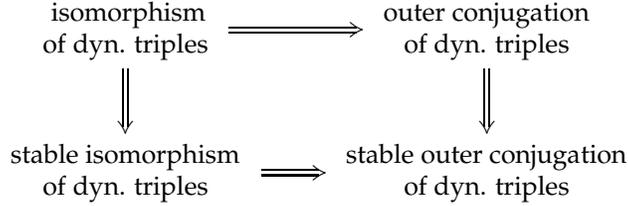
\xymatrix{
{\begin{array}{c} 
\textrm{ isomorphism} \\ 
\textrm{of dyn. triples} 
\end{array}}
\ar@{=>}[r]\ar@{=>}[d]&
{\begin{array}{c} 
\textrm{outer conjugation} \\ 
\textrm{of dyn. triples} 
\end{array}}
\ar@{=>}[d]\\
{\begin{array}{c} 
\textrm{stable isomorphism} \\ 
\textrm{of dyn. triples} 
\end{array}}
\ar@{=>}[r]&
{\begin{array}{c} 
\textrm{stable outer conjugation} \\ 
\textrm{of dyn. triples} 
\end{array}}
}
$$
The following  example shall illustrate
an   important feature of the notion of stably outer conjugation.
\begin{exa}
Let $(\rho,P,E)$ be a dynamical triple, and\label{PageTheExampleOfStableEquivalence}
let $\lambda_{_{G}}:G\to\U(L^2(G))$ be the left regular
representation of $G$. Then the two triples
$(\rho,P,E)$ and $(({\rm Ad}\circ \lambda_{_{G}})\otimes \rho,
\PU(L^2(G))\otimes P,E)$ are stably outer conjugate.
\end{exa}
\begin{pf}
The triple $(\rho,P,E)$ and its stabilisation  
$(\Eins\otimes\rho,\PU(L^2(G))\otimes P,E)$
are stably isomorphic and  the triples $(\Eins\otimes\rho,\PU(L^2(G))\otimes P,E)$
and $({\rm Ad}\circ\lambda_G\otimes\rho,\PU(L^2(G))\otimes P,E)$
are exterior equivalent by $c_i(u,g,z):=\lambda_{_{G}}(g)\otimes\Eins_\HH$.
\end{pf}

By ${\rm Dyn}$ \label{PageOfDyn} we denote the set valued functor
that sends a base space $B$ to the set of 
equivalence classes of stably outer conjugate 
dynamical triples over it, i.e.
$$
{\rm Dyn}(B):=\{\textrm{dynamical triples over } B\}\big/_\textrm{stable outer conj.}.
$$
In the same manner as we did for the functor ${\rm Par}$
we can fix a bundle $E\to B$ and define
$$
{\rm Dyn}(E,B):=\{\textrm{dynamical triples over } B \textrm{ with fixed }E\}
\big/_{\textrm{stable outer conj. over }{\rm id}_E}.
$$
Isomorphic bundles $E, E'$ lead to isomorphic
sets ${\rm Dyn}(E,B)\cong{\rm Dyn}(E',B)$, and
the bundle automorphisms ${\rm Aut}_B(E)$ act
on ${\rm Dyn}(E,B)$ by pullback. This yields a decomposition 
${\rm Dyn}(B)\cong \coprod_{[E]} ({\rm Dyn}(E,B)/{\rm Aut}_B(E))$.
\\


Our next goal is to find the link between dynamical triples
and the cohomology theory we introduce now. 
Let $\M^n$ be the abelian sheaf on $B$ defined 
by $\M^n(U):=C(U,{\rm Bor}(G^{\times n},L^\infty(G/N,\U(1)) ))$,
for $n=0,1,2,\dots$. 
Let $U_\bullet=\{ U_i| i\in I\}$ be an open cover of $B$ and
let $g\in \check Z^1(U_\bullet,\underline{G/N})$.
Note that $\M^n$ is a right $\underline{G/N}$-module, for all $n=0,1,2,\dots,$
by shifting the $G/N$-variable.
We consider the the double complex 
$$
\xymatrix{
\vdots&\vdots&\vdots& \\
\check C^2(U_\bullet,\M^0)\ar[r]^{d_*}\ar[u]^{\delta_g}&
\check C^2(U_\bullet,\M^1)\ar[r]^{d_*}\ar[u]^{\delta_g}&
\check C^2(U_\bullet,\M^2)\ar[r]^{\quad d_*}\ar[u]^{\delta_g}&\cdots\\
\check C^1(U_\bullet,\M^0)\ar[r]^{d_*}\ar[u]^{\delta_g}&
\check C^1(U_\bullet,\M^1)\ar[r]^{d_*}\ar[u]^{\delta_g}&
\check C^1(U_\bullet,\M^2)\ar[r]^{\quad d_*}\ar[u]^{\delta_g}&\cdots\\
\check C^0(U_\bullet,\M^0)\ar[r]^{d_*}\ar[u]^{\delta_g}&
\check C^0(U_\bullet,\M^1)\ar[r]^{d_*}\ar[u]^{\delta_g}&
\check C^0(U_\bullet,\M^2)\ar[r]^{\quad d_*}\ar[u]^{\delta_g}&\cdots,
}
$$
wherein the horizontal arrows $d_*$ are 
induced  by the boundary operator 
$d:{\rm Bor}(G^{\times n}, L^\infty(G/N,\U(1)) )\to
{\rm Bor}(G^{\times n+1},L^\infty(G/N,\U(1))) $ of group cohomology\footnote{See \ref{SecCechandGroupCoho}}
by acting point-wise on functions
and the vertical arrows $\delta_g$ are the twisted $\check{\rm C}$ech
coboundary operators.
In fact, all of the squares are commutative, hence we
obtain a resulting total complex
 $\big( C_{\rm tot}^\bullet(U_\bullet,\M^\bullet),\partial_g\big)$, i.e.
$C_{\rm tot}^p(U_\bullet,\M^\bullet):=\bigoplus_{p=k+l} \check C^k(U_\bullet,\M^l)$
and $\partial_g|_{C_{\rm tot}^p(U_\bullet,\M^\bullet)}:=\delta_g-(-1)^{p} d_*$ define
a cochain complex. The choice of the sign in $\partial_g$ will be convenient.
By $H^\bullet_{\rm tot}(U_\bullet,\M^\bullet,g)$ we denote the corresponding  
cohomology groups, and by passing to the  limit over all refinements 
of the open covering $U_\bullet$ we obtain 
$$
H^\bullet_{\rm tot}(B,\M^\bullet,g):={\rm lim}_{V_\bullet} 
H^\bullet_{\rm tot}(V_\bullet,\M^\bullet,g).
$$
We call this group the {\bf total cohomology} of $B$ with twist $g$.
In the same manner as explained on page \pageref{PageForLimit}
the limit does not depend on the covering $U_\bullet$ and is
independent of the choice of the refinement maps.

It is also similar to the discussion on twisted $\check{\rm C}$ech cohomology 
that the total cohomology groups are well defined objects for
the class $[g]$ of a cocycle $g$ up to an isomorphism.
\\

The connexion of dynamical triples and total cohomology has its origin 
in the local structure of triples as  we explain now.
Let $(\rho,P,E)$ be a dynamical triple. 
Let $U_\bullet=\{U_i|i\in I\}$ be an atlas for the underlying pair
with transition functions $g_{ij},\zeta_{ij}$
and continuous cocycles $\mu_i$ as in Proposition \ref{PropTheLocalStrOfDeckers}.
One should  realise at this point that when we
suppress the non-commutativity of $\PU(\HH)$ for a moment 
the  equations $\delta_g\zeta_{..}=\Eins, d(\mu_{i}(u))=\Eins$ and 
equation (\ref{EqLocStucOfDeckers}) are equivalent 
to $\partial_g(\zeta_{..},\mu_{.})=\Eins$. We lift 
the transition functions and the cocycles to Borel functions.
This will  define a 2-cocycle for the total cohomology of $B$;
in detail: 
Without restriction (see equation (\ref{EqUnitCocycleConditonForTransFDuex}))
we can assume that the 
atlas is chosen such that the transition functions can be lifted continuously to
Borel valued functions $\overline\zeta_{ij}$; they define a twisted $\check{\rm C}$ech
2-cocycle $\delta_g\overline\zeta=:\psi_{...}\in \check C^2(U_\bullet,\M^0)$.
Further we can assume (if necessary we refine the atlas once more) that
all charts are contractible.
Therefore we can apply Corollary \ref{CorExistenceOfBorelCocycle}
to each of the $\mu_i$
and obtain continuous functions 
$\overline\mu_i:U_i\to{\rm Bor}(G\times G/N,\U(\HH))$
lifting  $\mu_i$, for all $i\in I$.
By Lemma \ref{LemTheInclusionIsConti} we can also
pass to
$\overline\mu_i:U_i\to{\rm Bor}(G, L^\infty(G/N,\U(\HH)))$
which is important to as $L^\infty(G/N,\U(\HH))$ is a continuous 
$G/N$-module.
From the cocycle identity $d(\mu_i(u))=\Eins$ we see that
$d(\overline\mu_i(u))=:\omega_i(u)$ defines 
a  group cohomology 2-cocycle
$\omega_.\in\check C^0(U_\bullet,\M^2)$.
On twofold intersections we can define\label{EqThisIsTheDefiOfPhiOhYear}
$\phi_{ji}(u)(g,z):=\overline\mu_i(u)(g,z) \overline\zeta_{ji}(u)(z)^{-1}
\overline\mu_j(u)(g,g_{ji}(u)+z)^{-1}\overline\zeta_{ji}(u)(z+gN)$ which is due to
equation (\ref{EqLocStucOfDeckers}) $\U(1)$-valued,
i.e. $\phi_{..}\in \check C^1(U_\bullet,\M^1)$.
Now,  the three families
of functions $\psi_{...},\phi_{..}$ and $\omega_.$
satisfy the algebraic relations
\begin{eqnarray*}
\delta_g\psi&=&1\\
\delta_g\phi&=& d_*\psi\\
d_*\phi&=&\delta_g\omega\\
d_*\omega&=&1
\end{eqnarray*}
which is equivalent to  
\begin{equation}
\partial_g(\psi_{...},\phi_{..},\omega_.)=0\in C^3_{\rm tot}(U_\bullet,\M^\bullet),
\end{equation}
i.e. $(\psi_{...},\phi_{..},\omega_.)$ is a total 2-cocycle. 
Of course, one can verify this by direct computation,
but indeed it is implicitly clear, because, informally\footnote{
i.e. up to the non-comutativity of $\U(\HH)$},
we have defined
$(\psi_{...},\phi_{..},\omega_.):=\partial_g
(\overline\zeta_{..},\overline\mu_.)
\in C^2_{\rm tot}(U_\bullet,\M^\bullet)$.

\begin{prop}\label{PropWouldLikeToClassifyTriples}
The assignment 
$(\rho,P,E)\mapsto (\psi_{...},\phi_{..},\omega_.)$
constructed above defines a homomorphism of groups
$$
{\rm Dyn}(E,B)\to H^2_{\rm tot}(B,\M^\bullet,g_{..}).
$$
\end{prop}
\begin{pf}
We must check that the defined total cohomology class
is independent of all choices.
This is simple to verify for the choice of 
the atlas, and the choice of the lifts of 
the transition functions and cocycles.
As stably isomorphic pairs have the same
local description, it is also clear that 
stably isomorphic pairs define the same
total cohomology class.
In detail we give the calculation that exterior 
equivalent triples define the same class:

Let $\rho$ and $\rho'$ be exterior equivalent deckers
on $(P,E)$. Let $c_i:U_i\times G\times G/N\to\U(\HH)$
be such that (E0), (E1) and (E2) of Defintion \ref{DefiOfExteriorEquivalence}
are satisfied.
If $\mu_i$ and $\mu_i'$ are the cocycles for 
the deckers, they both satisfy (\ref{EqLocStucOfDeckers}) for 
the same family of transition functions $\zeta_{ji}$. 
They are related by 
$\mu_i'=\mu_i\ ({\rm Ad}\circ c_i)$, and if 
$\overline\mu_i$ is a lift for $\mu_i$, then
$\overline\mu_i':=\overline\mu_i\ c_i$ defines
a lift for $\mu'_i$.
Let $(\psi_{...},\phi_{..},\omega_.)$ be the cocycle obtained 
from $(\overline\zeta_{ji},\overline\mu_i)$.
Then we have 
\begin{eqnarray*}
&&\overline\mu'_i(u)(h+g,z)\\
&=&\overline\mu_i(u)(h+g,z)\ c_i(u,h+g,z)\\
&=&\overline\mu_i(u)(h,z+gN)\ \overline\mu_i(u)(h,z)\  c_i(u,h+g,z)\ \omega_i(u)(g,h,z)\\
&\stackrel{\rm (E2)}{=}&
\overline\mu_i(u)(h,z+gN)\ \overline\mu_i(u)(g,z)\ 
\mu_i(u)(g,z)^{-1}\big(c_i(u,h,z+gN)\big)\ c_i(g,z)\ \omega_i(u)(g,h,z)\\
&=&
\overline\mu_i(u)(h,z+gN)\   c_i(u,h,z+gN)\ \overline\mu_i(u)(g,z)\ c_i(g,z)\ 
\omega_i(u)(g,h,z)\\
&=&
\overline\mu'_i(u)(h,z+gN)\  \overline\mu'_i(u)(g,z)\  \omega_i(u)(g,h,z),
\end{eqnarray*}
so $\omega_i'=\omega_i$, and
\begin{eqnarray*}
&&\overline\mu'_i(u)(g,z)\\
&=&
\overline\mu_i(u)(g,z)\ c_i(u,g,z)\\
&=&\overline\zeta_{ji}(u)(z+gN)^{-1}
\overline\mu_j(u)(g,g_{ji}(u)+z)
\overline\zeta_{ji}(u)(z)\ c_i(u,g,z)\ \phi_{ji}(u)(g,z)\\
&\stackrel{\rm (E1)}{=}&\overline\zeta_{ji}(u)(z+gN)^{-1}
\overline\mu_j(u)(g,g_{ji}(u)+z)
 c_j(u,g,g_{ji}(u)+z)\ \overline\zeta_{ji}(u)(z)\phi_{ji}(u)(g,z)\\
&=&\overline\zeta_{ji}(u)(z+gN)^{-1}
\overline\mu'_j(u)(g,g_{ji}(u)+z)
  \overline\zeta_{ji}(u)(z)\phi_{ji}(u)(g,z),
\end{eqnarray*} 
so $\phi'_{ji}=\phi_{ji}$.
\end{pf}

\subsection{Dual Pairs and Triples}
\label{SecDualPairsAndTriples}
Of course, the whole discussion we made so far 
for $(G,N)$  can be done for 
$(\widehat G,N^\perp)$. I.e. we can
replace $G$ by its dual group $\widehat G:={\rm Hom}(G,\U(1))$ and $N$ by 
the annihilator $N^\perp:=\{\chi|\chi|_N=1\} \subset \widehat G$ of $N$
everywhere.
This is meaningful as $\widehat G$ is second countable,
$N^\perp$ is discrete and $\widehat G/N^\perp$ compact
(see \ref{SuperSectonGoops}).

\begin{defi}
Let $B$ be a base space.
\begin{enumerate}
\item[$i)$]
A {\bf dual pair} $(\widehat P,\widehat E)$ over $B$ with underlying 
Hilbert space $\HH$ is 
a sequence $\widehat P\to\widehat E\to B$, wherein 
$\widehat E\to B$ is a $\widehat G/N^\perp$-principal fibre
bundel and $\widehat P\to \widehat E$ a $\PU(\HH)$-principal
fibre bundle, such that the latter bundle
is already trivial over the fibres of $\widehat E\to B$.
\item[$ii)$] \label{PageOfAllTheDuals}
A {\bf dual decker } $\hat\rho$ is an action 
$\hat\rho: \widehat P\times \widehat G\to \widehat P$
that lifts the induced $\widehat G$ action on $\widehat E$
and  $\hat\rho(\_,\chi):\widehat P\to\widehat P$ is a bundle 
isomorphisms for all $\chi \in\widehat G$.
\item[$iii)$]
A {\bf dual dynamical triple} $(\hat\rho,\widehat P,\widehat E)$
over $B$ is a pair $(\widehat P,\widehat E)$ over $B$ equipped with
a dual decker $\hat\rho$
\end{enumerate}
\end{defi}

It is clear now how we define 
$\widehat{\rm Par}(B)$, \label{PageOfDualFunctors}
$\widehat{\rm Par}(\widehat E,B),$
$\widehat{\rm Dyn}(B)$,
$\widehat{\rm Dyn}(\widehat E,B)$
and how all statements 
we have achieved so far
translate to dual pairs
and triples.

\subsection{Topological Triples}
\label{SecTopologicalTriples}

We introduce topological triples 
built out of a pair and a dual one. 
They were introduced first in 
\cite{BS} under the name 
T-duality triples in the special
case $G=\R^n, N=\Z$. Our definition
won't be exactly the same as in
\cite{BS} - we comment on this 
in section \ref{SecTheCaseOfRandZ}.
\\\\
There is a canonical
$\U(1)$-principal fibre bundle over
$G/N\times \widehat G/N^\perp$
which is called Poincaré bundle.
We recall its definition. 
Let 
\begin{equation}\label{EqThePoincareBundle}
Q:= \big(G/N\times \widehat G\times \U(1)\big)/N^\perp,
\end{equation}
where the action of $N^\perp$ is defined  by
$(z,\chi,t)\cdot n^\perp:=(z,\chi+n^\perp,t\ \langle n^\perp,z\rangle^{-1})$.
Then the obvious map
$Q\to G/N\times \widehat G/N^\perp$ is a $\U(1)$-principal fibre bundle.
Indeed, $\U(1)$ acts freely and transitive in each fibre by multiplication in the third component and
local sections of $Q\to G/N\times \widehat G/N^\perp$
are given by 
$G/N\times V_a\ni(z, \hat z)\mapsto [(z, \hat s_a(\hat z),1)]\in Q$,
where, by Lemma \ref{LemBorelSectionsForG/N},
$\hat s_a:V_a\to \widehat G$ is a family  of local  sections of the
$N^\perp$-principal bundle $\widehat G\to \widehat G/N^\perp$, 
$V_a\subset \widehat G/N^\perp$.
On the overlap $G/N\times V_{ab}$ two such sections are 
 related by 
\begin{eqnarray*}
 [(z, \hat s_a(\hat z),1)]&=&
 [(z, \hat s_b(\hat z)-n^\perp_{ab}(\hat z),1)]\\
 &=&[(z, \hat s_b(\hat z),\langle n^\perp(\hat z),z\rangle )]\\
 &=&[(z, \hat s_b(\hat z),1)]\cdot \langle n^\perp(\hat z),z\rangle,
 \end{eqnarray*}
wherein  $n^\perp_{ab}(\hat z):=-\hat s_a(\hat z)+\hat s_b(\hat z)$ which
defines a family of transition functions 
$n^\perp_{ab}:V_{ab}\to N^\perp$
for the bundle $\widehat G\to \widehat G/N^\perp.$ 
 Thus we have found that
$\nu^\perp_{ab}:G/N\times V_{ab}\to \U(1)$
defined by $\nu^\perp_{ab}(z,\hat z):= \langle n^\perp_{ab}(\hat z),z\rangle$
are transition functions for the bundle $Q$.

Dually, there is a second $\U(1)$-bundle 
$R:=(G\times \widehat G/N^\perp\times \U(1))/N\to G/N\times\widehat G/N^\perp$ which has 
transition functions 
$\nu_{cd}:W_{cd}\times  \widehat G/N^\perp\to \U(1)$ 
defined by
$\nu_{cd}(z,\hat z):= \langle \hat z, n_{cd}(z)\rangle$
for an open cover $\{W_a\}$ of $G/N$ and 
transition functions 
$n_{cd}:W_{cd}\to N$ 
of $G\to G/N.$

We denote the $\check{\rm C}$ech classes
in $\check H^1(G/N\times \widehat G/N^\perp,\underline{\U(1)})$
which these bundles define by $[Q]$ and $[R]$.

\begin{defi}\label{DefiOfPoincareClass}
The class $\pi:=-[Q]$ constucted above is
called the {\bf Poincaré class} of $G/N\times\widehat G/N^\perp$.
\end{defi}

Of course, in this definition we made a choice,
but up to a sign there is none.

\begin{lem}
$[Q]$ and $[R]$ are inverses of each other, i.e.
$[Q]+[R]=0\in \check H^1(G/N\times \widehat G/N^\perp,\underline{\U(1)})$.
\end{lem}

\begin{pf}
We have to show that $\{\nu_{cd}\ \nu^\perp_{ab}: W_{cd}\times V_{ab}\to\U(1)\}_{(a,c),(b,d)}$
is a $\check{\rm C}$ech coboundary.
Let $\hat s_a:V_a\to \widehat G$ and $ s_c:W_c\to  G$
be  families of local  sections such that
$n_{cd}(z)=- s_c(z)+s_d(z)$ and 
$n^\perp_{ab}(\hat z)=-\hat s_a(\hat z)+\hat s_b(\hat z).$
We show that
$$
\{\nu_{cd}\ \nu^\perp_{ab}:  W_{cd}\times V_{ab}\to\U(1)\}_{(a,c),(b,d)}
=\delta\big( \{\langle  \hat s_a(..), s_c(\_)\rangle:
W_c\times V_a\to\U(1)\}_{(a,c)}\big),
$$
wherin $\delta$ is the usual $\check{\rm C}$ech coboundary operator.
It is $\langle  \hat z, n_{cd}(z)\rangle=\langle  \hat s_a(\hat z), n_{cd}(z)\rangle$
and $\langle n^\perp_{ab}(\hat z),z\rangle=\langle n^\perp_{ab}(\hat z),s_d(z)\rangle$,
because $n_{cd}(z)\in N$ and $n^\perp_{ab}(\hat z)\in N^\perp$; therein the right hand 
side is the pairing $G\times\widehat G\to \U(1)$.
Thus
\begin{eqnarray*}
\nu_{cd}(z,\hat z)\nu^\perp_{ab}(z,\hat z)&=&
\langle  \hat s_a(\hat z), n_{cd}(z)\rangle
\langle n^\perp_{ab}(\hat z),s_d(z)\rangle\\
&=&
\langle  \hat s_a(\hat z), -s_c(z)+s_d(z)\rangle
\langle -\hat s_a(\hat z)+\hat s_b(\hat z),s_d(z)\rangle\\
&=&
\langle \hat s_b(\hat z),s_d(z)\rangle
\langle  \hat s_a(\hat z), s_c(z)\rangle^{-1}
\end{eqnarray*}
which proves the lemma.
\end{pf}

We now turn to the definition of topological 
triples.
Let $P\to E\to B$ be a pair and 
let  $\widehat P\to\widehat E\to B$
be a dual pair with same underlying Hilbert space $\HH$.
We consider the following 
diagram of Cartesian squares 
\begin{equation}\label{DiagAnAlmostTDualityDiagram}
\xymatrix{
&
& P\times_B \widehat E
\ar[rd]\ar[ld]&
&E\times_B\widehat P\ar[rd]\ar[ld]
&\\
&P_{\rm s} \ar[rd]&&
E\times_B \widehat E
\ar[rd]\ar[ld]&& \widehat P\ar[ld]\\
&&E\ar[rd]&&\widehat E\ar[ld]&\\
&&&B&.&
}
\end{equation}
Assume that there is  a $\PU(\HH)$-bundle isomorphism
$\kappa: E\times_B\widehat P\to  P\times_B\widehat E
$
which fits into the above diagram 
(\ref{DiagAnAlmostTDualityDiagram}), i.e.
it fixes its base $E\times_B \widehat E$.
Let us choose a  chart $ U_i\subset B$
common for the pair and the dual pair and
trivialise (\ref{DiagAnAlmostTDualityDiagram}) locally.
For each $u\in U_i$  this induces 
an automorphism $\kappa_i(u)$ of the
trivial $\PU(\HH)$-bundle over \label{PageOfKappai}
$G/N\times \widehat G/N^\perp$,
$$
\xymatrix{
G/N\times \widehat G/N^\perp\times
\PU(\HH)\ar[d]&G/N\times \widehat G/N^\perp\times
\PU(\HH)\ar[l]_{\kappa_i(u)}\ar[d]\\
G/N\times\widehat G/N^\perp&G/N\times\widehat G/N^\perp\ar[l]_=.
}
$$
This automorphism defines a $\check{\rm C}$ech class
$[\kappa_i(u)]\in\check H^1(G/N\times \widehat G/N^\perp,\underline{\U(1)})$
(cp. Therem \ref{ThmClassifacationOfPUBundlesAndAutos}).

\begin{defi}
We say  $\kappa$ {\bf satisfies the Poincaré condition}
if  for each chart $U_i$ and each $u\in U_i$
the equality  $[\kappa_i(u)]= \pi+p_1^*a+ p_2^*b$ holds,
for the Poincaré class $\pi$ and some classes
$a\in \check H^1(G/N,{\underline{\U(1)}})$ and
$b\in \check H^1(\widehat G/N^\perp,{\underline{\U(1)}})$.
Here $p_1,p_2$ are the projections from $G/N\times\widehat G/N^\perp$
on the first and second factor.
\end{defi}
Note that in this definition the classes $a,b$ are just of minor importance.
They are manifestations of the freedom to choose another atlas as 
they vary under the change of  the local trivialisations. 
In fact, one can always modify the local trivialisations of the given atlas
such that $a$ and $b$ vanish.

\begin{defi}\label{DefiTopTriples}
A {\bf topological  triple} $\big(\kappa,(P,E),(\widehat P,\widehat E)\big)$
over $B$ is a pair $(P,E)$  and  dual pair $(\widehat P,\widehat E)$ over $B$
(with same underlying Hilbert space $\HH$) together with a commutative diagram
\begin{equation}\label{DiagTDualityDiagram}
\xymatrix{
& 
& P \times_B \widehat E
\ar[rd]\ar[ld]&
&E\times_B\widehat P\ar[rd]\ar[ld]\ar[ll]_\kappa
&\\
&P \ar[rd]&&
E\times_B \widehat E
\ar[rd]\ar[ld]&\quad\qquad& \widehat P\ar[ld]\\
&&E\ar[rd]&&\widehat E\ar[ld]&\\
&&&B,&&
}
\end{equation}
wherein all squares are Cartesian and $\kappa$ is an isomorphism that satisfies
the Poincaré condition.
\end{defi}

We call two topological triples 
$\big(\kappa,(P,E),(\widehat P,\widehat E)\big)$ 
and
$\big(\kappa',(P',E'),(\widehat P',\widehat E')\big)$
(with underlying Hilbert spaces $\HH,\HH'$ respectively)
{\bf equivalent} if there is  a morphisms of pairs $(\varphi,\vartheta,\theta)$ from
$(P,E)$ to $(P',E')$ and a morphism of 
dual pairs $(\hat\varphi,\hat \vartheta,\hat \theta)$
from $(\widehat P,\widehat E)$
to $(\widehat P',\widehat E')$ such that
the induced diagram
\begin{equation}\label{DiagMorphOfTopTriples}
\xymatrix{
P\times_B\widehat E\ar[d]^{\vartheta\times_{_{B}}\hat\theta}
&&E\times_B\widehat P\ar[ll]_\kappa \ar[d]^{\theta\times_{_{B}}\hat\vartheta}\\
\varphi^*P'\times_B\widehat E'&&E'\times_B\hat\varphi^*\widehat P'\ar[ll]_{\kappa'}
}
\end{equation}
is  commutative
up to homotopy, i.e. the 
$\check{\rm C}$ech class of the bundle automorphism
$(\theta\times_{_{B}}\hat\vartheta)^{-1}\circ
\kappa'^{-1}\circ(\vartheta\times_{_{B}}\hat\theta)\circ
\kappa$ in $\check H^1(E\times_B\widehat P,\underline{\U(1)})$
vanishes.
The triples are called {\bf stably equivalent} 
if the stabilised triples 
$\big(\Eins\otimes\kappa,(P_{\HH_1},E),(\widehat P_{\HH_1},\widehat E)\big)$ and
$\big(\Eins\otimes\kappa',(P'_{\HH_1},E'),(\widehat P'_{\HH_1},\widehat E')\big)$
are equivalent for some separable Hilbert space $\HH_1$. The meaning 
of the index $\HH_1$ is stabilisation as in equation (\ref{EqStablilisedPUBundleP}).
Stable equivalence will be the right choice of equivalence for us, and
we introduce the set valued functor ${\rm Top}$ which \label{PageOfTopFunc}
associates to a base space $B$ the set of stable equivalence classes
of topological triples, i.e.
$$\label{PageTopEB}
{\rm Top}(B):=\{\textrm{topological triples over }B\}\big/_\textrm{stable equivalence}.
$$
If we choose  a $G/N$-bundle $E\to B$ we can consider all the topological
triples with this bundle fixed and 
those stable equivalences for which the identity over $E$ can be extendend 
to a morphisms of the underlying pairs.
We define
$$
{\rm Top}(E,B):=\{\textrm{topological triples with fixed }E\to B\}
\big/_{\textrm{stable equivalence over }{\rm id}_E}. 
$$
The automorphisms ${\rm Aut}_B(E)$ of $E$ over 
the identity of $B$ act on ${\rm Top}(E,B)$ by pullback,
and there is a correspondence 
${\rm Top}(B)\cong\coprod_{[E]} ({\rm Top}(E,B)/{\rm Aut}_B(E))$.

\begin{rem}
We already stated  in Remark \ref{RemThefirstRemOnTriples} that
our notion of topological triples in the case of $G=\R^n$, $N=\Z^n$ and the notion of  T-duality triples as found in \cite{BS}
do not agree.
Nevertheless the two notions lead to the same isomorphism classes, but
we postpone this clarification 
to section \ref{SecTheCaseOfRandZ}, Lemma \ref{LemTopTriplesAreTDualityTriples}.
\end{rem}

By definition, the Poincar\'e class $\pi$ has a geometric 
interpretation in terms of the Poincar\'e bundle (\ref{EqThePoincareBundle}).
For our purposes it will be important that we can give an 
analytical description of $\pi$.

\begin{lem}\label{LemTheDefiOfThePoincareClass}
Choose a Borel section $\sigma:G/N\to G$ 
and an arbitrary section  $\hat\sigma:\widehat G/N^\perp\to \widehat G$
of the corresponding quotient maps.
Then
\begin{enumerate}
\item[$i)$] the map\footnote
{An element $f=f(\_)\in L^\infty(G/N,\U(1))$ 
{\sl is} a multiplication operator on the 
Hilbert space $L^2(G/N)$, see section \ref{SecTheUnitaryAndProjectiveUnitaryGroup}.
}
\begin{eqnarray*} 
 \kappa^\sigma:G/N\times \widehat G/N^\perp&\to&
\PU(L^2(G/N)\otimes\HH)\\
(z,\hat z)&\mapsto& {\rm Ad}\big( \underbrace{\langle\hat\sigma(\hat z),
\sigma(\_-z)-\sigma(\_)\rangle}\otimes\Eins_\HH\big)\\
&&\qquad\qquad\qquad\ \ =:\overline\kappa^\sigma(z,\hat z)\in L^\infty(G/N,\U(1))
\end{eqnarray*}
is continuous and independent of the choice of $\hat\sigma$. Therefore it defines 
(a bundle isomorphism of the trivial $\PU(L^2(G/N,\HH))$-bundle and)
a class 
$$[\kappa^\sigma]\in \check H^1(G/N\times\widehat G/N^\perp,\underline{\U(1)}),
$$
and this class is independent of the choice of $\sigma$;
\item[$ii)$]
$[\kappa^\sigma]=\pi$.
\end{enumerate} 
\end{lem} 
\begin{pf}
$i)$  Firstly, we observe that
$\widehat G\ni \chi\mapsto \langle\chi,\sigma(\_)\rangle\in\U(L^2(G/N))$
is (strongly) continuous. In fact, the sequential continuity of this map\footnote{
$\widehat G$ is first countable.} follows by dominated  convergence.
Therefore the composition
\begin{eqnarray*}
(z,\chi)&\mapsto& \lambda_{_{G/N}}(z)\circ\langle\chi,\sigma(\_)\rangle\circ
\lambda_{_{G/N}}(-z)\circ \langle\chi,-\sigma(\_)\rangle\\
&&=\langle\chi,\sigma(\_-z)-\sigma(\_)\rangle
\end{eqnarray*}
is continous. Here $\lambda_{_{G/N}}$ is the  left regular representation, i.e.
$\lambda_{G/N}(z)F(x):=F(-z+x), F\in L^2(G/N)$.
Now, if $n^\perp\in N^\perp$ we
obtain $\langle n^\perp,\sigma(\_-z)-\sigma(\_)\rangle=
\langle n^\perp,\sigma(-z)\rangle\in\U(1)\subset\U(L^2(G/N)$, because
the difference $\sigma(-z)-(\sigma(\_-z)-\sigma(\_))$ is in $N$.
Thus $(z,\chi)\mapsto{\rm Ad}(\langle\chi,\sigma(\_-z)-\sigma(\_)\rangle)$
factors through the quotient map
$G/N\times\widehat G\to G/N\times \widehat G/N^\perp$ 
which establishes  $\kappa^\sigma$ as stated above.
To see that the class of $\kappa^\sigma$ is independent of 
$\sigma$ it is sufficient  to recognise that $\kappa^n$, defined 
by the same formula as $\kappa^\sigma$,
is unitary implemented for any Borel function $n:G/N\to N$. 
But this is the case, for 
$(z,\chi)\mapsto \langle\chi,n(\_-z)-n(\_)\rangle\in\U(L^2(G/N))$
is continuous and factors through the quotient $G/N\times \widehat G/N^\perp.$
\\\\
$ii)$ Let $\widehat G/N^\perp\supset V_a\stackrel{\hat s_a}{\to}\widehat G$ be a family 
of local sections, so $\hat s_b(\hat z)-\hat s_a(\hat z)=:n_{ab}^\perp(\hat z)\in N^\perp$
defines a set of transition functions.
Let $W_a:=G/N\times V_a$, then
$\kappa_a:W_a\ni(z,\hat z)\mapsto \langle \hat s_a(\hat z),\sigma(\_-z)-\sigma(\_)\rangle
\in\U(L^2(G/N))$ defines locally a continuous unitary lift of $\kappa^\sigma$.
Therefore, on twofold intersections $W_{ba}\ni (z,\hat z)$ we have
$\kappa_b(z,\hat z)\kappa_a(z,\hat z)^{-1}=:\kappa_{ab}(z,\hat z)\in\U(1)$, and
the class of $\kappa^\sigma$ is by definition the class  $[\kappa_{..}]$
of the cocycle $\kappa_{..}$. We obtain
\begin{eqnarray*}
\kappa_b(z,\hat z)\kappa_a(z,\hat z)^{-1}&=&
\langle \hat s_b(\hat z)-\hat s_a(\hat z),\sigma(\_-z)-\sigma(\_)\rangle\\
&=&\langle n^\perp_{ab}(\hat z),\sigma(\_-z)-\sigma(\_)\rangle\\
&=&\langle n^\perp_{ab}(\hat z),\sigma(-z)\rangle\\
&=&\langle n^\perp_{ab}(\hat z),z\rangle^{-1}\in\U(1),
\end{eqnarray*}
wherein the last equality identifies 
$N^\perp$ with $\widehat{G/N}$.
Therefore the class of $\kappa^\sigma$ coincides with
the negative of the class of the bundle $Q$.
\end{pf}

Of course, the situation is symmetric and 
there is an analogous statement  involving 
the class of bundle $R$. In that case,
if we replace everything by its dual counterpart,
we  deal with the function
\begin{eqnarray}\label{EqTheDualAnalyticalExpression} 
\hat\kappa^{\hat\sigma}:G/N\times \widehat G/N^\perp&\to&
\PU(L^2(\widehat G/N^\perp)\otimes\HH)\\
(z,\hat z)&\mapsto& {\rm Ad}( \underbrace{\langle\hat\sigma(..-\hat z)-\hat\sigma(..),
\sigma(z)\rangle}\otimes\Eins_\HH)\nonumber\\
&&\quad\qquad\qquad\quad\ \ 
=:\overline{\hat\kappa}^{\hat\sigma}(z,\hat z)\in L^\infty(\widehat G/N^\perp,\U(1))\nonumber
\end{eqnarray}
with a Borel section $\hat\sigma:\widehat G/N^\perp\to\widehat G$.

\begin{lem}\label{LemThePoincClassForEverythindRep}
The class of $\hat\kappa^{\hat\sigma}$ is the negative of the Poincaré class,
i.e.
$
\pi=[\kappa^{\sigma}]=-[\hat\kappa^{\hat\sigma}].
$
\end{lem}
\begin{pf}
We show that 
$(z,\hat z)\mapsto \kappa^\sigma(z,\hat z)\otimes\hat\kappa^{\hat\sigma}(z,\hat z)$
is  already unitarily  implemented.
Indeed,
\begin{eqnarray*}
&&\kappa^\sigma(z,\hat z)\otimes\hat\kappa^{\hat\sigma}(z,\hat z)\\
&=&
{\rm Ad}(\langle\hat\sigma(\hat z),\sigma(\_-z)-\sigma(\_)\rangle
\langle\hat\sigma(..-\hat z)-\hat\sigma(..),\sigma(z)\rangle)\\
&=&
{\rm Ad}(\langle\hat\sigma(\hat z),\sigma(\_-z)-\sigma(\_)+\sigma(z)\rangle
\langle\hat\sigma(..-\hat z)-\hat\sigma(..),\sigma(z)\rangle)\\
&=&
{\rm Ad}(\langle-\hat\sigma(..-\hat z)+\hat\sigma(..),\sigma(\_-z)-\sigma(\_)+\sigma(z)\rangle
\langle\hat\sigma(..-\hat z)-\hat\sigma(..),\sigma(z)\rangle)\\
&=&
{\rm Ad}(\langle-\hat\sigma(..-\hat z)+\hat\sigma(..),\sigma(\_-z)-\sigma(\_)\rangle)\\
&=&{\rm Ad}(
\lambda_{_{\widehat G/N^\perp}}(\hat z)\lambda_{_{ G/N}}(z)
\langle-\hat\sigma(..),\sigma(\_)\rangle)
\lambda_{_{\widehat G/N^\perp}}(-\hat z) 
\langle\hat\sigma(..),\sigma(\_)\rangle)
\lambda_{_{ G/N}}(-z)\\
&&
\qquad\qquad
\lambda_{_{\widehat G/N^\perp}}(\hat z)\langle-\hat\sigma(..),-\sigma(\_)\rangle)
\lambda_{_{\widehat G/N^\perp}}(-\hat z) 
\langle\hat\sigma(..),-\sigma(\_)\rangle)\\
&&\in \PU(L^2(G/N,\HH)\otimes L^2(\widehat G/N^\perp,\HH)).
\end{eqnarray*}
The argument of  Ad is  a  continuous, unitary expression, 
since the left regular representations
$\lambda_{_{\widehat G/N^\perp}}$ and $\lambda_{_{ G/N}}$
are (strongly) continuous.
\end{pf}

\newpage

\section{T-Duality}
In the last sections we have introduced our main objects
(dynamical and topological  triples).
In the following sections we single out 
specific subclasses of those and
show that they are related to each other.
In addition, we show that their relations
are precisely those which are obtained  from the associated $C^*$-dynamical picture.
Most part of this $C^*$-algebraic structure 
been observed in context
of  continuous trace algebras alone in \cite[Thm 2.2, 
Cor. 2.5]{RR}, \cite[Cor. 2.1]{OR} and more recent
in \cite[Thm. 3.1]{MR} with applications to T-duality.

However, we establish a different approach by use of 
the local structure of the underlying objects.
In particular, we do not 
need any assumption about local compactness
of the underlying spaces. If we assume our base spaces $B$
to be locally compact, then the bundles $E,\widehat E$ over $B$
will be locally compact and 
spectra of  the involved continuous trace algebras; 
but the proof we present is independent of  such an assumption. 
Moreover, it has the advantage of being  explicit enough
to point out the connexion between
the $C^*$-algebraic approach 
and the topological approach to  T-duality.
\\\\
{\bf Notation:} 
In several proofs of the following sections 
we have to check the validity of local identities.
All of these are straight forward computations in 
general, but to keep the formulas readable we 
drop the base variable $u\in B$.
E.g. an identity like 
$$\zeta_{ji}(u)(z+gN)\ \mu_i(u)(g,z)=\mu_j(u)(g,g_{ji}(u)+z)\ \zeta_{ji}(u)(z)
$$
becomes
$$\zeta_{ji}(z+gN)\ \mu_i(g,z)=\mu_j(g,g_{ji}+z)\ \zeta_{ji}(z).
$$
We indicate this with the \label{PageOfNoU} label \ $^\nou$
before any of these computations.

\subsection{The Duality Theory  of Dynamical Triples}
\label{SecTheDualityTheoryOfDynTrip}

We start with the definition of an important subclass of dynamical triples.
\begin{defi}
A dynamical triple is called {\bf dualisable } if
its associated total cohomology class
(Proposition \ref{PropWouldLikeToClassifyTriples})
is of
the form  $[\psi_{...},\phi_{..},\omega_.=d_*\nu_.]$,
i.e. the pair permits a  sufficiently refined  atlas such that 
$\omega_.$ is in the image of the boundary operator $d_*$.
\end{defi}

Let us denote the set of dualisable  dynamical triples
over $B$ by \label{PageOfDualisableFunctor}
${\rm Dyn}^\dag(B)$, respectively ${\rm Dyn}^\dag(E,B)$
for fixed $E$. Similarly, 
$\widehat{\rm Dyn}^\dag(B)$ and $\widehat{\rm Dyn}^\dag(E,B)$
are the sets of dualisable dual dynamical triples.
We are going to present a construction that associates to
a  dualisable  dynamical triple $x$ a dual dualisable dynamical triple
$\widehat x$.
This  construction will be an honest 
map on the level of equivalence classes, and
there it turns out to be involutive, i.e. $[{\widehat{\widehat x}}]=[x]$.
\\\\
Let $(\rho,P,E)$ be a dualisable dynamical triple over
$B$. Choose a sufficiently refined atlas $\{U_i | i\in I\}$ in the
sense that the transition functions $\zeta_{ij}$ and
cocycles $\mu_i$ lift 
continuously to $\overline\zeta_{ij}$ and $\overline\mu_i$ .
Then let $(\psi_{...},\phi_{..},\omega_{.})$ be the associated 
2-cocycle. Since the triple is dualisable, we can assume that
$\omega_.$ is in the image of $d_*$, let $\omega_.:=d_*\nu_.$.
We can modify $\overline\mu_i$
by multiplying with  $\nu_i^{-1}$.
Therefore we can assume without restriction
that $d(\overline\mu_i(u))=\Eins$ and $\omega_i=1$,
for all $i\in I, u\in U_i$. $\overline\mu_i$ is then unique up
to a function $U_i\to Z^1_{\rm Bor}(G,L^\infty(G/N,\U(1)))$. 
\\\\
We are going to define a set of transition
functions  
$$
\hat a_{ij}:=\hat g_{ij}\times\hat \zeta_{ij}:U_{ji}\to\widehat G/N^\perp\ltimes
{\rm Map}(\widehat G/N^\perp,\PU(L^2(G/N)\otimes\HH))
$$
for a dual pair.
The vanishing of $\omega_{.}$ implies that
$\phi_{ij}(u)\in Z^1_{\rm Bor}(G,L^\infty(G/N;\U(1)))$, for all
$i,j, u$. This tells us $^\nou$
\begin{eqnarray*}
\phi_{ij}(g,hN+\_)\ \phi_{ij}(h,\_)&=&\phi_{ij}(g+h,\_)\\
&=&\phi_{ij}(h+g,\_)\\
&=&\phi_{ij}(h,gN+\_)\ \phi_{ij}(g,\_)\in L^\infty(G/N,\U(1))
\end{eqnarray*}
which implies for $g=n\in N$ and all $hN\in G/N$
that
$\phi_{kl}(u)(n,hN+\_)=\phi_{kl}(u)(n,\_)$ holds, hence
$\phi_{kl}(u)(n,\_)\in L^\infty(G/N,\U(1))$ is a constant.
Further,  $n\mapsto\phi_{kl}(u)(n,\_)$ is a homomorphism
and continuous as $N$ is discrete. Thus
there exists
\label{PageThePageOfHatGji} 
$\hat g_{ij} :U_{ji}\to \widehat G/N^\perp (\cong\widehat N)$
such that
\begin{equation}\label{EqTheDualTorusCoycle}
\langle \hat g_{ij}(u),n\rangle=(\phi_{ij}(u)(n,z))^{-1} 
\end{equation}
for all $u\in U_{ji}, n\in N$ and (almost) all $z\in G/N$.

\begin{prop}
 $\{\hat g_{ij} :U_{ji}\to \widehat G/N^\perp | i,j\in I\}$  is a
$\check{\rm C}$ech 1-cocycle.
\end{prop}

\begin{pf}
We have $\delta_g\phi=d_*\psi$, explicitely this reads
for $g\in G, u\in U_{ijk}$
$$
\phi_{jk}(u)(g,\_)
\phi_{ik}(u)(g,\_)^{-1}
\phi_{ij}(u)(g,g_{jk}(u)+\_)
=\psi_{ijk}(u)(\_+gN)
\psi_{ijk}(u)(\_)^{-1}.
$$
For $g=n\in N$ the right-hand side vanishes, hence
the $\check{\rm C}$ech cocycle equation follows.
\end{pf}

In the next lemma we state some properties of 
the lifted cocycles $\overline\mu_i$.
It will be a useful technical tool later.

\begin{lem}\label{LemATechLemmaForMu}
The maps 
\begin{itemize}
\item[$(i)$] 
$U_i\times G/N\ni(u,z)\mapsto \overline\mu_i(u)(\ .\ ,z)|_N\in \U(L^2(N)\otimes\HH),$
\item[$(ii)$] 
$U_i\times G/N\ni(u,z)\mapsto {\rm Ad}(\overline\mu_i(u)(\ .\ ,z))\in \PU(L^2(G)\otimes\HH)$,
\item[$(iii)$]
$U_i\times G/N\ni(u,z)\mapsto {\rm Ad}(\overline\mu_i(u)(-\sigma(\_),z))\in \PU(L^2(G/N)\otimes\HH)$,
\item[$(iv)$]
$U_{ji}\times G/N\ni(u,z)\mapsto {\rm Ad}(\phi_{ij}(u)(\ .\ ,z))\in \PU(L^2(G)),$
\item[$(v)$]
$U_{ji}\times G/N\ni(u,z)\mapsto {\rm Ad}(\phi_{ij}(u)(-\sigma(\_),z))\in \PU(L^2(G/N))$
\end{itemize}
are continuous; and for all $u\in U_{ji},z\in G/N$ the formula
\begin{eqnarray}
&&{\rm Ad}(\overline\zeta_{ji}(u)(-\_)\ \phi_{ji}(u)(-\sigma(\_),0)^{-1}\nonumber)\\
&=&{\rm Ad}(\overline\mu_j(u)(-\sigma(\_+z),g_{ji}(u)+z))\
 \zeta_{ji}(u)(z)\ {\rm Ad}(\overline\mu_i(u)(-\sigma(\_+z),z)^{-1})\nonumber\\
&&{\rm Ad}(\langle\hat\sigma(\hat g_{ji}(u)),\sigma(\_+z)-\sigma(\_)\rangle)\nonumber\\
&&\in\PU(L^2(G/N,\HH))\label{EqLemATechLemmaForMu}
\end{eqnarray}
holds.

\end{lem}
\begin{pf}
$(i)$ Let $n\in N, g\in G$ and $z\in G/N$, then
\begin{eqnarray*}
\overline\mu_i(u)(n,z+gN)\overline\mu_i(u)(g,z)&=&
\overline\mu_i(u)(n+g,z)\\
&=&\overline\mu_i(u)(g,z)\overline\mu_i(u)(n,z).
\end{eqnarray*}
This implies that
$\overline\mu_i(u)(n,z+gN)=\mu_i(u)(g,z)\big(\overline\mu_i(u)(n,z)\big)\in\U(\HH)$.
Therein the right hand side is continuous in $g$, hence the left hand side is 
continuous in $z':=gN$. As the left hand side is symmetric in $z$ and $z'$
it is continuos in $z$.

$(ii)$ Let $(u_\alpha,z_\alpha)\to (u,z)$ be a converging net. Let $x_\alpha:=z_\alpha-z$
and choose $g_\alpha\to 0\in G$ such that $g_\alpha N=x_\alpha$. Such $g_\alpha$
exist -- take a local section of the quotient (Lemma \ref{LemBorelSectionsForG/N}). 
Then 
\begin{eqnarray*}
{\rm Ad}(\overline\mu_i(u_\alpha)(\ .\ ,z_\alpha) )
&=&{\rm Ad}(\overline\mu_i(u_\alpha)(\ .\ ,z+x_\alpha) )\\
&=&{\rm Ad}(\overline\mu_i(u_\alpha)(\ .\ +g_\alpha,z) )
{\rm Ad}(\overline\mu_i(u_\alpha)(g_\alpha,z)^{-1} )\\
&=&{\rm Ad}(\overline\mu_i(u_\alpha)(\ .\ +g_\alpha,z) )
\ \mu_i(u_\alpha)(g_\alpha,z)^{-1} 
\end{eqnarray*}
The second factor converges to 
$\mu_i(u)(0,z)=\Eins$.
For the first factor, note that
$(u,g)\mapsto{\rm Ad}(\overline\mu_i(u)(\ .\ +g,z) )
={\rm Ad}(\lambda_G(g)\overline\mu_i(u)(\ .\ ,z) \lambda_G(-g))
\in \PU(L^2(G,\HH))$
is continuous. This because 
${\rm Bor}(G\times G/N,\U(\HH))\ni \overline\mu_i(u) \mapsto  
\overline\mu_i(u)(\ .\ ,z)\in L^\infty(G,\U(\HH))$
is continuous.

$(iii)$ We shall show that 
$\nu(\ .\ )\mapsto \nu(\sigma(\_))$ 
is a continuous map from the multiplication operators
 $L^\infty(G,\U(\HH))$ to $L^\infty(G/N,\U(\HH))$.
Indeed, let $f\in L^2(G/N)$ and let
$\chi$ be the characteristic function of $\sigma(G/N)\subset G$,
so $(g\mapsto f_\sigma(g):= f(gN)\chi(g))\in L^2(G)$.
Let $\nu_n(\ .\ )\to \nu(\ .\ )\in L^\infty(G,\U(\HH))$ be a converging sequence.
Then
\begin{eqnarray*}
\|\nu_n(\sigma(\_))f(\_)-\nu(\sigma(\_))f(\_)\|^2
&=&\int_{G/N}|\nu_n(\sigma(z))f(z)-\nu(\sigma(z))f(z)|^2\ dz\\
&=&\int_{G}|\nu_n(g)f_\sigma(g)-\nu(g)f_\sigma(g)|^2\ dg\\
&\to& 0,\ {\rm for\ } n\to \infty.
\end{eqnarray*}
$(iv),(v)$ The fourth and fifth  statement follow directly from $(i)$ and  $(ii)$ and from the 
definition of $\phi_{ji}$ on page \pageref{EqThisIsTheDefiOfPhiOhYear}.

Equation (\ref{EqLemATechLemmaForMu}) follows by threefold application of the definition of $\phi_{ji}$.
We let $n(x,y):=\sigma(x+y)-\sigma(x)-\sigma(y)\in N$, then $^\nou$
\begin{eqnarray*}
&&{\rm Ad}\Big(\overline\zeta_{ji}(-\_)\phi_{ji}(-\sigma(\_),0)^{-1}\Big)\\
&=&{\rm Ad}\Big(\overline\mu_j(-\sigma(\_),g_{ji})
	\overline\zeta_{ji}(0)\overline\mu_i(-\sigma(\_),0)^{-1}\Big)\\
&=&{\rm Ad}\Big(\overline\mu_j(-\sigma(\_),g_{ji})
	\overline\mu_j(\sigma(z),g_{ji})^{-1}
	\overline\zeta_{ji}(z)\\
&&	\overline\mu_i(\sigma(z),0)\overline\mu_i(-\sigma(\_),0)^{-1}\Big)\\
&=&{\rm Ad}\Big(\overline\mu_j(-\sigma(\_)-\sigma(z),g_{ji}+z)
	\overline\zeta_{ji}(z)\\
&&	\overline\mu_i(-\sigma(\_)-\sigma(z),z)^{-1}\Big)\\
&=&{\rm Ad}\Big(\overline\mu_j(-\sigma(\_+z),g_{ji}+z)
	\overline\mu_j(n(\_,z),g_{ji}+z)
	\overline\zeta_{ji}(z)\\
&&	\overline\mu_i(n(\_,z),z)^{-1}\overline\mu_i(-\sigma(\_+z),z)^{-1}\Big)\\
&=&{\rm Ad}\Big(\overline\mu_j(-\sigma(\_+z),g_{ji}+z)
	\overline\zeta_{ji}(z)\\
&&	\phi_{ji}(n(\_,z),z)^{-1}
	\overline\mu_i(-\sigma(\_+z),z)^{-1}\Big).
\end{eqnarray*}
Since 
$
{\rm Ad}(	\phi_{ji}(u)(n(\_,z),z)^{-1})=
{\rm Ad}(\langle \hat g_{ji}(u),n(\_,z)\rangle)=
{\rm Ad}(\langle\hat\sigma(g_{ji}(u)),\sigma(\_+z)-\sigma(\_)\rangle),
$
the assertion is proven.
\end{pf}

We now turn to the definition of the projective unitary 
transition functions $\hat\zeta_{ji}$ for the dual pair,
 but before we remark on the local definition we
 make.
 
\begin{rem}
The ad hoc definition of the transition functions of the dual pair
by formula (\ref{EqTheDualPUCocycle}) below may seem very unsatisfactory,
because one cannot even guess the origin of this formula.
In Theorem \ref{ThmTDualityGenaralCase}  we will see that 
the crossed product  $G\times_\rho \Gamma(E,F)$
(see appendix \ref{SecTheSecOnCrossedProducts})
of the associated  $C^*$-algebra of sections $\Gamma(E,F),
F:=P\times_{\PU(\HH)}\K(\HH),$
can be explicitly computed by a fibre-wise, 
modified Fourier transform. The behaviour of 
this transformation under the change of charts will
lead finally to formula (\ref{EqTheDualPUCocycle})
(and shows that the crossed product 
is  isomorphic to
the associated $C^*$-algebra of sections of the
dual pair).
 
However, 
we are in the pleasant situation that we 
can avoid the  $C^*$-algebraic apparatus
at this point and can formulate the theory 
in bundle theoretic terms only. 
\end{rem}

By Lemma \ref{LemATechLemmaForMu}, 
$(u,z)\mapsto {\rm Ad}(\phi_{ji}(u)(-\sigma(\_),z))$ is continuous hence
its restriction to $U_{ji}\times \{0\}$ is:
 $u\mapsto {\rm Ad}(\phi_{ji}(u)(-\sigma(\_),0))$. 
Then we let
\begin{eqnarray}
{\hat\zeta}_{ji}(u)(\hat z)
&:=&{\rm Ad}\Big(
(\overline \kappa^\sigma(-g_{ji}(u),\hat g_{ji}(u)+\hat z)\otimes\Eins_\HH)\
(\lambda_{_{G/N}}(-g_{ji}(u)\otimes\Eins_\HH)\nonumber\\
&&\qquad\qquad
\overline\zeta_{ji}(u)(-\_)\ 
(\phi_{ji}(u)(-\sigma(\_),0)^{-1}\otimes\Eins_\HH)\Big).
\label{EqTheDualPUCocycle}
\end{eqnarray}
Therein 
 $\lambda_{_{G/N}}$ is the left regular representation on $L^2(G/N)$,
and $\overline\kappa^\sigma$ is taken from Lemma \ref{LemTheDefiOfThePoincareClass}.
Indeed, this defines a continuous map 
$$
\hat\zeta_{ji}:U_{ji}\to {\rm Map}(\widehat G/N^\perp,\PU(L^2(G/N)\otimes\HH)),
$$  
but its  definition involved several choices, namely,
the atlas, the liftings $\overline\zeta_{ij},\overline\mu_{i}$ 
and  the section $\sigma$.

\begin{thm}
The family $\hat a_{..}=\{\hat a_{ij}=g_{ji}\times\zeta_{ji} | i,j\in I\}$ is a $\check{\rm  C}$ech cocycle
and its class $[\hat a_{..}]$ is independent of the choices involved.
\end{thm}

\begin{pf}
We have to show  that $\delta_{\hat g}\hat\zeta_{}=\Eins$.
We insert (\ref{EqTheDualPUCocycle}) and obtain $^\nou$
\begin{eqnarray*}
&&{\hat\zeta}_{ji}(\hat z){\hat\zeta}_{ki}(\hat z)^{-1}
{\hat\zeta}_{kj}(\hat g_{ji}+\hat z)\\
&=&
{\rm Ad}\Big(
(\langle\hat\sigma(\hat g_{ji}+\hat z),\sigma(\_+g_{ji})-\sigma(\_)\rangle\otimes\Eins)\\
&&(\lambda_{_{G/N}}(-g_{ji})\otimes\Eins)\ \overline\zeta_{ji}(-\_)\ 
(\phi_{ji}(-\sigma(\_),0)^{-1}\otimes\Eins)\\
&&(\phi_{ki}(-\sigma(\_),0)\otimes\Eins)\ \overline\zeta_{ki}(-\_)^{-1}\ 
(\lambda_{_{G/N}}(g_{ki})\otimes\Eins)\\
&&(\langle\hat\sigma(\hat g_{ki}+\hat z),\sigma(\_+g_{ki})-\sigma(\_)\rangle^{-1}\otimes\Eins)\\
&&(\langle\hat\sigma(\hat g_{kj}+\hat g_{ji}+\hat z),\sigma(\_+g_{kj})-\sigma(\_)\rangle\otimes\Eins)\\
&&(\lambda_{_{G/N}}(-g_{kj})\otimes\Eins)\ \overline\zeta_{kj}(-\_)
(\phi_{kj}(-\sigma(\_),0)^{-1}\otimes\Eins)\Big)\\
&=&
{\rm Ad}\Big(
(\langle\hat\sigma(\hat g_{ji}+\hat z)-\hat\sigma(\hat g_{ki}+\hat z),\sigma(\_+g_{ji})-\sigma(\_)\rangle\otimes\Eins)\\
&& (\psi_{kji}(-\_-g_{ji})\otimes\Eins)
(\phi_{ji}(-\sigma(\_+g_{ji}),0)^{-1}\otimes\Eins)\\
&&(\phi_{ki}(-\sigma(\_+g_{ji}),0)\otimes\Eins)
(\phi_{kj}(-\sigma(\_),0)^{-1}\otimes\Eins)
\Big).
\end{eqnarray*}
The argument of Ad has simplified to a multiplication operator,
so it remains to show that it is a constant in $\U(1)$.
By use of $\delta_g\phi=d_*\psi$ and $d_*\phi=1$ we continue 
\begin{eqnarray*}
\dots&=&{\rm Ad}\Big(
(\langle\hat\sigma(\hat g_{ji}+\hat z)-\hat\sigma(\hat g_{ki}+\hat z),\sigma(\_+g_{ji})-\sigma(\_)\rangle\otimes\Eins)\\
&&\psi_{kji}(0)\ (\phi_{kj}(-\sigma(\_+g_{ji}),g_{ji})\otimes\Eins)
(\phi_{kj}(-\sigma(\_),0)^{-1}\otimes\Eins)\Big)\\
&=&{\rm Ad}\Big(
(\langle\hat\sigma(\hat g_{ji}+\hat z)-\hat\sigma(\hat g_{ki}+\hat z),\sigma(\_+g_{ji})-\sigma(\_)\rangle\otimes\Eins)\\
&&\psi_{kji}(0)\ 
(\phi_{kj}(-\sigma(\_)+\sigma(\_+g_{ji}),0)^{-1}\otimes\Eins)\Big)\\
&=&
{\rm Ad}\Big(
\psi_{kji}(0)\ (\phi_{kj}(-\sigma(\_)+\sigma(\_+g_{ji}),0)^{-1}\otimes\Eins)\\
&&\langle\hat\sigma(\hat g_{ji}+\hat z)-\hat\sigma(\hat g_{ki}+\hat z),
\sigma(g_{ji})\rangle\\
&&(\langle\hat\sigma(\hat g_{ji}+\hat z)-\hat\sigma(\hat g_{ki}+\hat z),
\underbrace{\sigma(\_+g_{ji})-\sigma(\_)-\sigma(g_{ji})}\rangle\otimes\Eins)\Big),\\
&&\qquad\qquad\qquad\qquad\qquad\qquad\qquad\qquad\ \ 
 \in N
\end{eqnarray*}
this can be transformed by use of the definition of $\hat g_{ji}$
\begin{eqnarray}
\dots&=&{\rm Ad}\Big(
\psi_{kji}(0)\ (\phi_{kj}(-\sigma(\_)+\sigma(\_+g_{ji}),0)^{-1}\otimes\Eins)\nonumber\\
&&
\langle\hat\sigma(\hat g_{ji}+\hat z)-\hat\sigma(\hat g_{ki}+\hat z),
\sigma(g_{ji})\rangle\nonumber\\
&&(\langle -\hat g_{kj},
\sigma(\_+g_{ji})-\sigma(\_)-\sigma(g_{ji})\rangle\otimes\Eins)\Big)\nonumber\\
&=&{\rm Ad}\Big(
\psi_{kji}(0) \phi_{kj}(\sigma(g_{ji}),0)^{-1}\
\langle\hat\sigma(\hat g_{ji}+\hat z)-\hat\sigma(\hat g_{ki}+\hat z),
\sigma(g_{ji})\rangle\Big)\nonumber\\
&=&\Eins_{L^2(G/N)\otimes\HH}.
\end{eqnarray}
This is the desired result.

To check that all the choices involved have no effect on the
class of this cocycle is straight forward, and we skip the
tedious computation here. 
\end{pf}

Let us define a pair for $\hat g_{ij},\hat\zeta_{ij}$ explicitly. Let \label{PageDualOfEanfFdd}
$$
\widehat E:=\coprod_i(U_i\times \widehat G/N^\perp)\big/_\sim
$$
with relation $(i,u,\hat z)\sim (j,u,\hat g_{ji}(u)+\hat z)$ and 
$$
\widehat P:=\coprod_i(U_i\times \widehat G/N^\perp\times \PU(L^2(G/N)\otimes\HH))\big/_\sim
$$
with relation $(i,u,\hat z,U)\sim (j,u,\hat g_{ji}(u)+\hat z,\hat\zeta_{ji}(u)(\hat z)U)$.
$(\widehat P,\widehat E)$ is a dual pair over $B$ with underlying 
Hilbert space $\widehat \HH:=L^2(G/N)\otimes\HH$. 
Then the statement of the theorem is that we have constructed a map
\begin{eqnarray*}
{\rm Dyn}^\dagger(B)&\to&\widehat{\rm Par}(B),\\
 {[(\rho,P,E)]}&\mapsto& [(\widehat P,\widehat E)]
\end{eqnarray*}
for one easily checks that outer conjugate triples define
isomorphic duals.
We are going to  improve this statement in the next theorem.

Let us define a family of projective unitary 1-cocycles
by the following simple formula. Let $\chi\in \widehat G,
\hat z\in\widehat G/N^\perp, u\in U_i,$ then
we define
\begin{eqnarray}\label{EqTheDualDeckerCocycle}
\hat\mu_i(u)(\chi,\hat z):=
{\rm Ad}( \langle\chi,-\sigma(\_)\rangle\otimes\Eins_\HH)\in \PU(L^2(G/N)\otimes\HH).
\end{eqnarray}
Therein $\sigma:G/N\to G$ is the same Borel section
which we have used to define $\hat\zeta_{ji}$.
By dominated convergence\footnote{
$\widehat G$ is first countable.}, 
it is clear that 
$\chi\mapsto\langle\chi,-\sigma(\_)\rangle\in L^\infty(G/N,\U(1))$
is continuous, so $\hat\mu_i$ is.

\begin{thm}\label{ThmTheMainTheoremOfDualityTheoryOfDynTriples}
The family $\{\hat\mu_i |i\in I\}$ defines a dual decker $\hat\rho$
on the pair $(\widehat P,\widehat E)$, and we obtain 
a bijection  from the set of dualisable dynamical triples to the set of dualisable dual dynamical
 triples,
\begin{eqnarray*}
{\rm Dyn}^\dagger(B)&\to&\widehat{\rm Dyn}^\dagger(B).\\
{[\rho,P,E]}&\mapsto&[\hat\rho,\widehat P,\widehat E]
\end{eqnarray*}
In formulas 
 (\ref{EqTheDualTorusCoycle}),(\ref{EqTheDualPUCocycle}) and (\ref{EqTheDualDeckerCocycle}) we can replace everything by its dual counterpart, 
 i.e. we interchange the rôle of triples and dual triples,  then
 we obtain a map 
$${\widehat{\rm Dyn}}^\dagger(B)\to{\rm Dyn}^\dagger(B).
$$
Moreover, these two maps are natural and inverse to each other, so we 
have an equivalence of functors
$$
{\rm Dyn}^\dag \cong \widehat{\rm Dyn}^\dag.
$$
\end{thm}
\begin{pf}
To see that the cocycles $\hat\mu_i$ define a decker we have
to verify that
$$
\hat\mu_i(u)(\chi,\hat z)=\hat\zeta_{ji}(u)(\hat z+\chi N^\perp)^{-1}
\hat\mu_j(u)(\chi,\hat g_{ji}(u)+\hat z)\hat\zeta_{ji}(u)(\hat z).
$$
We just compute the right-hand side $^\nou$
\begin{eqnarray*}
&&\hat\zeta_{ji}(\hat z+\chi N^\perp)^{-1}
\hat\mu_j(\chi,\hat g_{ji}+\hat z)\hat\zeta_{ji}(\hat z)\\
&=&{\rm Ad}\Big(
\overline\zeta_{ji}(-\_)^{-1}\ (\phi_{ji}(-\sigma(\_),0)\otimes\Eins)\ 
(\lambda_{_{G/N}}(g_{ji})\otimes\Eins)\\
&&(\langle\hat\sigma(\hat z+\chi N^\perp+\hat g_{ji}),\sigma(\_)-\sigma(\_+g_{ji})\rangle\otimes\Eins)\\
&&(\langle\chi,-\sigma(\_)\rangle\otimes\Eins)\\
&&(\langle\hat\sigma(\hat z+\hat g_{ji}),-\sigma(\_)+\sigma(\_+g_{ji})\rangle\otimes\Eins)\\
&&
(\lambda_{_{G/N}}(-g_{ji})\otimes\Eins)\ \overline\zeta_{ji}(-\_)\ 
(\phi_{ji}(-\sigma(\_),0)^{-1}\otimes\Eins)
\Big)\\
&=&{\rm Ad}\Big(
(\langle\hat\sigma(\hat z+\chi N^\perp+\hat g_{ji}),\sigma(\_-g_{ji})-\sigma(\_)\rangle\otimes\Eins)\\
&&(\langle\chi,-\sigma(\_-g_{ji})\rangle\otimes\Eins)\\
&&(\langle\hat\sigma(\hat z+\hat g_{ji}),-\sigma(\_-g_{ji})+\sigma(\_)\rangle\otimes\Eins)
\Big)\\
&=&{\rm Ad}\Big(
(\langle\chi,-\sigma(\_)\rangle
\langle\hat\sigma(\hat z+\hat g_{ji}+\chi N^\perp)-\chi-\hat\sigma(\hat z+\hat g_{ji}),
\sigma(\_-g_{ji})-\sigma(\_)\rangle\otimes\Eins)
\Big)\\
&=&{\rm Ad}\Big(
(\langle\chi,-\sigma(\_)\rangle\otimes\Eins)\
\underbrace{
\langle\hat\sigma(\hat z+\hat g_{ji}+\chi N^\perp)-\chi-\hat\sigma(\hat z+\hat g_{ji}),
-\sigma(g_{ji})\rangle}
\Big)\\
&=&
\hat\mu_i(\chi,\hat z).
\qquad\qquad\qquad\qquad\qquad\qquad\qquad\quad 
=:(\hat\phi_{ji}(\chi,\hat z))^{-1}\equiv\Eins\ {\rm mod}\ \U(1)
\end{eqnarray*}
This establishes $\hat\rho$. \label{PageOfDualPhiTotalCocyc}
A careful look at this calculation shows that $\hat\phi_{..}$ indeed is
the second term of the total dual cocycle 
$(\hat\psi_{...},\hat\phi_{..},1)$ defined by our constructed dual
dynamical triple.
It is  easy to see that  another choice of the 
section $\sigma$ alters $\hat\mu_i$ and $\hat \zeta_{ji}$ 
precisely in such a way that they define an isomorphic pair.
So we have established a map
${\rm Dyn}^\dag(B)\to \widehat{\rm Dyn}^\dag(B)$
and, by replacing everything by its dual, a map in  opposite 
direction.

To prove the remaining assertions of the theorem, 
we shall apply our construction twice.
We will find that the double dual
$(\hat{\hat \rho},\widehat{\widehat P},\widehat{\widehat E})$
is isomorphic to $({\rm Ad}\circ\lambda_G\otimes\rho,\PU(L^2(G))\otimes P,E)$.
In particular, we already see from the definition of $\hat\phi_{ji}$
that the double dual bundle
$\widehat{\widehat E}$ has cocycle 
$\hat{\hat g}_{ji}(u) :=\hat\phi_{ji}(u)(\_,0)|_{N^\perp}=g_{ji}(u)$,
so there is $\theta:\widehat{\widehat E}\cong E$.
We compute the double dual cocycle 
$$
\hat{\hat\zeta}_{ji}:U_{ij}\to{\rm Map}(G/N,\PU(L^2(\widehat G/N^\perp)\otimes L^2(G/N)\otimes\HH)).
$$  
By definition, we have
by (\ref{EqTheDualPUCocycle}) and (\ref{EqTheDualAnalyticalExpression}) $^\nou$
\begin{eqnarray*}
&&\hat{\hat\zeta}_{ji}(z)\\
&=&
{\rm Ad}\Big(
(\overline{\hat\kappa}^{\hat\sigma}(g_{ji}+z,-\hat g_{ji})\otimes\Eins\otimes\Eins)
(\lambda_{_{\widehat G/N^\perp}}(-\hat g_{ji})\otimes\Eins\otimes\Eins)\\ 
&&\overline{\hat\zeta}_{ji}(-..)
(\hat\phi_{ji}(-\hat\sigma(..),0)^{-1}\otimes\Eins\otimes\Eins)\Big)\\
&=&{\rm Ad}\Big(
(\overline{\hat\kappa}^{\hat\sigma}(g_{ji}+z,-\hat g_{ji})\otimes\Eins\otimes\Eins)
(\lambda_{_{\widehat G/N^\perp}}(-\hat g_{ji})\otimes\Eins\otimes\Eins)\\
&& (\overline{\kappa}^{\sigma}(-g_{ji},\hat g_{ji}-..)\otimes\Eins)
(\Eins\otimes\lambda_{_{G/N}}(-g_{ji})\otimes\Eins)\\
&&(\Eins\otimes\overline{\zeta}_{ji}(-\_))
(\Eins\otimes \phi_{ji}(-\sigma(\_),0)^{-1}\otimes\Eins)
(\hat\phi_{ji}(-\hat\sigma(..),0)^{-1}\otimes\Eins\otimes\Eins)\Big).\\
\end{eqnarray*}
With equation (\ref{EqLemATechLemmaForMu})  and 
the definitions of $\overline\kappa^\sigma,\phi_{ji},\dots$ this 
reads
\begin{eqnarray*}
\dots&=&{\rm Ad}\Big(
	(\langle \hat\sigma(..+\hat g_{ji})-\hat\sigma(..),\sigma(g_{ji}+z)\rangle\otimes\Eins		\otimes\Eins)\\
&&	(\lambda_{_{\widehat G/N^\perp}}(-\hat g_{ji})\otimes\Eins\otimes\Eins)
	(\langle\hat\sigma(\hat g_{ji}-..),\sigma(\_+g_{ji})-\sigma(\_)\rangle\otimes\Eins)\\
&&	(\Eins\otimes\lambda_{_{G/N}}(-g_{ji})\otimes\Eins)
	(\Eins\otimes\overline\mu_j(-\sigma(\_+z),g_{ji}+z))\Big)\\
&&	(\Eins\otimes\Eins\otimes\zeta_{ji}(z))\ {\rm Ad}\Big((\Eins\otimes\overline\mu_i(- \sigma(\_+z),z))
	(\Eins\otimes\langle\hat\sigma(\hat g_{ji}),\sigma(\_+z)-\sigma(\_)\rangle\otimes\Eins)\\
&&	(\langle\hat\sigma(\hat g_{ji}-..)+\hat\sigma(..)-\hat\sigma(\hat g_{ji},
	-\sigma(g_{ji}+z)+\sigma(z)\rangle\otimes\Eins\otimes\Eins)\Big),
\end{eqnarray*}
and after some intermediate steps we find
\begin{eqnarray*}
\hat{\hat\zeta}_{ji}(u)(z)=
	\eta_j(u)(g_{ji}(u)+z)^{-1}\
	{\rm Ad}\Big(\lambda_{_{\widehat G/N^\perp}}(-\hat g_{ji}(u))\Big)
	\otimes\Eins\otimes
	\zeta_{ji}(u)(z)\
	\eta_i(u)(z),
\end{eqnarray*}
wherein 
\begin{eqnarray*}
\eta_i(u)(z)&:=&{\rm Ad}\Big(
(\langle\hat\sigma(..)+\hat\sigma(-..),\sigma(z)\rangle
\langle\hat\sigma(-..),\sigma(\_-z)-\sigma(\_)\rangle\otimes\Eins\otimes\Eins)\\
&&\qquad(\Eins\otimes\overline\mu_i(u)(-\sigma(\_),z)^{-1})
(\Eins\otimes\lambda_{_{G/N}}(z)\otimes\Eins)\Big).
\end{eqnarray*}
This already proves part of the second half of the theorem, for 
we see that the double dual is isomorphic to a pair with transition
functions 
$$
(u,z)\mapsto{\rm Ad}(\lambda_{_{\widehat G/N^\perp}}(-\hat g_{ji}(u)))\otimes\Eins\otimes
\zeta_{ji}(u)(z).
$$ 
For general reasons such a pair is isomorphic to a pair 
with transition function $\zeta_{ji}$,
for 
$$
(u,z)\mapsto \lambda_{_{\widehat G/N^\perp}}(-\hat g_{ji}(u))\otimes\Eins\otimes\Eins
\in\U(L^2(\widehat G/N^\perp)\otimes L^2(G/N)\otimes\HH)
$$
is  a unitary cocycle.
However, we can give a concrete isomorphism. Let us denote \label{PageFouriertrrrr}
by $\mathcal F:L^2(\widehat G/N^\perp)\to L^2(N)$ the Fourier 
transform, then we obtain a multiplication operator
\begin{eqnarray*}
\mathcal F\circ \lambda_{_{\widehat G/N^\perp}}(-\hat g_{ji}(u))\circ
\mathcal F^{-1}=\langle \hat g_{ji},\ .\ \rangle^{-1}=\phi_{ji}(u)(\ .\ ,0)|_N
\in \U(L^2(N)).
\end{eqnarray*}
Thus, by definition of $\phi_{ji}$ we have
\begin{eqnarray*}
&&{\rm Ad}(\lambda_{_{\widehat G/N^\perp}}(-\hat g_{ji}(u)))\otimes\Eins\otimes
\zeta_{ji}(u)(z)\\
&=&{\rm Ad}\Big( (\mathcal F^{-1}\otimes\Eins\otimes\Eins)
\circ \overline\mu_j(u)(\ .\ ,g_{ji}(u)+z)^{-1}|_N\circ
(\mathcal F\otimes\Eins\otimes\Eins)\Big)\\
&&(\Eins\otimes\Eins\otimes\zeta_{ji}(u)(z))\
{\rm Ad}\Big( (\mathcal F^{-1}\otimes\Eins\otimes\Eins)
\circ \overline\mu_i(u)(\ .\ ,z)|_N\circ(\mathcal F\otimes\Eins\otimes\Eins)\Big),
\end{eqnarray*}
and we let
$$
\vartheta_i(u)(z):= {\rm Ad}\Big( (\mathcal F^{-1}\otimes\Eins\otimes\Eins)
\circ \overline\mu_i(u)(\ .\ ,z)|_N\circ(\mathcal F\otimes\Eins\otimes\Eins)\Big)\ \eta_i(u)(z).
$$
Due to Lemma \ref{LemATechLemmaForMu} $(i)$ $\vartheta_i$ is continuous.
We obtain
$$
\hat{\hat\zeta}_{ji}(u)(z)=\vartheta_j(u)(g_{ji}(u)+z)^{-1}\ 
\big(\Eins\otimes\Eins\otimes \zeta_{ji}(u)(z)\big)\ \vartheta_i(u)(z),
$$
so  $\vartheta_i$ define local isomorphisms between 
$\widehat{\widehat P}$ and $\PU(L^2(\widehat G/N^\perp)\otimes L^2(G/N))\otimes P$ which fit together to a global isomorphism 
$$
(\Eins,\vartheta,\theta):(\widehat{\widehat P},\widehat{\widehat E})\to
(\PU(L^2(\widehat G/N^\perp)\otimes L^2(G/N))\otimes P,E).
$$
We have to compute the behaviour of the double dual decker  under  
this isomorphism, i.e. we have to compute
$$
\mu'_i(u)(g,z):=\vartheta_i(u)(z+gN)\ \hat{\hat\mu}_i(u)(g,z)) \vartheta_i(u)(z)^{-1},
$$
for $
\hat{\hat\mu}_i(u)(g,z)={\rm Ad}(\langle\hat\sigma(..),-g\rangle\otimes\Eins\otimes\Eins).
$
This yields\label{PageThisPageContainsExplicitMaterial} $^\nou$
\begin{eqnarray*}
\mu'_i(g,z)&=&{\rm Ad}\Big(
(\mathcal F^{-1}\otimes\Eins\otimes\Eins)
\circ \overline\mu_i(\ .\ ,z+gN)|_N\circ(\mathcal F\otimes\Eins\otimes\Eins)\\
&&(\Eins\otimes\lambda_{_{G/N}}(gN)\otimes\Eins)
(\langle\hat\sigma(..),{-g-\sigma(\_)+\sigma(\_+gN)}\rangle\otimes\Eins)\\
&&(\Eins\otimes\overline\mu_i(\sigma(\_+gN)-\sigma(\_),z))\\
&&(\mathcal F^{-1}\otimes\Eins\otimes\Eins)\circ \overline\mu_i(\ .\ ,z)^{-1}|_N\circ
(\mathcal F\otimes\Eins\otimes\Eins)\Big)\\
&=&{\rm Ad}\Big(
(\mathcal F^{-1}\otimes\Eins\otimes\Eins)\circ \overline\mu_i(\ .\ ,z+gN)|_N\\
&&(\Eins\otimes\lambda_{_{G/N}}(gN)\otimes\Eins) 
(\lambda_{_{N}}(g+\sigma(\_)-\sigma(\_+gN))\otimes\Eins)\\
&&(\Eins\otimes\overline\mu_i(\sigma(\_+gN)-\sigma(\_),z))\\
&& \overline\mu_i(\ .\ ,z)^{-1}|_N\circ(\mathcal F\otimes\Eins\otimes\Eins)\Big)\\
&=&{\rm Ad}\Big(
(\mathcal F^{-1}\otimes\Eins\otimes\Eins)\circ 
(\Eins\otimes\lambda_{_{G/N}}(gN)\otimes\Eins)\\
&&( \lambda_{_{N}}(g+\sigma(\_)-\sigma(\_+gN))\otimes\Eins)
\circ(\mathcal F\otimes\Eins\otimes\Eins)\\
&&(\Eins\otimes\Eins\otimes\overline\mu_i(g,z))\Big).
\end{eqnarray*}
The last equality is just the cocycle condition for $\overline\mu_i$.
We can make a further manipulation by use of 
the following isomorphism 
$$
\mathcal S:
 L^2(N)\otimes L^2(G/N)\cong L^2(N\times G/N)\to L^2(G).
$$ 
We let\label{PageThePageOfS}
$(\mathcal S(f))(g):=f(g-\sigma(gN),gN)$. Its inverse is given
by $(\mathcal S^{-1}(f))(n,z):=f(n+\sigma(z))$, and it is immediate to 
verify that 
$$
\mathcal S^{-1}\circ\lambda_{_{G}}(g)\circ\mathcal S=
\lambda_{_{G/N}}(gN) \lambda_{_{N}}(g+\sigma(\_)-\sigma(\_+gN)).
$$
This implies for the cocycle that
\begin{eqnarray*}
\mu_i'(u)(g,z)&=&
{\rm Ad}\Big(
(\mathcal F^{-1}\otimes\Eins\otimes\Eins)
\circ(\mathcal S\otimes\Eins)\circ
(\lambda_{_{G}}(g)\otimes \overline\mu_i(u)(g,z))\\
&&\circ
(\mathcal S\otimes\Eins)\circ(\mathcal F\otimes\Eins\otimes\Eins)\Big)
\end{eqnarray*}
To summarise, we have shown that there is an isomorphism of 
dynamical triples 
$$
\big((\mathcal S\otimes\Eins)\circ(\mathcal F\otimes\Eins\otimes\Eins),\vartheta, \theta\big):
(\hat{\hat\rho},\widehat{\widehat P},\widehat{\widehat E})\to
(({\rm Ad}\circ\lambda_{_G})\otimes\rho,\PU(L^2(G))\otimes P,E).
$$
So, as we discussed in Example \ref{PageTheExampleOfStableEquivalence},
the double dual triple is outer conjugate to the  triple we started with. 

We finally comment on the naturality of the defined 
maps. Let $f:B'\to B$ be a map of bases spaces.
Then we have a diagram
$$
\xymatrix{
{\rm Dyn}^\dag(B)\ar[r]\ar[d]^{f^*}&\widehat{\rm Dyn}^\dag(B)\ar[d]^{f^*}\\
{\rm Dyn}^\dag(B')\ar[r]&\widehat{\rm Dyn}^\dag(B').
}
$$
On the local level pullback with $f^*$  
is the purely formal substitution of $u\in B$
by $f(u')$ for $u'\in B'$ in all formulas, and it follows 
that the diagram commutes.

This proves the theorem.
\end{pf}

\subsection{The Relation to Topological T-Duality}
\label{SecTheRelationToTopologicalTDuality}
So far we have found that dualisable dynamical triples and 
dualisable dual dynamical triples have the same isomorphism 
classes. 
We  demonstrate that the theory developed so far is
intimately  connected to topological T-duality.
\\\\
Let $(\rho,P,E)$ be dualisable and $\widehat P,\widehat E$ as 
above.
Let \label{PageOfStableP} $P_{\rm top}:=\PU(L^2(G/N))\otimes P.$

\begin{thm}\label{ThmTheTopologicalisationMapTau}
The two pairs $(P_{\rm top},E)$ and  $(\widehat P,\widehat E)$ 
with underlying Hilbert space $\widehat \HH=L^2(G/N)\otimes\HH$
span a  
topological triple $(\kappa^{\rm top},(P_{\rm top},E),(\widehat P,\widehat E))$, 
and we obtain a map
$$
\tau(B):{\rm Dyn}^\dagger(B)\to{\rm Top}(B).
$$
Moreover, this map is natural, i.e. it defines  a natural transformation of functors
$$
\tau:{\rm Dyn}^\dag\to {\rm Top}.
$$
\end{thm}

\begin{pf}
We have to show that there
exists (a natural choice of) 
an isomorphism $\kappa^{\rm top}:E\times_B\widehat P\to P_{\rm top}\times_B \widehat E$
which satisfies the Poincaré condition.
We are going to define  $\kappa^{\rm top}$ by local isomorphisms of
the locally trivialised pairs.
Let
\begin{eqnarray*}
\kappa^{\rm top}_i:U_i\to {\rm Map}( G/N\times\widehat G/N^\perp,\PU(L^2(G/N)\otimes\HH))
\end{eqnarray*}
be given by the  formula
\begin{eqnarray}
&&\kappa^{\rm top}_i(u)(z,\hat z):={\rm Ad}\Big(  
\underbrace{(
\overline\kappa^\sigma(z,\hat z)\otimes\Eins)\ 
\overline\mu_i(u)(-\sigma(\_),z)^{-1}\ 
(\lambda_{_{G/N}}(z)\otimes\Eins)}\Big)\nonumber\\
&&
\qquad\qquad\qquad
\qquad\qquad\qquad
\qquad\qquad\qquad
\quad
=:\overline\kappa^{\rm top}_i(u)(z,\hat z)\label{EqTheLocalDefinitionOfKappaTop}
\end{eqnarray}
 $\kappa^{\rm top}_i$ is continuous, for
the first factor  is exactly $\kappa^\sigma$ of 
 Lemma \ref{LemTheDefiOfThePoincareClass} which defines
the Poincaré class, thus it is continuous. The left regular 
representation $\lambda_{_{G/N}}$ is continuous and
by Lemma \ref{LemATechLemmaForMu} $(iii)$
$(u,z)\mapsto {\rm Ad}(\overline\mu_i(u)(-\sigma(\_),z) )$
is continuous. 
The functions 
$(u,z,\hat z)\mapsto 
\Eins\otimes \zeta_{ji}(u)(z)$
and 
$(u,z,\hat z)\mapsto 
\Eins\otimes \hat\zeta_{ji}(u)(\hat z)$
are transition functions for the bundles
$P_{\rm dyn}\times_B\widehat E$ and $E\times_B\widehat P$, and
the functions $\kappa_i^{\rm top}$ will define a 
global isomorphism if and only if
$$
\kappa^{\rm top}_j(u,g_{ji}(u)+z,\hat g_{ji}(u)+\hat z)\ \hat\zeta_{ji}(u)(\hat z)\  \kappa^{\rm top}_i(u,z,\hat z)^{-1}=\Eins\otimes\zeta_{ji}(u)(z),
$$
for $(u,z,\hat z)\in U_{ij}\times G/N\times \widehat G/N^\perp.$
By construction, the Poincaré condition 
will be satisfied automatically.

We are able to do this calculation on the projective unitary level directly, but
later we need the  same calculation on the unitary level, so we
prepare ourselves first with a choice of lifts $\overline{\hat\zeta}_{ji}$.

Recall that $u\mapsto {\rm Ad}(\phi_{ji}(u)(-\sigma(\_),0))$ is continuous.
If we pass once more to a refined atlas which again we denote by $\{U_i\}$, 
then there exist continuous lifts
$$
\overline\phi_{ji}:U_{ji}\to L^\infty(G/N,\U(1))
$$ 
such that
${\rm Ad}(\overline\phi_{ji}(u))={\rm Ad}(\phi_{ji}(u)(-\sigma(\_),0))$;
by those we define lifted transition functions functions
$\overline{\hat\zeta}_{ji}:U_{ji}\to L^\infty(\widehat G/N^\perp,L^\infty( G/N,\U(\HH)))$ by
\begin{eqnarray}
\overline{\hat\zeta}_{ji}(u)(\hat z)&:=&
(\overline\kappa^\sigma(-g_{ji}(u),\hat g_{ji}(u)+\hat z)\otimes\Eins_\HH)\ 
(\lambda_{_{G/N}}(-g_{ji}(u))\otimes\Eins_\HH)\nonumber\\
&&\qquad\qquad
\overline\zeta_{ji}(u)(-\_)\ 
(\overline\phi_{ji}(u)^{-1}\otimes\Eins_\HH),\label{EqTheDualLiftedPUCocycle}
\end{eqnarray}
for  $u\in U_{ij},\hat z\in \widehat G/N^\perp.$

Now, we calculate $^\nou$
\begin{eqnarray}
&&\overline\kappa^{\rm top}_j(g_{ji}+z,\hat g_{ji}+\hat z)\  \overline{\hat\zeta}_{ji}(\hat z)\  \overline\kappa^{\rm top}_i(z,\hat z)^{-1}\nonumber\\
&=&
(\langle \hat\sigma(\hat g_{ji}+\hat z),
-\sigma(\_)+\sigma(\_-g_{ji}-z)\rangle\otimes\Eins)\ 
\overline\mu_j(-\sigma(\_),g_{ji}+z)^{-1}\nonumber\\
&&(\lambda_{_{G/N}}(g_{ji}+z)\otimes\Eins)\nonumber\\ 
&&(\langle \hat\sigma(\hat g_{ji}+\hat z), -\sigma(\_)+\sigma(\_+g_{ji})\rangle\otimes\Eins)\
(\lambda_{_{G/N}}(-g_{ji})\otimes\Eins)\nonumber\\
&&\overline\zeta_{ji}(-\_)\ (\overline\phi_{ji}^{-1}\otimes\Eins)\nonumber\\
&&
(\lambda_{_{G/N}}(-z)\otimes\Eins)\ 
\overline\mu_i(-\sigma(\_),z)\ 
(\langle \hat \sigma(\hat z) ,\sigma(\_)-\sigma(\_-z)\rangle\otimes\Eins) \nonumber\\
&=&
(\langle \hat\sigma(\hat g_{ji}+\hat z)-\hat\sigma(\hat z),
-\sigma(\_)+\sigma(\_-z)\rangle\otimes\Eins)\nonumber\\
&&\overline\mu_j(-\sigma(\_),g_{ji}+z)^{-1}
\overline\zeta_{ji}(z-\_)\
(\overline\phi_{ji}(\_-z)^{-1}\otimes\Eins)
\overline\mu_i(-\sigma(\_),z)\ 
 \nonumber\\
&=&
(\langle \hat\sigma(\hat g_{ji}+\hat z)-\hat\sigma(\hat z),
-\sigma(\_)+\sigma(\_-z)+\sigma(z)\rangle\otimes\Eins)\ 
\langle \hat\sigma(\hat g_{ji}+\hat z)-\hat\sigma(\hat z),
-\sigma(z)\rangle\nonumber\\
&&\overline\mu_j(-\sigma(\_))(u,g_{ji}+z)^{-1}
\overline\zeta_{ji}(z-\_)\
(\overline\phi_{ji}(\_-z)^{-1}\otimes\Eins)
\overline\mu_i(-\sigma(\_),z)\ 
 \nonumber\\
&=&
(\langle \hat g_{ji},
-\sigma(\_)+\sigma(\_-z)+\sigma(z)\rangle\otimes\Eins)\ 
\langle \hat\sigma(\hat g_{ji}+\hat z)-\hat\sigma(\hat z),
-\sigma(z)\rangle\nonumber\\
&&\overline\mu_j(-\sigma(\_),g_{ji}+z)^{-1}
\overline\zeta_{ji}(z-\_)\
(\overline\phi_{ji}(\_-z)^{-1}\otimes\Eins)
\overline\mu_i(-\sigma(\_),z).\nonumber
\end{eqnarray} 
The variable $\hat z$ only occurs in  
$\langle \hat\sigma(\hat g_{ji}+\hat z)-\hat\sigma(\hat z),
-\sigma(z)\rangle$.
To the remaining expressinon we can apply  
equation (\ref{EqLemATechLemmaForMu}) 
of Lemma \ref{LemATechLemmaForMu} which 
is valid on the unitary level up to a $\U(1)$-valued perturbation 
$\alpha'_{ji}:U_{ji}\to L^\infty(G/N,\U(1))$.
We obtain
\begin{eqnarray}
\dots&=&
\underbrace{\langle \hat\sigma(\hat g_{ji}(u)+\hat z)-\hat\sigma(\hat z),-\sigma(z)\rangle
\ \alpha'_{ji}(u)(z)^{-1}}\
(\Eins\otimes\overline\zeta_{ji}(u)(z)).\label{EqThisIsTheEquationForAlphaIntheConCase}\\
&&\qquad\qquad\qquad\qquad\ =:\alpha_{ji}(u)(z,\hat z)^{-1}\in \U(1)\nonumber
\end{eqnarray} 
This establishes the existence of $\kappa^{\rm top}.$

It remains to  check that the constructed 
equivalence  class of the topological triple only depends on the  
equivalence class of the dynamical triple,
but this not difficult to establish.
We just mention that if we start with an exterior equivalent 
decker which differs locally by 
$c_i:U_i\times G\times G/N\to\U(\HH)$
from the cocycle $\mu_i$, then
$\kappa^{\rm top}$ differs by the 
null-homotopic bundle automorphism
which is locally of the form
 $c_i(u,-\sigma(\_+z),z)$. 
 
Finally, we remark that $\tau$ is natural.
Indeed, let $f:B'\to B$ be a map of bases spaces.
Then there is a commutative diagram
$$
\xymatrix{
{\rm Dyn}^\dag(B)\ar[r]^{\tau(B)}\ar[d]^{f^*}&{\rm Top}(B)\ar[d]^{f^*}\\
{\rm Dyn}^\dag(B')\ar[r]^{\tau(B')}&{\rm Top}(B'),
}
$$
because  in our local construction 
pullback with $f$ is just pullback of 
the underlying locally defined  functions by $f$  
which corresponds to the formal substitution of $u\in B$
by $f(u')$ for $u'\in B'$ in all formulas. 
\end{pf}

Our aim is to construct an inverse of 
$\tau(B):{\rm Dyn}^\dagger(B)\to {\rm Top}(B)$.
In particular, we wish to construct from the data of 
a topological triple  
$$
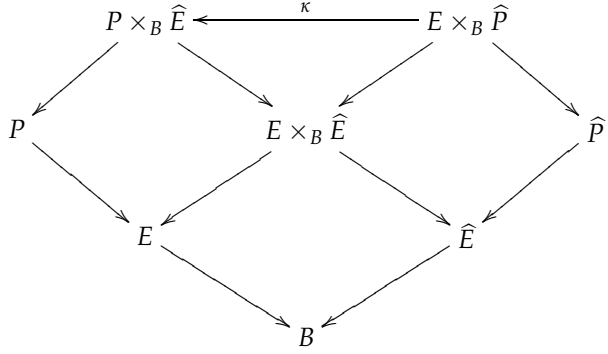
\xymatrix{
&P\times_B\widehat E\ar[dr]\ar[dl]&&E\times_B \widehat P\ar[dl]\ar[dr]\ar[ll]_\kappa&\\
P\ar[dr]&&E\times_B\widehat E\ar[dr]\ar[dl]&& \widehat P\ar[dl]\\
&E\ar[dr]&&\widehat E\ar[dl]&\\
&&B&&
}
$$
a decker $\rho$ on (a pair stably isomorphic to) $(P,E)$.
The space  $E\times_B \widehat P$ clearly has a $G/N$-action, and
by use of $\kappa$ we can define a $G/N$-action on $P\times_B\widehat E$.
The leading idea is that  the Poincaré condition enables us to 
push forward this action down to a  $G$-action on $P$.
We will see that this fails in general
as an obstruction against a decker occurs.
\\\\
Let $(\kappa,(P,E),(\widehat P,\widehat E))$ be a topological
triple over $B$, and
let $g_{ij},\zeta_{ij}$ and $\hat g_{ij},\hat\zeta_{ij}$ be transition
functions for the pairs $(P,E)$ and $(\widehat P,\widehat E)$.
Without restriction we let  $\widehat \HH=L^2(G/N)\otimes\HH$ be the underlying Hilbert space.
The existence of $\kappa$ implies the existence 
of locally defined $\kappa_i:U_i\to{\rm Map}( G/N\times\widehat G/N^\perp,\PU(\widehat\HH))$ such that 
$$
\kappa_j(u)(g_{ji}(u)+z,\hat g_{ji}(u)+\hat z)\  \hat\zeta_{ji}(u)(\hat z)\  \kappa_i(u)(z,\hat z)^{-1}= \zeta_{ji}(u)(z).
$$
Since $\kappa$ satisfies the Poincaré condition each $\kappa_i$ 
can be written in the  form
$$\kappa_i(u)(z,\hat z)=
\kappa^a_i(u)(z)\ \kappa^{\sigma}(z,\hat z)\ {\rm Ad}(v_i(u,z,\hat z))\ \kappa_i^b(u)(\hat z),
$$
wherein $\kappa^\sigma$ is from Lemma \ref{LemTheDefiOfThePoincareClass} and 
$\kappa^a_i, \kappa^b_i$ are some continuous projective unitary functions and 
$v_i:U_i\times G/N\times \widehat G/N^\perp\to\U(\widehat\HH)$ is a continuous unitary function. 
We introduce some short hands to get rid of $\kappa^a_i$ and $\kappa^b_i$.
We let $\zeta'_{ji}(u)(z):=\kappa_{j}^a(u)(g_{ji}(u)+z)^{-1}
\zeta_{ji}(u)(z)\kappa_i^a(u)(z)$ and analogously  \label{PageOfPrimedTransis}
$\hat\zeta'_{ji}(u)(\hat z):=\kappa_{j}^b(u)(g_{ji}(u)+z)
\hat\zeta_{ji}(u)(z)\kappa_i^b(u)(z)^{-1}$.
Thus we have $^\nou$
\begin{eqnarray}
\zeta'_{ji}(z)&=&\kappa^\sigma(g_{ji}+z,\hat g_{ji}+\hat z)
{\rm Ad}(v_j(g_{ji}+z,\hat g_{ji}+\hat z))\nonumber\\
&&\hat\zeta'_{ji}(\hat z){\rm Ad}(v_i(z,\hat z))^{-1}
\kappa^\sigma(z,\hat z)^{-1}.\label{EqATopTripIsATopTripIsATopTrip}
\end{eqnarray}
Let us choose, if necessary
after a refinement of the given atlas,
lifts $\overline\zeta'_{ji}:U_{ij}\to L^\infty(G/N,\U(\widehat\HH))$
and $\overline{\hat\zeta'}_{ji}:U_{ij}\to L^\infty(\widehat G/N^\perp,\U(\widehat\HH))$
of the transition functions (Lemma \ref{LemLiftForLocTrivFB}).
Then we obtain $^\nou$
\begin{eqnarray}
\overline\zeta'_{ji}(z)&=&
(\overline\kappa^\sigma(g_{ji}+z,\hat g_{ji}+\hat z)\otimes\Eins)\
v_j(u,g_{ji}+z,\hat g_{ji}+\hat z)\nonumber\\
&&\overline{\hat\zeta'}_{ji}(\hat z)\ v_i(z,\hat z)^{-1}\
(\overline\kappa^\sigma(z,\hat z)^{-1}\otimes\Eins)\ \alpha_{ji}(z,\hat z),
\label{EqTheDefiOfAlpha}
\end{eqnarray}
for some continuous $\alpha_{ji}:U_{ij}\to
L^\infty(G/N\times \widehat G/N^\perp,\U(1))$ and 
$\overline\kappa^\sigma(z,\hat z):=
\langle\hat\sigma(\hat z),\sigma(\_-z)-\sigma(\_)\rangle\in L^\infty( G/N,\U(1))
\subset\U(\widehat\HH)$;
$\sigma$ and $\hat\sigma$ are both chosen to be Borel. 
By direct computation, it follows that 
\begin{equation}\label{EqAlphaIstWieIso}
(\delta_{g\times\hat g}\alpha_{..})_{kji}(u)(z,\hat z)=
\psi_{kji}(u)(z)\hat\psi_{kji}(u)(\hat z)^{-1},
\end{equation}
for the twisted $\check C$ech cocycles $\psi_{...}\cdot\Eins:=
\delta_g\overline\zeta'_{..}$ and $\hat\psi_{...}\cdot \Eins:=
\delta_{\hat g}\overline{\hat\zeta'}_{..}$.
It is our task now to use the Poincar\'e condition of the isomorphism $\kappa$ 
to deduce a more concrete formula for $\alpha_{ji}$ from which we can extract
the existence of a decker. It will be sufficient to investigate the structure of 
$$
A_{ji}(u,z):={\rm Ad}\big(\alpha_{ji}(u)(z,..)\big)
\in \PP L^\infty(\widehat G/N^\perp,\U(1)).
$$
We let \label{PageOfAandBeta}
$$
\beta_{ji}(u,z):=
{\rm Ad}(\langle\hat\sigma(..+\hat g_{ji}(u))-\hat\sigma(..),\sigma( g_{ji}(u)+z)\rangle)
\in \PP L^\infty(\widehat G/N^\perp,\U(1))
$$
Note for the next lemma that $\beta_{ji}(u,z)$
equals 
$\hat\kappa^{\hat\sigma}(g_{ji}(u)+z,\hat g_{ji}(u))^{-1}$
up to conjugation by the unitary implemented 
${\rm Ad}(\lambda_{_{\widehat G/N^\perp}}(\hat g_{ji}(u)))$, so the $\check{\rm C}$ech classes $[\beta_{ji}]$ and 
$[(u,z)\mapsto \hat\kappa^{\hat\sigma}(g_{ji}(u)+z,\hat g_{ji}(u))^{-1}]$ 
 in $\check H^1(U_{ji}\times G/N,\underline{\U(1)})$ defined 
by these functions agree.
(Cp. Lemma \ref{LemThePoincClassForEverythindRep})

The Hilbert space $\HH$ occurs in the following to force 
that certain operators  which are defined on different tensor factors of  
$\HH\otimes L^2(\widehat G/N^\perp)\otimes\widehat\HH$ commute.

\begin{lem}\label{LemNatuerlichesMineralwasser}
$A_{ji}$  is continuous, and there are  continuous functions\footnote{
$\UA(\HH)$ is a commutative, contractible subgroup of 
$\U(\HH)$ such that $\U(1)\cdot\Eins\subset\UA(\HH)$.
See section \ref{SecTheUnitaryAndProjectiveUnitaryGroup},
page \pageref{PageOfUA}.}

$$
w_{ji}:U_{ji}\to{\rm Map}( G/N, L^\infty(\widehat G/N^\perp,\UA(\HH))))
\subset {\rm Map}( G/N, \U(\HH\otimes L^2( \widehat G/N^\perp)))
$$
and 
$$
\gamma_{ji}:U_{ji}\to \PP\UA(\HH)
$$ 
such that
$$
\Eins_\HH\otimes A_{ji}(u,z)=\Eins_\HH\otimes\beta_{ji}(u,z)\ {\rm Ad}(w_{ji}(u)(z))\ (\gamma_{ji}(u)\otimes\Eins_{L^2(\widehat G/N^\perp)}).
$$
\end{lem}
\begin{pf}
 Equation (\ref{EqTheDefiOfAlpha}) implies
\begin{eqnarray*}
\Eins\otimes\zeta'_{ji}(u)(z)&=&
{\rm Ad}\big( \Eins_{L^2(\widehat G/N^\perp)}\otimes\overline \zeta_{ji}(u)(z)\big)\\
&=&{\rm Ad}\Big((\overline\kappa^\sigma(g_{ji}(u)+z,\hat g_{ji}(u)+..)\otimes\Eins_\HH)
v_j(u,g_{ji}(u)+z,\hat g_{ji}(u)+..)\\
&&\overline{\hat\zeta'}_{ji}(u)(..)v_i(u,z,..))^{-1}
(\overline\kappa^\sigma(z,..)^{-1}\otimes\Eins_\HH)\Big)\ (A_{ji}(u)(z)\otimes\Eins_{\widehat\HH})\\
&&\in \PU(L^2(\widehat G/N^\perp)\otimes\widehat\HH)).
\end{eqnarray*}
The terms  inside the bracket are continuous unitary functions on
$U_{ji}\times G/N$.
Thus we have equality of the  $\check{\rm C}$ech classes 
$[\zeta'_{ji}]=[\Eins\otimes\zeta'_{ji}]=[A_{ji}\otimes\Eins]=[A_{ji}]\in \check H^1(U_{ji}\times G/N,\underline{\U(1)})$.
We compute the $\check{\rm C}$ech class of $\Eins\otimes\zeta'_{ji}$ by equation
(\ref{EqATopTripIsATopTripIsATopTrip}) for  $\hat z=0$. 
Then $\kappa^\sigma(z,0)=\Eins$, hence has trivial
$\check{\rm C}$ech class.
The map $(u,z)\mapsto \hat\zeta'_{ji}(u,0)$  has   $\check{\rm C}$ech class
$[(u,z)\mapsto \hat\zeta'_{ji}(u,0)]$. By Remark \ref{RemAbelianU}, 
let $ \gamma_{ji}:U_{ji}\to\PP \UA(\HH)$ be such that 
$[(u,z)\mapsto \hat\zeta'_{ji}(u,0)]=[(u,z)\mapsto \gamma_{ji}(u)]
={\rm pr}_{U_{ji}}^*[\gamma_{ji}]$.
Then
\begin{eqnarray*}
[\Eins\otimes \zeta'_{ji}]&=&[(u,z)\mapsto\kappa^\sigma(g_{ji}(u)+z,\hat g_{ji}(u))]
+[u\mapsto \hat\zeta'_{ji}(u,0)]\\
&=&\big(
(g_{ji}\circ {\rm pr}_{U_{ji}}+{\rm pr}_{G/N})\times (\hat g_{ji}\circ{\rm pr}_{U_{ji}})\big)^*\pi
+{\rm pr}_{U_{ji}}^*[\gamma_{ji}]\\
&=&[\beta_{ji}]+{\rm pr}_{U_{ji}}^*[\gamma_{ji}]. 
\end{eqnarray*}
This implies that $A_{ji}\otimes\Eins_\HH$ equals $\beta_{ji}\otimes\gamma_{ji}$ 
up to ${\rm Ad}$ of  a continuous unitary function $w_{ji}:U_{ji}\times G/N\to \U(L^2(\widehat G/N^\perp,\HH))$.
$w_{ji}$ takes  values in the subgroup 
$L^\infty(\widehat G/N^\perp,\UA(\widehat\HH))$ only, since
$A_{ji},\beta_{ji}$ and $\gamma_{ji}$ are maps  into the subgroup
$\PP L^\infty(\widehat G/N^\perp,\UA(\widehat\HH)) 
\subset \PU(L^2(\widehat G/N^\perp,\widehat\HH))$.

Finally, we just interchange the order of the Hilbert spaces 
in the tensor product $L^2(\widehat G/N^\perp)\otimes\HH$
to $\HH\otimes L^2(\widehat G/N^\perp)$. 
This will be convenient later.
\end{pf}

\begin{rem}
By passing to a refined atlas, we can always
achieve that  $\gamma_{ji}=\Eins$, for 
we can apply Lemma \ref{LemLiftForLocTrivFB} to
the functions $U_{ji}\ni u\mapsto\hat\zeta'_{ji}(u)(0)$.
In that case we have $w_{ji}(u)(z)(\hat z)\in \U(1)$.
\end{rem}

It is important to note that although $w_{ji}$ depends 
on the choice of $\gamma_{ji}$ the term
$(d_*w_{..})_{ji}(u)=d(w_{ji}(u))$ does not;
here $d$ is the (first)
boundary operator of group cohomology
\begin{eqnarray*}
{\rm Map}(G/N, L^\infty(\widehat G/N^\perp,\UA(\HH)))
&\stackrel{d}{\longrightarrow} &
{\rm Map}(G,{\rm Map}(G/N, L^\infty(\widehat G/N^\perp,\UA(\HH))))\\
w&\mapsto &\big( (h,z)\mapsto w(z+hN)(..)w(z)(..)^{-1}\big)
\end{eqnarray*}
The sheaf of continuous functions into the range (or source) of $d$ 
is a $\underline{G/N\times \widehat G/N^\perp}$-module in the obvious 
way.
So the following statement is meaningful 
for the twisted $\check{\rm C}$ech boundary
operator $\delta_{g_{..}\times \hat g_{..}}$. 

\begin{lem}
 ${\rm Ad}\Big((\delta_{g_{..}\times{\hat g_{..}}} d_*w_{..})_{kji}(u)(h,z)(..)\Big)=\Eins$, 
for all $u\in U_{kji},z\in G/N, h\in G.$
\end{lem}
\begin{pf}
Since $d_*$ and $\delta_{g_{..}\times{\hat g_{..}}} $ commute, we have to show that
the expression ${\rm Ad}((\delta_{g_{..}\times{\hat g_{..}}} w_{..})_{kji}(u)(z)(..))$ is independent of $z$.
Because  ${\rm Ad}((\delta_{g\times\hat g}\alpha_{..})_{kji}(u)(z,..))=
{\rm Ad}(\hat\psi_{kji}(u)(..)^{-1})$ and $\gamma_{ji}(u)$
are independent of $z$ we have to compute $\delta_{g_{..}\times{\hat g_{..}}}\beta_{..}$.
This yields $^\nou$
\begin{eqnarray*}
&&(\delta_{g_{..}\times{\hat g_{..}}} \beta_{..})_{kji}(z)\\
&=&{\rm Ad}\big(
\langle \hat\sigma(..+\hat g_{ji})-\hat\sigma(..),-\sigma(g_{ji}+z)\rangle\\
&&\langle \hat\sigma(..+\hat g_{ki})-\hat\sigma(..),\sigma(g_{ki}+z)\rangle\\
&&\langle \hat\sigma(..+\hat g_{ji}+\hat g_{kj})-\hat\sigma(..+\hat g_{ji}),-\sigma(g_{ki}+z)\rangle\big)\\
&=& {\rm Ad}(
\langle \hat\sigma(..+\hat g_{ji})-\hat\sigma(..),\sigma(g_{ki})-\sigma(g_{ji})\rangle),
\end{eqnarray*}
and the lemma is proven.
\end{pf}

It follows that
\begin{equation}\label{EqThisIsVarphiTheObstructionCocycleAgainstDeckers}
\varphi_{kji}(u)(h,z):=(\delta_{g_{..}\times{\hat g_{..}}} d_*w_{..})_{kji}(u)(h,z)(..)\in\U(1)
\end{equation}
satisfies $\delta_g\varphi=1$ and $d_*\varphi_{ji}=1$, i.e. it
defines a twisted $\check{\rm C}$ech 2-class
\begin{equation}
[\varphi_{...}]\in\check H^2(B,\underline{Z^1_{\rm cont}(G,{\rm Map}(G/N,\U(1)))},g_{..}).
\end{equation}

\begin{lem}\label{LemTheObstructionLemmaAgainstADecker}
The construction of the class $[\varphi_{...}]$ defines a natural map
$${\rm Top}(E,B)\to \check H^2(B,\underline{Z^1_{\rm cont}(G,{\rm Map}(G/N,\U(1)))},g_{..}).
$$
\end{lem}

\begin{pf}
We must show that the constructed class only depends on the 
isomorphism class of the triple. 
Let us check all choices in reversed order of their appearance.

Firstly,
the choice of $w_{ji}$ is only unique up to 
a continuous scalar function $U_{ji} \times G/N\to \U(1)$,
so $\varphi_{...}=\delta_{g\times\hat g}d_*w_{..}$  changes 
by a boundary.

Secondly,
the choice of $\gamma_{ji}$, is only unique up to 
a unitary implemented function
$U_{ji}\to \UA(\HH)$, so $w_{ji}$ may change by this function.
But this function is independent 
of $z\in G/N$,
 so $d_*w_{ji}$ does not change.

Thirdly,
if we choose  different lifts of the transition functions $\zeta_{ji}$, 
this would change $\alpha_{ji}$ by a function
$U_{ji}\to L^\infty(G/N,\U(1))$, so $A_{ji}$ does not change.

Fourthly,
if we choose different lifts of the transition functions 
$\hat\zeta_{ji}$, then 
$\alpha_{ji}$ is changed by a function $U_{ji}\to L^\infty (\widehat G/N^\perp,\U(1))$,
so $A_{ji}$, thus $w_{ji}$,
changes by a function independent of $z\in G/N$,
 so $d_*w_{ji}$ does not change.

Fifthly, 
if we choose another atlas  for our construction,
we can take a common refinement and the normalisation 
procedure $\zeta_{ji}\mapsto\zeta'_{ji},\hat\zeta_{ji}\mapsto\hat\zeta_{ji}'$,
leads to the same equations as above, except
we restricted us to the refined atlas, so the class of
$\varphi_{...}$ is not changed.

Sixthly,
if we start with a topological triple  which is
isomorphic to $(\kappa,(P,E),(\widehat P,\widehat E)$,
then, locally, the isomorphisms of the underlying pairs 
have the same effect as a change of the atlas which
does not change the class of $\varphi_{...}$.
However, we must take care of the homotopy commutativity 
of diagram (\ref{DiagMorphOfTopTriples}).
So
if  $\kappa'$ is homotopic to $\kappa$ , then 
they differ locally by a continuous unitary 
function $v'_i:U_i\times G/N\times\widehat G/N^\perp\to 
\U(\widehat \HH)$, and
equation  (\ref{EqTheDefiOfAlpha}) becomes 
\begin{eqnarray*}
&&v'_j(u,g_{ji}(u)+z,\hat g_{ji}(u)+\hat z)^{-1}
\overline\zeta'_{ji}(u)(z)
v'_i(u,z,\hat z)\\
&=&
(\overline\kappa^\sigma(g_{ji}(u)+z,\hat g_{ji}(u)+\hat z)\otimes\Eins)\
v_j(u,g_{ji}(u)+z,\hat g_{ji}(u)+\hat z)\nonumber\\
&&\overline{\hat\zeta'}_{ji}(u)(\hat z)\ v_i(u,z,\hat z)^{-1}\
(\overline\kappa^\sigma(z,\hat z)^{-1}\otimes\Eins)\ \alpha'_{ji}(u)(z,\hat z).
\end{eqnarray*}
We must investigate how $\alpha'_{ji}$ is related to $\alpha_{ji}$.
The family $v'_i$ defines the bundle automorphism
$\kappa'\circ\kappa^{-1}$ on $P\times_B\widehat E$ so 
\begin{equation}\label{Eqsjhdkfhljkdfh}
\zeta'_{ji}(u)(z)\big(v'_i(u,g_{ji}(u)+z,\hat g_{ji}(u)+\hat z)\big)
=v'_j(u,z,\hat z)\ \alpha''_{ji}(u)(z,\hat z),
\end{equation}
for a scalar $\alpha''_{ji}:U_i\times G/N\times\widehat G/N^\perp\to\U(1)$. 
It follows from the three equations 
for $\alpha_{ji},\alpha'_{ji}$ and $\alpha''_{ji}$ that
$\alpha'_{ji}=\alpha_{ji}\ \alpha''_{ji}$.
So $w_{ji}$ changes by $\alpha''_{ji}$.
But $\alpha''_{..}$ is a cocycle $\delta_{g\times \hat g}\alpha''_{..}=1$
which is easily computed by its definition 
(\ref{Eqsjhdkfhljkdfh}),
so $\delta_{g\times\hat g}w_{..}$ and $\varphi_{...}$ do not
change.

Finally,
we just remark that the defined map is natural with respect to pullback.
I.e. if $f:B'\to B$ is a map of base spaces, then
there is a commutative diagram
$$
\xymatrix{
{\rm Top}(E,B)\ar[r]\ar[d]^{f^*}&
\check H^2(B,Z^1_{\rm cont}(G,{\rm Map}(G/N,\U(1))),g_{..})\ar[d]^{f^*}\\
{\rm Top}(f^*E,B')\ar[r]&
\check H^2(B',Z^1_{\rm cont}(G,{\rm Map}(G/N,\U(1))),f^*g_{..}).
}
$$

This proves the lemma.
\end{pf}

The quotient map $G\to G/N$ induces a map on
the twisted $\check{\rm C}$ech groups
$$
\xymatrix{
\check H^2(B,\underline{Z^1_{\rm cont}(G/N,{\rm Map}(G/N,\U(1)))},g_{..})
\ar[d]\\
\check H^2(B,\underline{Z^1_{\rm cont}(G,{\rm Map}(G/N,\U(1)))},g_{..}),
}
$$
and the map of Lemma
\ref{LemTheObstructionLemmaAgainstADecker}
has an obvious factorisation 
$$
\xymatrix{
{\rm Top}(E,B)\ar[drr]_{
\rm Lemma\
\ref{LemTheObstructionLemmaAgainstADecker}\qquad\qquad}
\ar[rr]
&&\check H^2(B,\underline{Z^1_{\rm cont}(G/N ,{\rm Map}(G/N,\U(1)))},g_{..})
\ar[d]\\
&& \check H^2(B,\underline{Z^1_{\rm cont}(G,{\rm Map}(G/N,\U(1)))},g_{..}),
}
$$
since the definition (\ref{EqThisIsVarphiTheObstructionCocycleAgainstDeckers}) 
 of $\varphi_{...}$ implies
\begin{equation}\label{EqTheStricnessObstructionFactors}
\varphi_{kji}(u)(h+n,z)=\varphi_{kji}(u)(h,z),
\end{equation} 
for all $n\in N$.

\begin{defi}\label{DefiOfStrictAndAlmoStrict}
A topological T-duality triple over $B$ is called {\bf almost strict} 
if its $\check{\rm C}$ech class $[\varphi_{...}]$ defined by 
Lemma 
\ref{LemTheObstructionLemmaAgainstADecker} vanishes.
The triple is called {\bf strict} if its $\check{\rm C}$ech class 
defined by the  horizontal map in the diagram above vanishes already.
\end{defi}

We denote by ${\rm Top}^{\rm as}(B)$
respectively   ${\rm Top}^{\rm s}(B)$ the  set of all
almost strict respectively strict topological triples over $B$;
so we have obvious inclusions \label{PageOfAllTheTops}
$$
{\rm Top}^{\rm s}(B)\subset{\rm Top}^{\rm as}(B)\subset{\rm Top}(B).
$$

\begin{rem}
For the class 
of a topological triple 
$[(\kappa,(P,E),(\widehat P,\widehat E))]\in{\rm Top}(B)$
its class in ${\rm Top}(E,B)$ is 
only well-defined up to the action of
${\rm Aut}_B(E)$ on ${\rm Top}(E,B)$.
However, the the vanishing of the 
obstruction class in Definition \ref{DefiOfStrictAndAlmoStrict} is independent 
of the possible choices, so 
${\rm Top}^{\rm as}(B)$  and  ${\rm Top}^{\rm s}(B)$
are well-defined.
\end{rem}

We will see   
that strict and almost strict play a major rôle in our 
theory. 
The following two lemmata give a first feeling. 

\begin{lem}\label{LemImTauSubsetStrict}
The image of the map $\tau(B)$ from dualisable dynamical
to topological triples is contained 
in the set of strict topological triples,
$$
{\rm im}(\tau(B))\subset {\rm Top}^{\rm s}(B).
$$
\end{lem}

\begin{pf}
In equation (\ref{EqThisIsTheEquationForAlphaIntheConCase}) 
we already computed $\alpha_{ji}$ for a topological triple which is
constructed out of a dynamical one. 
(In (\ref{EqThisIsTheEquationForAlphaIntheConCase}) we did not 
normalise $\zeta_{ji}\mapsto \zeta_{ji}'$, but this does not change 
$\alpha_{ji}$.) The result is
 $\alpha_{ji}(u)(z,\hat z)=
\langle\hat\sigma(\hat g_{ji}(u)+\hat z)-\hat\sigma(\hat z),
\sigma(z)\rangle\alpha'_{ji}(u)(z)$.
So in this case we have
\begin{eqnarray*}
A_{ji}(u)(z)&=&
{\rm Ad}(\langle\hat\sigma(\hat g_{ji}(u)+..)-\hat\sigma(..),\sigma(z)\rangle)\\
&=&
{\rm Ad}(\langle\hat\sigma(\hat g_{ji}(u)+..)-\hat\sigma(..),
\sigma(g_{ji}(u)+z)-\sigma(g_{ji}(u)\rangle)\\
&=&\beta_{ji}(u,z)\ \underbrace{{\rm Ad}(\langle\hat\sigma(\hat g_{ji}(u)+..)-\hat\sigma(..),
-\sigma(g_{ji}(u)\rangle)},\\
&&\qquad\qquad\qquad\qquad
\textrm{independent of }z\in G/N
\end{eqnarray*}
i.e. we can choose $w_{ji}$ such that $w_{ji}(u)(z)$ does not depend on $z$.
Of course, the choice of $w_{ji}$ is only determined up to
a scalar $U_{ji}\times G/N\to \U(1)$,  but in any case 
$\chi_{ji}\cdot \Eins:=d_*w_{ji}$ defines a function
$U_{ji}\to Z^1_{\rm cont}(G,{\rm Map}(G/N,\U(1)))$, so 
by construction $\varphi_{kji}=(\delta_g\chi_{..})_{kji}$ is a boundary.
\end{pf}

The next (technical) lemma
will  be  the crucial point in the construction of 
a decker from the data of a topological triple.

\begin{lem}\label{LemTheCrucialLemmaToConstructDelta}
Assume $(\kappa,(P,E),(\widehat P,\widehat E))$ is almost
strict. Then we can find a (sufficiently refined) atlas $\{U_i|i\in I\}$
such that for  $w_{ji}$  from above
there exists a family $m_i:U_i\to Z^1_{\rm cont}(G,{\rm Map}(G/N, 
L^\infty(\widehat G/N^\perp,\UA(\HH) )))
$ such that
\begin{eqnarray}\label{EqInClaim1Forbbbbr}
{\rm Ad}(d_* w_{ji}(u)(h,z)(..))&=&{\rm Ad}((\delta_ {g_{..}\times{\hat g_{..}}} m_.)_{ji}(u)(h,z)(..))\\
&&\in\PP L^\infty(\widehat G/N^\perp,\UA(\HH))\nonumber
\end{eqnarray}
\end{lem}

\begin{pf}
The proof uses a standard Zorn's lemma argument
as it can be found in \cite{Di}.

 First we note that the space
$Z:=Z^1_{\rm cont}(G,{\rm Map}(G/N, L^\infty(\widehat G/N^\perp,\UA(\HH) )))$ is contractible. In fact, if $H:[0,1]\times\UA(\HH)\to\UA(\HH)$ 
is a contraction  with each $H(t,\_)$  a group homomorphism, 
then the push-forward $H_\#:[0,1]\times Z\to Z$ preserves the
cocycle relation and is a contraction.

We can assume that the atlas is sufficiently refined such that,
firstly,
$\varphi_{...}$ from (\ref{EqThisIsVarphiTheObstructionCocycleAgainstDeckers})
 is a boundary, i.e.
there exist $\chi_{ji}:U_{ji}\to Z^1_{\rm cont}(G,{\rm Map}(G/N,\U(1)))$ 
such that $(\delta_g\chi)_{kji}=\varphi_{kji}$, and secondly,
as $B$ is a paracompact Hausdorff space, 
even the closed cover $\bigcup_{i\in I}\overline U_i=B$ is locally finite and 
that $w_{ji},\chi_{ji}$ are well-defined on $\overline U_{ji}$, for all $j,i\in I$.
We let
\begin{eqnarray}
M&:=&\Big\{ (J, m_.)\ |\ J\subset I, {\rm for\ all\ }j\in J:
 m_j:\overline U_j\to Z
{\rm\ such \ that\ }\nonumber\\
&&{\rm for\ all\ }i,j\in J { \rm\ and\ all\ }
u\in U_{ji},z\in G/N{\rm \ we\ have}\nonumber\\
&&  m_i(u)(h,z)(..)=
d_*w_{ji}(u)(h)(z)(..)^{-1} \chi_{ji}(u)(h,z)\nonumber \\
&&\qquad\qquad\qquad\qquad m_j(u)(h,g_{ji}(u)+z)(\hat g_{ji}(u)+..)
\Big\}\label{EqTheDefiOfMi}
\end{eqnarray}
For each $i\in I$, $(\{i\},\{\Eins\})\in M$, because 
we can assume that $w_{ii}=\Eins,\chi_{ii}=1$, so $M$ is non-empty.
We define a partial order on $M$ by
$(J, m_.)\le(J', m_.)$ if and only if 
$J\subset J'\subset I$ and $ m_j= m'_j,$
for all $j\in J$. 
For each chain in $M$ the union of the index and cocycle sets 
is  an upper bound, hence by Zorn's lemma there exists a
maximal element $(J, m_.)$.
Assume $J\not= I$, so there is some $a\in I\backslash J$.
Let $R:=\bigcup_{j\in J}(\overline U_j\cap \overline U_{a})\subset \overline U_{a}$.
For $u\in R$ we define
$\tilde m_{a}(u)(h,z)(..):=
d_*w_{ja}(u)(h,z)(..)^{-1}
 m_{j}(u)(h,g_{ja}(u)+z)(\hat g_{ja}(u)+..)
\chi_{ja}(u)(h,z)$, if 
$u\in \overline U_j$.
Due to $\delta_g\chi=\delta_{g\times\hat g}d_*w$, 
this definition is independent of $j\in J$.
We end up with a diagram  
$$
\xymatrix{
R\ar[r]^{\tilde m_{a}}\ar[d]_\cap&Z\\
\overline U_{a}\ar@{.>}[ru]_{ m_{a}}&.
}
$$
Since our  cover  is locally finite, $R$ is closed, but 
$Z$ is contractible, therefore
an extension $m_{a}$ exists \cite[Lem. 4]{DD}.
This contradicts the maximality of $(J, m_.)$,
so $J=I$. Finally, equation (\ref{EqInClaim1Forbbbbr})
holds, since $\chi_{ji}(u)(h,z)\in\U(1)$. 
\end{pf}

The constructed family $\{m_i\}$ is the last ingredient 
to write down an explicit formula for a decker $\rho^{\rm dyn}$ on
(a stabilisation of) the pair $(P,E)$.
It is not hard to guess that the cocycles $m_i$ will be an
essential part of the cocycles $\mu_i^{\rm dyn}$ which we
will define to implement the decker $\rho^{\rm dyn}$ locally.
But unfortunately 
the constructed family $\{m_i\}$ is by no means unique,
and this non-uniqueness is the origin of the following discussion. 
As a matter of fact this discussion will simplify drastically when
we consider the special case of $G=\R^n$ with lattice $N=\Z^n$
in the next section below.
However, now we continue with the discussion of almost strict triples from above
and work out the general framework.
\\\\
Let $\varphi_{...}$ be the twisted $\check{\rm C}$ech cocycle 
from equation (\ref{EqThisIsVarphiTheObstructionCocycleAgainstDeckers}).
The triple under consideration is assumed to be almost strict, so
(after refining the atlas $U_\bullet$) we have
$\varphi_{kji}=(\delta_g\chi_{..})_{kji}$ 
for a chain $\chi_{..}\in \check C^1(U_{\bullet},\underline{ Z^1_{\rm cont}(G,{\rm Map}(G/N,\U(1))},g_{..})$ 
as in the proof above.
Obviously, this chain $\chi_{..}$ is only well defined up to a cocycle 
$\chi^1_{..}\in\check Z^1(U_{\bullet},\underline{ Z^1_{\rm cont}(G,{\rm Map}(G/N,\U(1))},g_{..}).$

A choice of $\chi_{..}$ determines the family 
$\{m_i\}$ in  equation (\ref{EqTheDefiOfMi}) not completely 
but up to
a family $n_i:U_i\to Z^1_{\rm cont}(G,{\rm Map}(G/N,L^\infty(\widehat G/N^\perp,\UA(\HH))))$ which satisfies $(\delta_{g\times\hat g}n_.)_{ji}=\Eins.$
We already mentioned that $m_i$ will be  part of 
the cocycles $\mu_i^{\rm dyn}$. 
Then it will turn out that the two families 
 $\{m_i\}$ and $\{m_i\ n_i\}$  define exterior equivalent 
deckers.
However, this is not the case when make we another choice of
of $\chi_{..}$, say  $\chi_{..}\ \chi^1_{..}$ for $\chi^1_{..}$ as above.
Then $\{m_i\}$ must be replaced by $\{m_i\ m_i^1\}$ 
for $m_i^1$ being such that 
$(\delta_{g\times \hat g}m^1_.)_{ji}=\chi^1_{ji}\cdot\Eins,
$
and we will find that
the corresponding deckers are exterior equivalent if and only 
if the class 
$[\chi^1_{..}]\in
\check H^1(B,\underline{Z^1_{\rm cont}(\!G\!/\!N,{\rm Map}\!(\!G\!/\!N,\U(\!1\!)\!)\!)},g_{..})$
vanishes. In other words, for each class $[\chi^1_{..}]$
we obtain a different class of dynamical triples.
It is then the obvious question, whether  the different 
dynamical triples still have something in common.
Or if $x$ is the class of almost strict topological triple
and $x^{\rm dyn}$ is the class of one of the possible dynamical triples indicated,
is there  any relation between $x$ and $\tau(B)(x^{\rm dyn})$?
Can we describe the difference, in particuar, when do they equal?
A partial answer of this question is already given in Lemma \ref{LemImTauSubsetStrict} which shows that $x$ and $\tau(B)(x^{\rm dyn})$
only can equal if $x$ is strict.

Consider the short exact sequences 
$$
\xymatrix{
0\ar[d]&&1\ar[d]\\
N^\perp\ar[d]\ar[rr]&&
Z^1_{\rm cont}(\!G\!/\!N,{\rm Map}\!(\!G\!/\!N,\U(\!1\!)\!)\!)
\ar[d]\\
\widehat G\ar[d]\ar[rr]^
{\hspace{-1cm}\subset} &&
Z^1_{\rm cont}(\!G,{\rm Map}\!(\!G\!/\!N,\U(\!1\!)\!)\!)
\ar[d]^{\rm restr.}\\
\widehat G/N^\perp\ar[d]\ar[rr]^\cong&&
\widehat N\ar[d]\\
0&&0.
}
$$
They both induce long exact sequences in (twisted)
$\check{\rm C}$ech cohomology. The relevant part for us is 
\begin{equation}\label{DiagOfExatSequences}
\xymatrix{
&
Q_{_{G/N}}:=\check H^1(B,\underline{Z^1_{\rm cont}(\!G\!/\!N,{\rm Map}\!(\!G\!/\!N,\U(\!1\!)\!)\!)},g_{..})
\ar[d]^q\\
\check H^1(B,\underline{\widehat G})\ar[d]\ar[r] &
Q_{_{G}}:=\check H^1(B,\underline{Z^1_{\rm cont}(\!G,{\rm Map}\!(\!G\!/\!N,\U(\!1\!)\!)\!)},g_{..})
\ar[d]\\
\check H^1(B,\underline{\widehat G/N^\perp})\ar[d]\ar[r]^\cong&
\check H^1(B,\underline{\widehat N},g_{..})\\
\check H^2(B,\underline{N^\perp}).
&
}
\end{equation}
Due to equation (\ref{EqTheStricnessObstructionFactors})
restriction $\chi_{ji}(u)|_N$ defines a class $[\chi_{..}|_N]\in \check H^1(B,\underline{\widehat N}, g_{..})$ and by diagram
(\ref{DiagOfExatSequences}) a class in $\check H^1(B,\underline{\widehat G/N^\perp})$, i.e the class of a $\widehat G/N^\perp$-principal fibre bundle 
$\widehat E_{\chi_{..}}\to B$.
Exactness of the columns in diagram (\ref{DiagOfExatSequences})
implies that this class is only well-defined up
to the quotient group
$Q_{_{G}}/ q(Q_{_{G/N}}),$
because the class  varies with the choice of $\chi_{..}$.
The bundle $\widehat E_{\chi_{..}}$ will be connected to 
the description of $\tau(B)(x^{\rm dyn})$. Namely, if
$x=[(\kappa,(P,E),(\widehat P,\widehat E))]$ and 
$x^{\rm dyn}=[(\rho^{\rm dyn},P_{\rm dyn},E)]$ ($P_{\rm dyn}$ will
be a stabilisation of $P$ with a certain Hilbert space), then
the result will be
$$
\tau(B)(x^{\rm dyn})=[(\kappa^{\rm top},{(P_{\rm dyn}}_{\rm top},E),
(\widehat {P_{\rm dyn}},\widehat E\oplus \widehat E_{\chi_{..}}))].
$$
If the topological triple is strict,
then by definition we can choose $\chi_{..}$ such that 
$ \chi_{ji}(u)|_N=1$, so
the class $[\chi_{..}|_N]$ vanishes and the bundle $\widehat E_{\chi_{..}}$ is trivialisable.
In that case 
we give a description of
$\widehat {P_{\rm dyn}}$ 
in the proof of the theorem below.
If the topological triple is in the image 
of $\tau(B)$, then we even know more
as the next theorem makes  precise.
Therein $\mathcal P(M)$ denotes the power set of a set $M$, 
and we denote the image of $\tau(B)$ by \label{PageOfTopImAndP}
${\rm Top}^{\rm im}(B)$, so
$$
{\rm Top}^{\rm im}(B)\subset
{\rm Top}^{\rm s}(B)\subset
{\rm Top}^{\rm as}(B)\subset
{\rm Top}(B).
$$

\begin{thm}\label{ThmTheMainTheoremOfThisThesisNumberThree}
The map $\tau(B)$ is injective, and there are three maps
\begin{eqnarray*}
&\delta^{\rm as}(B)&:{\rm Top}^{\rm as}(B)\to \mathcal P({\rm Dyn}^\dag(B))\\
&\delta^{\rm s}(B)&:{\rm Top}^{\rm s}(B)\to \mathcal P({\rm Dyn}^\dag(B))\\
&\delta^{\rm im}(B)&:{\rm Top}^{\rm im}(B)\to {\rm Dyn}^\dag(B)
\end{eqnarray*}
with the following  properties:
\begin{itemize}
\item[(a)]
For each $x\in{\rm Top^{as}}(B)$ the set $\delta^{\rm as}(B)(x)\subset
{\rm Dyn}^\dag(B)$ is a $Q_{_G}$-torsor. 
\item[(b)]
If $x\in{\rm Top^{s}}(B)$, then
 $\delta^{\rm s}(B)(x)\subset \delta^{\rm as}(B)(x)$ is a 
 $q(Q_{_{G/N}})$-subtorsor, and
 for each $x^{\rm dyn}\in \delta^{\rm s}(B)(x)$
the (class of the) $\widehat G/N^\perp$-bundle
$\widehat{E_{\rm dyn}}$
 of the (class of the) dual pair of $\tau(B)(x^{\rm dyn})$
is (the class of) the $\widehat G/N^\perp$-bundle $\widehat E$ of $x$.

\item[(c)]
If $x$ is in the image of $\tau(B)$, then 
$\delta^{\rm im}(B)(x)\in\delta^{\rm s}(B)(x)$, and
$\delta^{\rm im}(B)$ is the inverse of $\tau(B)$. 
\end{itemize}
Moreover,  $\delta^?(B)$ is natural in the sense that it extends to a natural transformation
of functors
$$
\delta^?: {\rm Top^{?}}\to \mathcal P\circ{\rm Dyn}^\dag,\qquad ?={\rm as,\ s,\ im}.
$$
\end{thm}

\begin{pf}
The proof consists of several steps which are rather technical, so 
we first give an overview of the proof:
\\\\
In Step 1 we construct from the data of an almost 
strict topological triple a decker on the underlying pair 
(after stabilisation).
The local construction will depend on many choices
and it is the statement of Step 2 that almost all of these 
choices do not interfere with the equivalence class of triple defined.
An exception is the choice of $\chi_{..}$ (above)
which will cause the $Q_{_G}$-torsor and
the $q(Q_{_{G/N}})$-subtorsor structure in $(a),(b)$.
Then in Step 3  we  compute the composition
$\delta^{\rm im}(B)\circ\tau(B)$ and find
that this is the identity on ${\rm Dyn}^\dag(B)$,
hence $\tau(B)$ is injective and $\delta^{\rm im}(B)$ is surjective.
In Step 4 we compute the reverse
composition $\tau(B)\circ\delta^{\rm s}(B)$, 
which will lead us to the result stated in $(b)$.
In particular this calculation shows that
$\tau(B)\circ\delta^{\rm im}(B)$
is the identity on the image of $\tau(B)$.
As both compositions $\delta^{\rm im}(B)\circ\tau(B)$ and 
$\tau(B)\circ\delta^{\rm im}$ are the identity maps the
statement of $(c)$ is clear then.
In Step 5 we finally  comment on the naturality of the maps. 
\\\\
In the hole of the proof we maintain the notation introduced above.
\\\\
{\bf  Step 1:} \label{PageOfStep1}
For an almost strict topological triple $(\kappa,(P,E),(\widehat P,\widehat E))$ we construct a dualisable dynamical triple $(\rho,P_{\rm dyn},E_{\rm dyn})$.
\\\\
The idea is simple. All we have to do is to write down an explicit formula for
a family of cocycles which satisfies equation (\ref{EqLocStucOfDeckers})
for the transition functions of the pair $(P,E)$.

Assume the topological triple has underlying Hilbert space $\widehat \HH=L^2(G/N)\otimes \HH$.
Let $\zeta_{ji},\zeta'_{ji},v_i, w_{ji},\varphi_{kji},\dots$ be as above. We can assume that the atlas
$U_\bullet$ is sufficiently refined such that
$\varphi_{...}$ is a boundary. 
So we choose $\chi_{..}$ as above such that
$\varphi_{...}=\delta_g\chi_{..}$,
then let $m_i$ be defined by (\ref{EqTheDefiOfMi}) in the last lemma.

Let us consider 
$
U_{ji}\ni u\mapsto \Eins_\HH\otimes \lambda_{_{\widehat G/N^\perp}}(\hat g_{ji}(u))
\in\U(\HH\otimes L^2(\widehat G/N^\perp)),$
this defines a unitary implemented $\check{\rm C}$ech cocycle, therefore
it is a boundary, i.e. there are $l^0_i:U_i\to\U( \HH\otimes L^2(\widehat G/N^\perp))$ such
that $l^0_j(u)^{-1}\ l^0_i(u)=\Eins_\HH\otimes \lambda_{_{\widehat G/N^\perp}}(\hat g_{ji}(u))$.
Note that this need not be true for $\lambda_{_{\widehat G/N^\perp}}(\hat g_{ji}(u))$ alone
as $\widehat G/N^\perp$ may be finite, $L^2(\widehat G/N^\perp)$ finite 
dimensional and $\U(L^2(\widehat G/N^\perp))$ not contractible.
We conclude that there exists a family \label{PageOfLambdaili}
$\lambda_i:U_i\to{\rm Map}(G/N, \PU(\HH\otimes L^2(\widehat G/N^\perp)\otimes
L^2(G/N)\otimes \HH))$ such
that this family permits lifts
$l_i:U_i\to {\rm Map} (G/N, \U(\HH\otimes L^2(\widehat G/N^\perp)\otimes
L^2(G/N)\otimes \HH))$ such
that 
\begin{eqnarray}\label{EqTheEqForOverlineLambda}
&&l_j(u)(g_{ji}(u)+z)^{-1} (\Eins\otimes\Eins\otimes\overline\zeta_{ji}(u)(z)) 
l_i(u)(z)\nonumber\\
&=&
\Eins\otimes \lambda_{_{\widehat G/N^\perp}}(\hat g_{ji}(u)) \otimes\overline\zeta_{ji}(u)(z).
\end{eqnarray}
For example we may take $\lambda_i(u)(z)={\rm Ad}(l^0_i(u)\otimes\Eins\otimes\Eins)$, but 
it will be important that we allow ourselves to have the freedom of 
a more flexible  form of $\lambda_i$.
One should also note that  equation (\ref{EqTheEqForOverlineLambda})
is stated for unitary and not for projective unitary operators.

We define a family (of local isomorphisms ) 
$\vartheta_i: U_i\to {\rm Map}(G/N,
\PU(\HH\otimes L^2(\widehat G/N^\perp)\otimes L^2(G/N)\otimes\HH))$ by
\begin{eqnarray*}
\vartheta_i(u)(z)
&:=&
\lambda_i(u)(z)\ 
(\Eins_\HH\otimes\Eins_{L^2(\widehat  G/N^\perp)}\otimes\kappa^a_i(u,z))\\ 
&&{\rm Ad}\Big( 
(\Eins_\HH\otimes \overline\kappa^\sigma(z,..)\otimes\Eins_\HH)\ 
(\Eins_\HH\otimes v_i(u,z,..))\Big)
\end{eqnarray*}
which we use to define
another family (of cocycles)
$$
\mu^{\rm dyn}_i:U_i\to Z^1_{\rm cont}(G,{\rm Map}(G/N,\PU(\HH\otimes 
L^2(\widehat G/N^\perp)\otimes L^2(G/N)\otimes\HH))
$$
by
\begin{eqnarray}
\mu^{\rm dyn}_i(u)(h,z)
&:= &\vartheta_i(u)(z+hN)\nonumber\\ 
&&{\rm Ad}\Big((m_i(u)(h,z)\otimes\Eins_{L^2(G/N)}\otimes\Eins_{\HH})\nonumber\\
&&(\Eins_\HH\otimes \langle\hat\sigma(..),h\rangle\otimes\Eins_{L^2(G/N)}\otimes\Eins_{\HH})\Big)
\nonumber\\
&& \vartheta_i(u)(z)^{-1}.\label{EqTheDefiOfTheCocycleFromATopTrip}
\end{eqnarray}
The stabilised pair $(P_{\rm dyn},E_{\rm dyn}):=(\PU(\HH\otimes L^2(\widehat G/N^\perp))\otimes P,E)$ \label{PageOfPDynStbilll}
has transition functions $g_{ji},\zeta^{\rm dyn}_{ji}:=\Eins_\HH\otimes\Eins_{L^2(\widehat G/N^\perp)}\otimes\zeta_{ji}$, and we claim that
on $(P_{\rm dyn},E_{\rm dyn})$ the family $\{\mu^{\rm dyn}_i\}_{i\in I}$ defines a dualisable \label{pageOfDefiOfDualisableDecker}
decker $\rho^{\rm dyn}$. 
To verify this we have to check that
$$
\mu^{\rm dyn}_j(u)(h,g_{ji}(u)+z)^{-1}\  \zeta^{\rm dyn}_{ji}(u)(z+hN)\ 
\mu^{\rm dyn}_i(u)(h,z)=
\zeta^{\rm dyn}_{ji}(u)(z).
$$
Although lengthy, this is a straight forward calculation. 
Indeed, we have $^\nou$
\begin{eqnarray*}
&&\mu^{\rm dyn}_j(h,g_{ji}+z)^{-1}\ \zeta^{\rm dyn}_{ji}(z+hN)\
\mu^{\rm dyn}_i(h,z)\\
&=&
\vartheta_j(g_{ji}+z)
{\rm Ad}((m_j(h,g_{ji}+z)(..)^{-1}\otimes\Eins\otimes\Eins)(\Eins\otimes\langle\hat\sigma(..),-h\rangle\otimes\Eins\otimes\Eins))\\
&&\vartheta_j(g_{ji}+z+hN)^{-1} 
(\Eins\otimes\Eins\otimes \zeta_{ji}(z+hN))
\vartheta_i(z+hN)\\
&&{\rm Ad}((m_i(h,z)(..)\otimes\Eins\otimes\Eins)(\Eins\otimes \langle\hat\sigma(..),h\rangle\otimes\Eins\otimes\Eins))
\vartheta_i(z)^{-1}\\
&=&
\vartheta_j(g_{ji}+z)
{\rm Ad}((m_j(h,g_{ji}+z)(..)^{-1}\otimes\Eins\otimes\Eins) 
(\Eins\otimes\langle\hat\sigma(..),-h\rangle\otimes\Eins\otimes\Eins))\\
&&{\rm Ad}\Big( (\Eins\otimes v_j(g_{ji}+z+hN,..)^{-1}) 
(\Eins\otimes \overline\kappa^\sigma(g_{ji}+z+hN,..)^{-1}\otimes\Eins) \Big)\\
&&(\Eins\otimes{\rm Ad}(\lambda_{_{\widehat G/N^\perp}}(\hat g_{ji}))
\otimes \zeta'_{ji}(z+hN))\\
&&{\rm Ad}\Big(
(\Eins\otimes \overline\kappa^\sigma(z+hN,..)\otimes\Eins)
(\Eins\otimes v_i(z+hN,..))\Big)\\
&&{\rm Ad}((m_i(h,z)(..)\otimes\Eins\otimes\Eins)
(\Eins\otimes \langle\hat\sigma(..),h\rangle\otimes\Eins\otimes\Eins))
\vartheta_i(z)^{-1}\\
&=&
\vartheta_j(g_{ji}+z)
{\rm Ad}((m_j(h,g_{ji}+z)(..)^{-1}\otimes\Eins\otimes\Eins)
(\Eins\otimes \langle\hat\sigma(..),-h\rangle\otimes\Eins\otimes\Eins))\\
&&{\rm Ad}(\Eins\otimes\lambda_{_{\widehat G/N^\perp}}(g_{ji})\otimes\Eins\otimes\Eins)\\
&&\Eins\otimes\Big[{\rm Ad}\Big( v_j(g_{ji}+z+hN,\hat g_{ji}+..)^{-1}\ 
(\overline\kappa^\sigma(z+hN,\hat g_{ji}+..)^{-1}\otimes\Eins)\Big)\\
&& \zeta'_{ji}(z+hN))\ {\rm Ad}\Big((\overline\kappa^\sigma(z+hN,..)\otimes\Eins)\ v_i(z+hN,..)\Big)\Big]\\
&&{\rm Ad}((m_i(h,z)(..)\otimes\Eins\otimes\Eins)
(\Eins\otimes\langle\hat\sigma(..),h\rangle\otimes\Eins\otimes\Eins))
\vartheta_i(z)^{-1}.
\end{eqnarray*}
The expression inside the squared brackets can be 
rewritten by equation (\ref{EqTheDefiOfAlpha}). This reads
\begin{eqnarray*}
\dots&=&
\vartheta_j(g_{ji}+z)
{\rm Ad}((m_j(h,g_{ji}+z)(..)^{-1}\otimes\Eins\otimes\Eins)
(\Eins\otimes \langle\hat\sigma(..),-h\rangle\otimes\Eins\otimes\Eins))\\
&&{\rm Ad}(\Eins\otimes\lambda_{_{\widehat G/N^\perp}}(g_{ji})\otimes\Eins\otimes\Eins)\\
&&\Eins\otimes\Big[ {\rm Ad}(\overline{\hat\zeta'}_{ji}(..))\ 
(A_{ji}(z+hN)\otimes\Eins\otimes\Eins)\Big]\\
&&{\rm Ad}((m_i(h,z)(..)\otimes\Eins\otimes\Eins)
(\Eins\otimes \langle\hat\sigma(..),h\rangle\otimes\Eins\otimes\Eins))
\vartheta_i(z)^{-1}.
\end{eqnarray*}
If we insert furthermore the results of the previous lemmata, 
we obtain
\begin{eqnarray}
\dots&=&
\vartheta_j(z)
{\rm Ad}((m_j(h,g_{ji}+z)(..)^{-1}\otimes\Eins\otimes\Eins)
(\Eins\otimes \langle\hat\sigma(..),-h\rangle\otimes\Eins\otimes\Eins))\nonumber\\
&&{\rm Ad}(\Eins\otimes\lambda_{_{\widehat G/N^\perp}}(g_{ji})\otimes\Eins\otimes\Eins)\label{EqStep21}\\
&& {\rm Ad}(\Eins\otimes\overline{\hat\zeta'}_{ji}(..))\ 
(\Eins\otimes\beta_{ji}(z+hN)\otimes\Eins\otimes\Eins)
 ({\rm Ad}(w_{ji}(z+hN)(..))\otimes\Eins\otimes\Eins)\nonumber\\
&&(\gamma_{ji}\otimes\Eins\otimes\Eins)
{\rm Ad}((m_i(h,z)(..)\otimes\Eins\otimes\Eins)
(\Eins\otimes \langle\hat\sigma(..),h\rangle\otimes\Eins\otimes\Eins))
\vartheta_i(z)^{-1}\nonumber\\
&=&
\vartheta_j(z)\ {\rm Ad}(\Eins\otimes\lambda_{_{\widehat G/N^\perp}}(g_{ji})\otimes\Eins\otimes\Eins)\nonumber\\
&&{\rm Ad}\Big(
(\Eins\otimes\overline{\hat\zeta'}_{ji}(..))
(\Eins\otimes \langle\hat\sigma(..+\hat g_{ji})-\hat\sigma(..),\sigma(g_{ji}+z+hN)-h\rangle\otimes\Eins\otimes\Eins)\nonumber\\
&&(w_{ji}(z)(..)\otimes\Eins\otimes\Eins)\Big)
(\gamma_{ji}\otimes\Eins\otimes\Eins) \vartheta_i(z)^{-1}.\label{EqStep22}
\end{eqnarray}
As $\sigma(g_{ji}+z+hN)-h$ and $\sigma(g_{ji}+z)$ differ by some element
in $N$, we obtain the identity 
\begin{eqnarray}
&&\langle\hat\sigma(..+\hat g_{ji})-\hat\sigma(..),\sigma(g_{ji}+z+hN)-h\rangle\nonumber\\
&=&
\langle\hat\sigma(..+\hat g_{ji})-\hat\sigma(..),\sigma(g_{ji}+z+hN)-h-
\sigma(g_{ji}+z)\rangle\nonumber\\
&&\langle\hat\sigma(..+\hat g_{ji})-\hat\sigma(..),\sigma(g_{ji}+z)\rangle\nonumber\\
&=&
\langle \hat g_{ji}(u),\sigma(g_{ji}+z+hN)-h-
\sigma(g_{ji}+z)\rangle\label {EqTheScalarTermofStep22}\\
&&\langle\hat\sigma(..+\hat g_{ji})-\hat\sigma(..),\sigma(g_{ji}+z)\rangle.\nonumber
\end{eqnarray}
Therein the first factor is a scalar and the second defines $\beta_{ji}$.
So the main calculation continues 
\begin{eqnarray}
\dots&=&
\vartheta_j(z)\ {\rm Ad}(\Eins\otimes\lambda_{_{\widehat G/N^\perp}}(g_{ji})\otimes\Eins\otimes\Eins)
{\rm Ad}(
(\Eins\otimes\overline{\hat\zeta'}_{ji}(..)))\label{EqStep23}\\
&&(\Eins\otimes\beta_{ji}(z)\otimes\Eins\otimes\Eins)
({\rm Ad}(w_{ji}(z)(..))\otimes\Eins\otimes\Eins)
(\gamma_{ji}\otimes\Eins\otimes\Eins) \vartheta_i(z)^{-1},\nonumber
\end{eqnarray}
and again by equation (\ref{EqTheDefiOfAlpha})
\begin{eqnarray}
\dots&=&
\vartheta_j(z)\ {\rm Ad}(\Eins\otimes\lambda_{_{\widehat G/N^\perp}}(g_{ji})\otimes\Eins\otimes\Eins)\nonumber\\
&&\Eins \otimes \Big[
{\rm Ad}\Big( v_j(g_{ji}+z,\hat g_{ji}+..)^{-1}\ 
(\overline\kappa^\sigma(z,\hat g_{ji}+..)^{-1}\otimes\Eins)\Big)\nonumber\\
&& (\Eins\otimes\zeta'_{ji}(z))
{\rm Ad}\Big((\overline\kappa^\sigma(z,..)\otimes\Eins)\ v_i(z,..)\Big)\Big]
 \vartheta_i(z)^{-1}\label{EqStep24}\\
&=&
\vartheta_j(z)\ \Eins \otimes {\rm Ad}\Big( v_j(g_{ji}+z,..)^{-1}\ 
(\overline\kappa^\sigma(z,..)^{-1}\otimes\Eins)\Big)\nonumber\\
&&{\rm Ad}(\Eins\otimes\lambda_{_{\widehat G/N^\perp}}(g_{ji}))
\otimes\zeta'_{ji}(z)\nonumber\\
&&\Eins\otimes{\rm Ad}\Big((\overline\kappa^\sigma(z,..)\otimes\Eins)\ v_i(z,..)\Big)
 \vartheta_i(z)^{-1}\nonumber\\
&=&
\Eins\otimes\Eins\otimes\zeta_{ji}(z)\nonumber\\
&=&\zeta^{\rm dyn}_{ji}(z).\nonumber
\end{eqnarray}
This shows that the $\mu^{\rm dyn}_i$s define a decker $\rho^{\rm dyn}$ on
$(P_{\rm dyn},E_{\rm dyn})$.
By construction $\rho^{\rm dyn}$ is dualisable, because 
$(h,z)\mapsto (m_i(u)(h,z)(..)\otimes\Eins\otimes\Eins)
(\Eins\otimes \langle\hat\sigma(..),h\rangle\otimes\Eins\otimes\Eins)$ is 
a continuous  and unitary implemented 1-cocycle.
\\\\
{\bf Step 2:}
The construction of Step 1
defines a map $\delta^{\rm as}(B):{\rm Top}^{\rm as}(B)\to \mathcal P({\rm Dyn}^\dag(B))$
and $\delta^{\rm as}(B)(x)$ is a $Q_{_G}$-torsor, for each $x\in{\rm Top^{as}}(B).$
If $x$ is strict then there is a distinguished $q(Q_{_{G/N}})$-subtorsor $X^{\rm dyn}$,
and if $x$ is in the image of $\tau(B)$ we single out a specific  element
$x^0\in X^{\rm dyn}$.
We then just define
$\delta^{\rm s}(B)(x):=X^{\rm dyn}$ and
$\delta^{\rm im}(B)(x):=x^0$ in the respective cases.
\\\\
We have to show that all choices involved do not change the class
of the dynamical triple $(\rho^{\rm dyn}, P_{\rm dyn},E_{\rm dyn})$.
That the choice of the atlas has no effect on the class of the 
constructed dynamical triple is rather obvious.
It is less obvious for the choices of 
$\lambda_i$ (or $l_i$), $m_i$ and the homotopy class of 
$\kappa$, i.e. the choice of an isomorphic
topological triple.
We convince ourselves that three other choices of these
define exterior equivalent deckers.

Firstly, 
if we choose another $\lambda_i$, say
$\lambda'_i={\rm Ad}(l'_i\otimes\Eins\otimes\Eins)$,
then by equation (\ref{EqTheEqForOverlineLambda})
the family 
$\nu_i(u)(z):=\lambda_i(u)(z) {\lambda'}_i^{-1}(u)(z)$ defines an automorphism 
of $P_{\rm dyn}$, and the  $\check{\rm C}$ech class of this automorphism 
 $[\nu_.]\in\check H^1(B,\underline{\U(1)})$
vanishes. Thus we are precisely in the situation of Example \ref{ExaOfExteriorEquivalence}.

Secondly,
equation (\ref{EqTheDefiOfMi}) shows that
$m_i$ is unique up to 
functions $n_i:U_i\to Z^1_{\rm cont}(G,{\rm Map}(G/N,L^\infty(\widehat G/N^\perp,\UA(\HH))))$
such that $\delta_{g\times \hat g}n_.=\Eins$.
We change the atlas of the the constructed dynamical
triple such that we have transition functions
$\zeta^1_{ji}(u)(z):=\vartheta_j(u)(g_{ji}(u)+z)^{-1}
\zeta_{ji}(u)(z)\vartheta_i(u)(z)$
and cocycles
$\mu_i^1(u)(h,z):={\rm Ad}((m_i(u)(h,z)\otimes\Eins_{L^2(G/N)}\otimes\Eins_{\HH})
(\Eins_\HH\otimes \langle\hat\sigma(..),h\rangle\otimes\Eins_{L^2(G/N)}\otimes\Eins_{\HH})).
$
Assume  $m_i$ is changed by $n_i$. Then
because of the commutativity of $\UA(\HH)$
we have $n_i(u)(h+g,z)=\mu_i^1(u)(g,z)^{-1}(n_i(u)(h,z+gN)) n_i(u)(g,z)$;
and because of the $\lambda_i$-terms in $\vartheta_i$ we
have that
$\zeta^1_{ji}(u)(z)(n_i(u)(h,z)(..))
=n_i(u)(h, z)(..-\hat g_{ji}(u))
=n_j(u)(h,g_{ji}(u)+z)(..)$.
So $c_i:=n_i$ defines an exterior equivalence
(Definition \ref{DefiOfExteriorEquivalence}).

Thirdly,
if $\kappa'$ is homotopic to $\kappa$,
the $\check{\rm C}$ech class of the bundle automorphism
 $\kappa'\circ\kappa^{-1}$ vanishes,
and the change of $\mu_i^{\rm dyn}$ caused by 
this automorphism is  again 
covered by Example \ref{ExaOfExteriorEquivalence}.
\\\\
Let us now discuss the choice of $\chi_{..}$.
We always have the freedom to change 
it by a cocycle $\chi^1\in \check Z^1(U_\bullet,\underline{
Z^1_{\rm cont}(G,{\rm Map}(G/N,\U(1)))},g_{..})$.
The consequence is a change of 
$m_i$ by $m^1_i:U_i\to Z^1_{\rm cont}(G,{\rm Map}(G/N,L^\infty(\widehat G/N^\perp,\UA(\HH))))$ such that $\delta_{g\times \hat g}m^1_.=\chi^1_{..}$.
It is clear that if $\chi^1_{..}=\delta_g\chi^2_.$ is a boundary,
then the deckers $\rho^{\rm dyn}$ and ${\rho^{\rm dyn}}^1$ which 
correspond to $m_i$ and $m_i\  m_i^1$ are exterior equivalent, 
for we can define $n_i:=\chi^2_i m^1_i$, so $\delta_{g\times \hat g}n_.=\Eins$
which leads us to the case we already discussed above.
If the class $[\chi^1_{..}]$ does not vanish, then it follows 
that the corresponding deckers are not exterior equivalent 
and the corresponding triples are not stably outer conjugate.

Thus it follows that for each $[\chi^1_{..}]\in Q_{_G}$ we get
a different dynamical triple, and we define 
$\delta^{\rm as}(B)(x)\subset{\rm Dyn^\dag}(B)$
to be the set of these dynamical triples.

It is obvious that $Q_{_G}$ acts freely  and transitively
on $\delta^{\rm as}(B)(x)$, i.e. it is a $Q_{_G}$-torsor.
If the triple $x$ is strict, then we can choose
$\chi_{..}$ such that
$\chi_{ji}(u)|_N=1$, and this property is preserved by the
action of $q(Q_{_{G/N}})\subset Q_{_G}$. 
So we singled out a specific $q(Q_{_{G/N}})$-subtorsor 
$X^{\rm dyn}\subset\delta^{\rm as}(B)(x)$, and we 
define $\delta^{\rm s}(B)(x):=X^{\rm dyn}.$
If $x$ is in the image of $\tau(B)$, then we saw in 
in the proof of Lemma \ref{LemImTauSubsetStrict} 
that $w_{ji}$ satisfies $d_*w_{ji}(u)(z)\in\U(1)\cdot \Eins$, so
it is meaningful to define $\chi_{ji}:=d_*w_{ji}$.
Then we let 
$x^0\in X^{\rm dyn}$ be the element which corresponds
to this particular choice of $\chi_{..}$, and we put $\delta^{\rm im}(B)(x):=x^0$.
\\\\
{\bf Step 3:}
$\delta^{\rm im}(B)\circ\tau(B)={\rm id}_{{\rm Dyn}^\dag(B)}$,
so $\tau(B)$ is injective and $\delta(B)$ is surjective.
\\\\
The formal calculation is similar to to what we did in the of 
the proof of Theorem \ref{ThmTheMainTheoremOfDualityTheoryOfDynTriples}.

Let $(\rho, P, E)$ be a dualisable dynamical triple having 
transition functions $g_{ji},\zeta_{ji}$ and cocycles $\mu_i$.
Recall the definition of $\tau(B)$ in particular of the topological triple
$(\kappa^{\rm top},(P_{\rm top},E),(\widehat P,\widehat E))$ out
of which we must compute 
$(\rho^{\rm dyn},{P_{\rm top}}_{\rm dyn},E)$.

${P_{\rm top}}_{\rm dyn}=
\PU(\HH\otimes L^2(\widehat G/N^\perp)\otimes L^2(G/N))\otimes P$ 
is stably isomorphic to $P$
and we shall compute $\mu^{\rm dyn}_i$.
We first give one possible, explicit formula for $\lambda_i(u)(z)$.
Let $\hat{\mathcal F}:L^2(\widehat G/N^\perp)\to L^2(N)$ \label{PageOfInvFourier}
be the inverse Fourier  transform, then 
$$
\hat {\mathcal F}\circ\lambda_{_{\widehat G/N^\perp}}(\hat g_{ji}(u))\circ
\hat{\mathcal F}^{-1}=\langle\hat g_{ji}(u),\ .\ \rangle^{-1}\in L^\infty(N,\U(1)),
$$
and the definition of $\hat g_{ji}$ implies
\begin{eqnarray*}
&&{\rm Ad}(\lambda_{_{\widehat  G/N^\perp}}(\hat g_{ji}(u)))
\otimes\zeta^{\rm top}_{ji}(u)(z)\\
&=&
{\rm Ad}\Big((\hat {\mathcal F}^{-1}\otimes\Eins\otimes\Eins)\circ
\overline\mu_j(u)(\ .\  ,g_{ji}(u)+z)^{-1}|_N\circ (\hat{\mathcal F}\otimes\Eins\otimes\Eins)\Big)\\
&&\Eins\otimes\zeta^{\rm top}_{ji}(u)(z)
{{\rm Ad}\Big((\hat {\mathcal F}^{-1}\otimes\Eins\otimes\Eins)\circ
\overline\mu_i(u)(\ .\ ,z)|_N\circ (\hat{\mathcal F}\otimes\Eins\otimes\Eins)\Big)}.
\end{eqnarray*}
Therefore we have 
\begin{eqnarray*}
l_i(u)(z)&=&
((\Eins\otimes\hat {\mathcal F}^{-1}\otimes\Eins\otimes\Eins)\circ
(\Eins\otimes\overline\mu_i(u)(\ .\ ,z)|_N)\circ 
(\Eins\otimes\hat{\mathcal F}\otimes\Eins\otimes\Eins)\\
&&\in\U(\HH\otimes L^2(\widehat G/N^\perp)\otimes L^2(G/N)\otimes \HH)
\end{eqnarray*}
which is continuous by Lemma \ref{LemATechLemmaForMu} $(i)$,
hence $\lambda_i(u)(z)={\rm Ad}(l_i(u)(z))$.
From the local definition of $\kappa^{\rm top}$ in equation 
(\ref{EqTheLocalDefinitionOfKappaTop}) we read off that
$$
\kappa_i(u)(z,\hat z)=\kappa^a_i(u,z)\ \kappa^\sigma(z,\hat z)\ {\rm Ad}(v_i(u,z,\hat z))$$
for $\kappa^a_i(u,z):={\rm Ad}(\overline\mu_i(u)(-\sigma(\_),z)^{-1})$ and
$v_i(u,z,\hat z):=\lambda_{_{G/N}}(z)\otimes\Eins$.
In the proof of Lemma \ref{LemImTauSubsetStrict} we already observed that in the case 
of a topological triple constructed out of dynamical one 
we have $d_*w_{ji}(u)(z)\in\U(1)\cdot\Eins$, and by the definition 
of $\delta^{\rm im}(B)$ we  choose $\chi_{..}:=d_*w_{..}$.
Therefore we can choose  $m_i=\Eins$. As a consequence the first $\HH$-slot of the tensor product 
$\HH\otimes L^2(\widehat G/N^\perp)\otimes L^2(G/N)\otimes \HH$ will 
contain the identity operator $\Eins$ only.
We have all ingredients together to compute $\mu_i^{\rm dyn}$ $^\nou$.
\begin{eqnarray*}
\mu_i^{\rm dyn}(h,z)&=&
\vartheta_i(z+hN)
(\Eins\otimes\langle\hat\sigma(..),h\rangle\otimes\Eins\otimes\Eins)
\vartheta_i(z)^{-1}\\
&=&
{\rm Ad}\Big((\Eins\otimes\hat {\mathcal F}^{-1}\otimes\Eins\otimes\Eins)\circ
\Eins\otimes\overline\mu_i(\ .\ ,z+hN)|_N\circ 
(\Eins\otimes\hat{\mathcal F}\otimes\Eins\otimes\Eins)\\
&&(\Eins\otimes\Eins\otimes \overline\mu_i(-\sigma(\_),z+hN)^{-1})\\
&&(\Eins\otimes\overline\kappa^\sigma(z+hN,..))
(\Eins\otimes\Eins\otimes\lambda_{_{G/N}}(z+hN)\otimes\Eins)))\\
&&(\Eins\otimes\langle\hat\sigma(..),h\rangle\otimes\Eins\otimes\Eins)\\
&&(\Eins\otimes\Eins\otimes\lambda_{_{G/N}}(-z)\otimes\Eins)
(\Eins\otimes\overline\kappa^\sigma(z+hN,..)^{-1} ))\\
&&(\Eins\otimes\Eins\otimes\overline\mu_i(-\sigma(\_),z))\\
&&(\Eins\otimes \hat {\mathcal F}^{-1}\otimes\Eins\otimes\Eins)\circ
\Eins\otimes (\overline\mu_i(\ .\ ,z)|_N)^{-1}\circ 
(\Eins\otimes\hat{\mathcal F}\otimes\Eins\otimes\Eins)\Big).
\end{eqnarray*}
Since $(\hat{\mathcal F}\otimes\Eins)\circ\langle\hat\sigma(..),h-\sigma(\_+hN)+\sigma(\_)\rangle\circ(\hat{\mathcal F}^{-1}\otimes\Eins)=\lambda_{_{N}}(h-\sigma(\_+hN)+\sigma(\_))
\in \U(L^2(\widehat G/N^\perp)\otimes L^2(G/N))$ and 
by the cocycle identity for $\overline\mu_i$ this transforms to
\begin{eqnarray*}
\dots&=&
{\rm Ad}\Big((\Eins\otimes\hat {\mathcal F}^{-1}\otimes\Eins\otimes\Eins)\circ
(\Eins\otimes\overline\mu_i(\sigma(\_)+\ .\ ,z+hN-\_))\\
&&(\Eins\otimes \lambda_{_{G/N}}(hN)\otimes\Eins\otimes\Eins)
(\Eins\otimes\lambda_{_{N}}(h-\sigma(\_+hN)+\sigma(\_)))\otimes\Eins)\\
&&(\Eins\otimes
(\overline\mu_i(\sigma(\_)+\ .\ ,z-\_))^{-1})\circ 
(\Eins\otimes\hat{\mathcal F}\otimes\Eins\otimes\Eins)\Big)\\
&=&
{\rm Ad}\Big((\Eins\otimes\hat {\mathcal F}^{-1}\otimes\Eins\otimes\Eins)\circ\\
&&(\Eins\otimes \lambda_{_{G/N}}(hN)\otimes\Eins\otimes\Eins)
(\Eins\otimes\lambda_{_{N}}(h-\sigma(\_+hN)+\sigma(\_)))\otimes\Eins)\\
&&
(\Eins\otimes\overline\mu_i(\sigma(\_+hN)+\ .\ +h-\sigma(\_+hN)+\sigma(\_),z-\_))\\
&&(\Eins\otimes
(\overline\mu_i(-\sigma(\_)-\ .\ ,z)))\\
&&\circ 
(\Eins\otimes\hat{\mathcal F}\otimes\Eins\otimes\Eins)\Big)\\
&=&
{\rm Ad}\Big((\Eins\otimes\hat {\mathcal F}^{-1}\otimes\Eins\otimes\Eins)\circ\\
&&(\Eins\otimes \lambda_{_{G/N}}(hN)\otimes\Eins\otimes\Eins)
(\Eins\otimes\lambda_{_{N}}(h-\sigma(\_+hN)+\sigma(\_)))\otimes\Eins)\\
&&(\Eins\otimes\Eins\otimes\Eins\otimes
(\overline\mu_i(h,z)))\circ 
(\Eins\otimes\hat{\mathcal F}\otimes\Eins\otimes\Eins)\Big).
\end{eqnarray*}
Let $\mathcal S: L^2(N)\otimes L^2(G/N)\to L^2(G)$ be the isomorphism
introduced on page \pageref{PageThePageOfS}. There we discussed its 
behaviour with respect to the left regular representation on $L^2(G)$. 
This leads us finally to 
\begin{eqnarray*}
\mu^{\rm dyn}_i(u)(h,z)&=&
{\rm Ad}\Big(
(\Eins\otimes\hat {\mathcal F}^{-1}\otimes\Eins\otimes\Eins)\circ
(\Eins\otimes {\mathcal S}^{-1}\otimes\Eins)\circ\\
&&(\Eins\otimes \lambda_{_{G}}(h)\otimes\overline\mu_i(u)(h,z)))\\
&&\circ 
(\Eins\otimes {\mathcal S}\otimes\Eins)\circ
(\Eins\otimes\hat{\mathcal F}\otimes\Eins\otimes\Eins)\Big).
\end{eqnarray*}
Thus we have shown that 
$(\rho,P,E)$ and $(\rho^{\rm dyn},{P_{\rm top}}_{\rm dyn},E)$ 
are outer conjugate.
\\\\
{\bf Step 4:} We compute
$\tau(B)(\delta^{\rm s})(B)(x))$ which will
complete the statement of $(b)$.
In particular this calculation shows that 
$\tau(B)\circ\delta^{\rm im}(B)=
{\rm id}_{{\rm im}(\tau(B))}.$
\\\\
Let $(\kappa,(P,E),(\widehat P,\widehat E))$ be a topological
triple with underlying Hilbert space $\widehat \HH=L^2(G/N)\otimes\HH$.
For the first steps we only need the assumption that the triple is 
almost strict.  
So let $(\rho^{\rm dyn},P_{\rm dyn},E_{\rm dyn})$ be the dynamical triple
which we constructed in Step 1. As we saw  $E_{\rm dyn}$ is nothing but $E$
itself.
From this dynamical triple we are going to  construct 
the topological triple 
$$(\kappa^{\rm top},({P_{\rm dyn}}_{\rm top},E_{\rm dyn}),(\widehat{P_{\rm dyn}},\widehat{ E_{\rm dyn}}))$$ 
(according to the construction of $\tau(B)$)
which we then have to compare  with
$(\kappa,(P,E),(\widehat P,\widehat E))$.

To compute the dual pair $(\widehat{P_{\rm dyn}},\widehat{ E_{\rm dyn}})$
we with determination of  $\phi^{\rm dyn}_{ji}$ the second term
of the total cocycle $(\psi^{\rm dyn}_{...},\phi^{\rm dyn}_{..},1)$ of
the constructed dynamical triple $(\rho^{\rm dyn},P_{\rm dyn},E_{\rm dyn})$.
By equation (\ref{EqTheDualTorusCoycle}), the restriction of $\phi^{\rm dyn}_{ji}$ to $N$ defines 
the (the cocycle of) bundle $\widehat{E_{\rm dyn}}$.

The definition of the cocycles $\mu^{\rm dyn}_i$ of the decker $\rho^{\rm dyn}$
contains $\lambda_i$. We  make use of the possibility to choose 
$\lambda_i(u)(z)={\rm Ad}(l^0_i(u)\otimes\Eins\otimes\Eins)$ as we explained in
Step 1. In Step 2 we saw that this choice does not have any effect 
on the level of equivalence classes.
Assume that the atlas is sufficiently refined and consists of contractible 
charts such that the following continuous  lifts exist.
Lifts of the transition functions
$\overline\zeta_{ji}^{\rm dyn}=\Eins\otimes\Eins\otimes\overline\zeta_{ji}$
and lifts
$\overline\kappa^a_i:U_i\to L^\infty(G/N,\U(\widehat \HH))$
of  $\kappa^a$
which we use to define unitary lifts $\overline\mu_i^{\rm dyn}$
of $\mu^{\rm dyn}_i$ according to (\ref{EqTheDefiOfTheCocycleFromATopTrip})
in the obvious way, i.e. we drop the Ad. 
Let $\delta_{ji}:U_{ji}\to L^\infty(G/N,\U(1))$ be such that
$$
\alpha_{ji}(u)(z,\hat z)= \delta_{ji}(u)(z)\
\langle\hat\sigma(\hat g_{ji}(u)+\hat z)-\hat\sigma(\hat z),\sigma(g_{ji}(u)+z\rangle\
 w_{ji}(u)(z)(\hat z)\ \overline\gamma_{ji}(u), 
 $$
for
a lift $\overline\gamma_{ji}:U_{ji}\to\U(\HH)$ of $\gamma_{ji}$;
therein all notation is as above. 
Then by definition \label{PageOfPhiDyn}
\begin{eqnarray*}
&&\phi^{\rm dyn}_{ji}(u)(h,z)\\
&=&
\overline\mu_i^{\rm dyn}(u)(h,z)\ \overline\zeta^{\rm dyn}_{ji}(u)(z)^{-1}\
\overline\mu_j^{\rm dyn}(u)(h,g_{ji}(u)+z)^{-1}\ \overline\zeta_{ji}^{\rm dyn}(u)(z+hN),
\end{eqnarray*}
and if we repeat 
the calculation 
of Step 1 on the unitary level, 
there are four equalities which  
must be modified by $\U(1)$-valued functions. Namely,
equation (\ref{EqStep21}) by  $\delta(u)(z+hN)$,
equation (\ref{EqStep22}) by $\chi_{ji}$ from equation (\ref{EqTheDefiOfMi}),
equation (\ref{EqStep23}) by the scalar term of (\ref{EqTheScalarTermofStep22}) and
equation (\ref{EqStep24}) by $\delta(u)(z)^{-1}$.
We  finally find
\begin{eqnarray*}
\phi^{\rm dyn}_{ji}(u)(h,z)
&=&\langle \hat g_{ji}(u),\sigma(g_{ji}(u)+z+hN)-h-\sigma(g_{ji}(u)+z)\rangle\\
&&\chi_{ji}(u)(h,z)\ (d\delta_{ji}(u))(h,z).
\end{eqnarray*}
It follows that $\widehat {E_{\rm dyn}}\cong \widehat E$ 
if and only if 
the topological triple we started with is strict and we
have chosen $\chi_{..}$ such that
$\chi_{ji}(u)|_N=1$.
Indeed, the cocycles of these bundles are
$\hat g^{\rm dyn}_{ji}(u):=(\phi^{\rm dyn}_{ji}(u)(\ .\ ,z)|_N)^{-1}
=\hat g_{ji}(u)\ \chi_{ji}(u)(-\ .\ ,z)|_N$.
This proves $(b)$.

To complete the computation of the dual
we have to  compute
$$\widehat{\zeta_{ji}^{\rm dyn}}:U_{ji}\to
{\rm Map}(\widehat G/N^\perp,\PU(L^2(G/N)\otimes\HH\otimes L^2(\widehat G/N^\perp)\otimes
L^2(G/N)\otimes\HH)).
$$
by equation (\ref{EqTheDualPUCocycle}).
The reader should not be confused by the two
different $L^2(G/N)$ factors in the tensor product.
The first is due to the stabilisation in the definition
of the dual  (\ref{EqTheDualPUCocycle}), 
and the second is due to the Hilbert space
we started with which is $L^2(G/N)\otimes\HH$.
We use the symbols $_\smile$ and $\_$ to distinguish   
multiplication operators on the two Hilbert spaces;
$_\smile$ for the first factor, due to the definition
of the dual and $\_$ for the second factor as we did
all the time.
The dual transition functions are given by $^\nou$
\begin{eqnarray*}
\widehat{\zeta_{ji}^{\rm dyn}}(\hat z)
&=&{\rm Ad}\Big(
(\overline\kappa^\sigma(-g_{ji},\hat g^{\rm dyn}_{ji}+\hat z)\otimes\Eins
\otimes\Eins\otimes\Eins\otimes\Eins)\\
&&
(\lambda_{_{G/N}}(-g_{ji})\otimes\Eins\otimes\Eins\otimes\Eins\otimes\Eins)\\
&&(\overline\zeta^{\rm dyn}_{ji}(-_\smile ))
(\phi_{ji}^{\rm dyn}(-\sigma(_\smile),0)^{-1}\otimes\Eins\otimes\Eins\otimes\Eins\otimes\Eins)\Big)\\
&=&{\rm Ad}\Big(
(\overline\kappa^\sigma(-g_{ji},\hat g^{\rm dyn}_{ji}+\hat z)\otimes\Eins
\otimes\Eins\otimes\Eins\otimes\Eins)\\
&&
(\lambda_{_{G/N}}(-g_{ji})\otimes\Eins\otimes\Eins\otimes\Eins\otimes\Eins)\\
&&
(\Eins\otimes l_j^0\otimes\Eins\otimes\Eins)
[\Eins\otimes \lambda_{_{\widehat G/N^\perp}}(\hat g_{ji})\otimes
(\overline\zeta_{ji}(-_\smile )]\\
&&
(\Eins\otimes {l^0_i}^{-1}\otimes\Eins\otimes\Eins)
(\phi_{ji}^{\rm dyn}(-\sigma(_\smile),0)^{-1}\otimes\Eins\otimes\Eins\otimes\Eins\otimes\Eins)\Big).
\end{eqnarray*}
By equation (\ref{EqTheDefiOfAlpha})  this is $^\nou$
\begin{eqnarray*}
\dots&=&{\rm Ad}\Big(
(\Eins\otimes l_j^0\otimes\Eins\otimes\Eins)\\
&&(\overline\kappa^\sigma(-g_{ji},\hat g^{\rm dyn}_{ji}+\hat z)\otimes\Eins
\otimes\Eins\otimes\Eins\otimes\Eins)
(\lambda_{_{G/N}}(-g_{ji})\otimes\Eins\otimes\Eins\otimes\Eins\otimes\Eins)\\
&&
(\Eins\otimes\Eins\otimes\overline\kappa^a_j(g_{ji}-_\smile))
(\Eins\otimes\Eins\otimes\overline\kappa^\sigma(g_{ji}-_\smile,\hat g_{ji}+\hat z)\otimes\Eins)\\
&&
(\Eins\otimes\Eins\otimes v_j(g_{ji}-_\smile,\hat g_{ji}+\hat z))
(\Eins\otimes\Eins\otimes\overline\kappa^b_j(\hat g_{ji}+\hat z))\\
&&
(\Eins\otimes\Eins\otimes\lambda_{_{\widehat G/N^\perp}}(\hat g_{ji})
\otimes\overline{\hat\zeta}_{ji}(\hat z))
(\alpha_{ji}(-_\smile,\hat z)\otimes\Eins\otimes\Eins\otimes\Eins\otimes\Eins)
\\ && 
(\Eins\otimes\Eins\otimes\overline\kappa^b_i(\hat z)^{-1})
(\Eins\otimes\Eins\otimes v_i(-_\smile,\hat z)^{-1})\\
&&
(\Eins\otimes\Eins\otimes\overline\kappa^\sigma(-_\smile,\hat z)^{-1}\otimes\Eins)
(\Eins\otimes\Eins\otimes\overline\kappa^a_i(-_\smile)^{-1})\\
&&
(\phi_{ji}^{\rm dyn}(-\sigma(_\smile),0)^{-1}\otimes\Eins\otimes\Eins\otimes\Eins\otimes\Eins)\\
&&(\Eins\otimes {l^0_i}^{-1}\otimes\Eins\otimes\Eins)\Big),
\end{eqnarray*}
and if we insert the formulas for $\overline\kappa^\sigma,\alpha_{ji}$ and 
$\phi^{\rm dyn}_{ji}$ from above, we end up with\nopagebreak  $^\nou$
\begin{eqnarray*}
\dots&=&
{\rm Ad}\Big(
(\Eins\otimes l_j^0\otimes\Eins\otimes\Eins)\\
&&(\langle\hat\sigma(\hat g^{\rm dyn}_{ji}+\hat z),\sigma(_\smile+g_{ji})-\sigma(_\smile)
\rangle\otimes\Eins
\otimes\Eins\otimes\Eins\otimes\Eins)
\\&&
(\lambda_{_{G/N}}(-g_{ji})\otimes\Eins\otimes\Eins\otimes\Eins\otimes\Eins)\\
&&
(\Eins\otimes\Eins\otimes\overline\kappa^a_j(g_{ji}-_\smile))\\ 
&&
(\Eins\otimes\Eins\otimes
\langle\hat\sigma(\hat g_{ji}+\hat z),\sigma(\_-g_{ji}+_\smile)-\sigma(\_)
\rangle\otimes\Eins)\\
&&
(\Eins\otimes\Eins\otimes v_j(g_{ji}-_\smile,\hat g_{ji}+\hat z))
(\Eins\otimes\Eins\otimes\overline\kappa^b_j(\hat g_{ji}+\hat z))\\
&&
(\Eins\otimes\Eins\otimes\lambda_{_{\widehat G/N^\perp}}(\hat g_{ji})
\otimes\overline{\hat\zeta}_{ji}(\hat z))
( \delta_{ji}(-_\smile)\otimes\Eins\otimes\Eins\otimes\Eins\otimes\Eins)\\
&&(\langle\hat\sigma(\hat g_{ji}+\hat z)
-\hat\sigma(\hat z),\sigma(g_{ji}-_\smile)\rangle
\otimes\Eins\otimes\Eins\otimes\Eins\otimes\Eins)\\
&&
(w_{ji}(-_\smile)(\hat z)\otimes\Eins\otimes
\Eins\otimes\Eins)
(\Eins\otimes  \overline\gamma_{ji}\otimes
\Eins\otimes\Eins\otimes\Eins)\\ 
&& 
(\Eins\otimes\Eins\otimes\overline\kappa^b_i(\hat z)^{-1})
(\Eins\otimes\Eins\otimes v_i(-_\smile,\hat z)^{-1})\\
&&
(\Eins\otimes\Eins\otimes
\langle\hat\sigma(\hat z),-\sigma(\_+_\smile)+\sigma(\_)\rangle
\otimes\Eins)\\
&&(\Eins\otimes\Eins\otimes\overline\kappa^a_i(-_\smile)^{-1})\\
&&
([\langle \hat g_{ji},\sigma(g_{ji}-_\smile)+\sigma(_\smile)-\sigma(g_{ji})\rangle\
\chi_{ji}(-\sigma(_\smile),0)\\ 
&&(d\delta_{ji})(-\sigma(_\smile),0)\  ]^{-1}
\otimes\Eins\otimes\Eins\otimes\Eins\otimes\Eins)\\
&&(\Eins\otimes {l^0_i}^{-1}\otimes\Eins\otimes\Eins)\Big).
\end{eqnarray*}
It follows  from the definition 
of $\overline\gamma_{ji}$ that
${\rm Ad}(w_{ji}(u)(0)(\hat z)^{-1})=\gamma_{ji}(u)$,
so  after some further manipulation 
with all the $\langle\dots\rangle$-expressions this 
finally transforms to $^\nou$
 \begin{eqnarray*}
\dots&=&
{\rm Ad}\Big(
(\Eins\otimes l_j^0\otimes\Eins\otimes\Eins)
\\
&&
(\Eins\otimes\Eins\otimes\overline\kappa^a_j(-_\smile))\\ 
&&
(\Eins\otimes\Eins\otimes
\langle\hat\sigma(\hat g_{ji}+\hat z),-\sigma(_\smile)+\sigma(\_+_\smile)-\sigma(\_)
\rangle\otimes\Eins)\\
&&
(\Eins\otimes\Eins\otimes v_j(-_\smile,\hat g_{ji}+\hat z))
(\Eins\otimes\Eins\otimes\overline\kappa^b_j(\hat g_{ji}+\hat z))\\
&&\lambda_{_{G/N}}(-g_{ji})\otimes\Eins\otimes\lambda_{_{\widehat G/N^\perp}}(\hat g_{ji})
\otimes\overline{\hat\zeta}_{ji}(\hat z))\\
&&( w_{ji}(-_\smile)(\hat z)\ w_{ji}(0)(\hat z)^{-1}\ \chi_{ji}(-\sigma(_\smile),0)^{-1}
\otimes \Eins\otimes\Eins\otimes\Eins)\\
&&(\Eins\otimes\Eins\otimes\overline\kappa^b_i(\hat z)^{-1})
(\Eins\otimes\Eins\otimes v_i(-_\smile,\hat z)^{-1})\\
&&
(\Eins\otimes\Eins\otimes
\langle\hat\sigma(\hat z),\sigma(_\smile)-\sigma(\_+_\smile)+\sigma(\_)\rangle
\otimes\Eins)\\
&&(\Eins\otimes\Eins\otimes\overline\kappa^a_i(-_\smile)^{-1})\\
&&
(\Eins\otimes {l^0_i}^{-1}\otimes\Eins\otimes\Eins)\\
&&(\langle\hat\sigma(\hat g^{\rm dyn}_{ji}+\hat z)
-\hat\sigma(\hat g_{ji}+\hat z),-\sigma(_\smile-g_{ji})+\sigma(_\smile)\rangle
\otimes\Eins\otimes\Eins\otimes\Eins\otimes\Eins)
\Big).
\end{eqnarray*}
So far we did not use that  the topological triple 
we started with is strict and that we want to compute 
the composition $\tau(B)\circ\delta^{\rm s}(B)$.
 We use this now, so
the $G$-slot of $\chi_{ji}$ 
factors through $G/N$ and in particular
$\hat g^{\rm dyn}_{ji}=\hat g_{ji}$,
so we find 
\begin{eqnarray*}
\dots&=&
\eta'_j(u,\hat g_{ji}(u)+\hat z)\
(\xi_{ji}(u)(\hat z)\otimes\hat\zeta_{ji}(u)(\hat z))\
\eta'_i(u)(z)^{-1}
\end{eqnarray*}
wherein we have used the short hands 
\begin{eqnarray*}
\xi_{ji}(u)(\hat z)&:=&
{\rm Ad}\Big(
[(\lambda_{_{G/N}}(-g_{ji}(u))\otimes\Eins)
 w_{ji}(u)(-_\smile)(\hat z)\ w_{ji}(u)(0)(\hat z)^{-1}\\ 
 &&\chi_{ji}(u)(-\sigma(_\smile),0)^{-1}]\ 
 \otimes\lambda_{_{\widehat G/N^\perp}}(\hat g_{ji}(u))\Big)\\
&&\in\U(L^2(G/N)\otimes\HH\otimes L^2(\widehat G/N^\perp))
\end{eqnarray*}
and
\begin{eqnarray*}
\eta'_i(u)(z)&:=&
{\rm Ad}\Big(
(\Eins\otimes l^0_i(u)\otimes\Eins\otimes\Eins)
(\Eins\otimes\Eins\otimes\overline\kappa^a_i(u)(-_\smile))\\ 
&&
(\Eins\otimes\Eins\otimes
\langle\hat\sigma(\hat z),-\sigma(_\smile)+\sigma(\_+_\smile)-\sigma(\_)
\rangle\otimes\Eins)\\
&&
(\Eins\otimes\Eins\otimes v_i(u,-_\smile,\hat z))
(\Eins\otimes\Eins\otimes\overline\kappa^b_i(u)(\hat z))
\Big)
\\
&&\in\PU(L^2(G/N)\otimes\HH\otimes L^2(\widehat G/N^\perp)\otimes 
L^2(G/N)\otimes\HH).
\end{eqnarray*}
At this point we can read off that the bundle defined by 
$\xi_{ji}$ is trivialisable if $\chi_{..}^{-1}d_*w_{ji}=\Eins$
which is  the datum we must consider when we 
compute $\tau(B)\circ\delta^{\rm im}(B)$ (cp. Step 2).
In this case let 
$x_i:U_i\to {\rm Map}(\widehat G/N^\perp,\U(L^2(G/N)\otimes \HH\otimes L^2(\widehat G/N^\perp)))$ be such that
${\rm Ad}(x_j(u)(\hat g_{ji}(u)+\hat z) x_i(u)(\hat z)^{-1})=\xi_{ji}(u)(\hat z)$.
Then we have shown  that $\eta_i(u)(\hat z):=\eta'_i(u)(\hat z)\ {\rm Ad}(x_i(u)(\hat z)\otimes\Eins\otimes\Eins)$
defines a family of local isomorphisms
which fit together to a global isomorphism of principal bundles 
$\eta:{\widehat P}_{\rm s}\to\widehat {P_{\rm dyn}}$.
The subscript $\rm s$ denotes stabilisation with respect
to the Hilbert space $L^2(G/N)\otimes\HH\otimes L^2(\widehat G/N^\perp)$.
We claim that the topological triple 
$(\kappa,(P,E),(\widehat P,\widehat E))$
and the constructed triple 
$(\kappa^{\rm top},({P_{\rm dyn}}_{\rm top},E),(\widehat {P_{\rm dyn}},\widehat E))$
are equivalent; $\kappa^{\rm top}$ is defined by equation 
(\ref{EqTheLocalDefinitionOfKappaTop}) out of $\mu^{\rm dyn}_i$.
This will prove the identity $\tau(B)\circ\delta(B)={\rm id}_{{\rm im}(\tau(B)}.$
 We show that  the diagram
$$
\xymatrix{
P_{\rm s}\times_B\widehat E\ar[dd]^=&&E\times_B\widehat P_{\rm s}
\ar[ll]_{\Eins\otimes\Eins\otimes\Eins\otimes\kappa}\ar[dd]^{{\rm id}\times_{_B}\eta}\\
&&\\
{P_{\rm dyn}}_{\rm top}
\times_B\widehat E&&E\times_B\widehat {P_{\rm dyn}}\ar[ll]_{\kappa^{\rm top}}
}
$$
commutes up to homotpy. Locally, this means that
there exist continuous maps
$V_i:U_i\times G/N\times \widehat G/N^\perp\to\PU(L^2(G/N)\otimes\HH\otimes L^2(\widehat G/N^\perp)\otimes L^2(G/N)\otimes\HH)$ such that
$\kappa^{\rm top}_i(u)(z,\hat z)\ \eta_i(u)(\hat z)=
(\Eins\otimes\Eins\otimes\Eins\otimes\kappa_i(u)(z,\hat z))\ {\rm Ad}(V_i(u,z,\hat z))$.
To prove this we will use that the projective unitary group is homotopy comutative in 
the sense of Corollary \ref{CorTheHomotopyCommutativityOfPU}.
We have $^\nou$
\begin{eqnarray*}
&&\kappa^{\rm top}_i(z,\hat z)\ \eta_i(\hat z)\\
&=&
{\rm Ad}\Big(
(\langle\hat\sigma(\hat z),\sigma(_\smile-z)-\sigma(_\smile)\rangle
\otimes\Eins\otimes\Eins\otimes\Eins\otimes\Eins)\\
&&(\mu^{\rm dyn}_i(-\sigma(_\smile),z))^{-1}
(\lambda_{_{G/N}}(z)
\otimes\Eins\otimes\Eins\otimes\Eins\otimes\Eins)\Big) \eta_i(\hat z)\\
&=&
{\rm Ad}\Big(
(\langle\hat\sigma(\hat z),\sigma(_\smile-z)-\sigma(_\smile)\rangle
\otimes\Eins\otimes\Eins\otimes\Eins\otimes\Eins)\\
&&
(\Eins\otimes l^0_i\otimes\Eins\otimes\Eins)
(\Eins\otimes\Eins\otimes\Eins\otimes\overline \kappa^a_i(z))\\
&&(\Eins\otimes\Eins\otimes
\langle\hat\sigma(..),\sigma(\_-z)-\sigma(\_)\rangle
 \otimes\Eins)\\
&&(\Eins\otimes\Eins\otimes v_i(z,..))\\
&&(\Eins\otimes\langle\hat\sigma(..),\sigma(_\smile)\rangle\otimes\Eins\otimes\Eins)\\
&&(m_i(-\sigma(_\smile),z)(..)^{-1}\otimes\Eins\otimes\Eins)\\
&&(\Eins\otimes v_i(z-_\smile,..)^{-1})\\
&&(\Eins\otimes
\langle\hat\sigma(..),-\sigma(\_-z+_\smile)+\sigma(\_)\rangle
 \otimes\Eins)\\
&&(\Eins\otimes\Eins\otimes\Eins\otimes\overline \kappa^a_i(z-_\smile)^{-1})
(\Eins\otimes {l^0_i}^{-1}\otimes\Eins\otimes\Eins)\\
&&(\lambda_{_{G/N}}(z)
\otimes\Eins\otimes\Eins\otimes\Eins\otimes\Eins)\\
&&(\Eins\otimes l^0_i\otimes\Eins\otimes\Eins)
(\Eins\otimes\Eins\otimes\overline\kappa^a_i(-_\smile))\\ 
&&
(\Eins\otimes\Eins\otimes
\langle\hat\sigma(\hat z),-\sigma(_\smile)+\sigma(\_+_\smile)-\sigma(\_)
\rangle\otimes\Eins)\\
&&
(\Eins\otimes\Eins\otimes v_i(-_\smile,\hat z))
(\Eins\otimes\Eins\otimes\overline\kappa^b_i(\hat z))\\
&&(x_i(\hat z)\otimes\Eins\otimes \Eins)
\Big).
\end{eqnarray*}
We see that the terms $\kappa^a_i(u)(z)$ and $\kappa_i^b(u)(z)$  occur.
All terms with $l^0_i,m_i,v_i,x_i$ are continuous unitary maps
The only bracket $\langle\dots\rangle$ terms not continuous 
as unitary maps are the two terms in the first and third last line.
If we collect them, we shall better pay attention 
to $\lambda_{_{G/N}}(z)$ which acts
in the $_\smile$-variable. Then these two terms equal
\begin{eqnarray*}
&&\langle\hat\sigma(\hat z),
\sigma(\_+_\smile-z)-\sigma(\_)-\sigma(_\smile)\rangle\\
&=&
\langle\hat\sigma(\hat z),
\sigma(\_-z)-\sigma(\_)\rangle\\
&&\langle\hat\sigma(\hat z),
-\sigma(\_-z)+\sigma(\_+_\smile-z)-\sigma(_\smile)\rangle\\
&=&
(\Eins\otimes\overline \kappa^\sigma(z,\hat z))\
\langle \hat z,
-\sigma(\_-z)+\sigma(\_+_\smile-z)-\sigma(_\smile)\rangle\\
&&\in\U(L^2(G/N)\otimes L^2(G/N))
\end{eqnarray*}
and the second factor is continuous as a unitary map.
Therefore we have
\begin{eqnarray*}
\kappa^{\rm top}_i(u)(z,\hat z)\ \eta_i(u)(\hat z)&=\ldots=&
(\Eins\otimes\Eins\otimes\Eins\otimes\kappa^a_i(u)(z))\\
&&(\Eins\otimes\Eins\otimes\Eins\otimes\kappa^\sigma(z,\hat z))\\ 
&&(\Eins\otimes\Eins\otimes\Eins\otimes\kappa^b_i(u)(z))\
{\rm Ad}(V'_i(u,z,\hat z))\\
&=&
(\Eins\otimes\Eins\otimes\Eins\otimes\kappa_i(u)(z,\hat z)\  {\rm Ad}(V_i(u,z,\hat z)),
\end{eqnarray*}
for suitable continuous unitary maps $V'_i,V_i$,
and we are done.
\\\\
{\bf Step 5:} Naturality.
\\\\
In the same way as we already stated for $\tau$,
we  have a commutative diagram
$$
\xymatrix{
{\rm Top^?}(B)\ar[d]^{f^*}\ar[rr]^{\delta^?(B)}&&\mathcal P({\rm Dyn^\dag}(B))\ar[d]^{f^*}&\\
{\rm Top^?}(B')\ar[rr]^{\delta^?(B')}&&\mathcal P({\rm Dyn^\dag}(B')),&?={\rm as,\ s,\ im},
}
$$
for each continuous map between base spaces $f:B'\to B$ 
\end{pf}

\begin{defi}
Let $P\to E\to B$ be a pair over the base space $B$.
We say that this pair has an {\bf extension} 
to a topological  (resp. almost strict topological, strict topological, dynamical or 
dualisable dynamical) triple if 
the class $[(P,E)]\in {\rm Par}(B)$ is in the image 
of the forgetful map from topological 
(resp. almost strict topological, strict topological, dynamical or 
dualisable dynamical)
 triples over $B$ to pairs over $B$.
\end{defi}

\begin{cor}\label{CorTheCorOfTheMainTheorem}
Let $(P,E)$ be a pair over $B$.
Then the following are equivalent:
\begin{itemize}
\item[$(i)$]
$(P,E)$ has an extension
to a strict  topological triple;
\item[$(ii)$]
$(P,E)$ has an extension
to an almost strict  topological triple;
\item[$(iii)$]
$(P,E)$ has an extension
to a dualisable dynamical triple.
\end{itemize}
\end{cor}
\begin{pf}
This is an immediate consequence of the previous theorem.
\end{pf}

\subsection{The Case of $G=\R^n$ and $N=\Z^n$}\label{SecTheCaseOfRandZ}
\label{SecTheCaseOfRAndZ}

So far we kept our analysis completely general in 
the sense that we did not specify the groups $G,N$.
We now turn to the important case of $G=\R^n$
with lattice $N=\Z^n$, $n=1,2,3,\dots$.
In the whole of this section we use the notation 
$\T^n:=\R^n/\Z^n$ \label{PageOfTorus}
for the torus  and $\hat\T^n:=\widehat{\R^n}/{\Z^n}^\perp$
for the dual torus. There should be no confusion 
to decide between $\hat\T^n$ and the dual group
$\widehat{\T^n}\cong\Z^n$.
\\\\
The first thing we should check is that in case
of $G=\R^n$, $N=\Z^n$ our definition of topological triples agrees with the one 
introduced in \cite{BS}.

The definition of T-duality triples in \cite[Def. 2.8]{BS}
differs in two points from what we stated in Definiton \ref{DefiTopTriples}.
The first point is that they use the language of twists \cite[A.1]{BS}
instead of $\PU(\HH)$-principal fibre bundles to model 
T-duality diagrams.
But the category of $\PU(\HH)$-principal bundles with homotopy 
classes of bundle isomorphisms as morphisms is a model of twists,
and our notion of stable equivalence of topological triples
leads to the same equivalence classes as twists modulo
isomorphism. This, because the notion of equivalence we use 
for topological triples requires the  commutativity of 
diagram (\ref{DiagMorphOfTopTriples}) only up to homotopy.
To explain the second point let us consider the filtration of $H^3(E,\Z)$
associated to the Leray-Serre spectral sequence
$$
0\subset F^3H^3(E,\Z)\subset  F^2H^3(E,\Z)\subset  F^1H^3(E,\Z)\subset H^3(E,\Z).
$$ 
By definition, an element $h\in H^3(E,\Z)$ is in the subgroup $F^{k+1}H^3(E,\Z)$ 
if for any $k$-dimensional CW-complex $C$ and any map $f:C\to B$ $h$ is in the kernel of the induced map $f^*:H^3(E,\Z)\to H^3(C\times_B E,\Z)$.
In our definition of a topological triple $(\kappa,(P,E),(\widehat P,\widehat E))$
we require that the class $[P]\in H^3(E,\Z)$
of the bundle $P\to E$  lies in
the subgroup $F^1H^3(E,\Z)$
which is equivalent to the requirement of the triviality of $P$ over
the fibres of $E\to B$ (see Definition \ref{DefiPairOverB}); analogously for $\widehat P$.
In \cite{BS} the definition of T-duality triples requires 
that the class $[P]$ is even in the second step of the filtration
$[P]\in F^2H^3(E,\Z)$; analogously for $\widehat P$.
The following lemma states that in the case of $G=\R^n$ and $N=\Z^n$ these two conditions are equivalent.

\begin{lem}\label{LemTopTriplesAreTDualityTriples}
Let $G=\R^n$ and $N=\Z^n$,  and let $(\kappa,(P,E),(\widehat P,\widehat E))$ be a topological triple over $B$, then
$$[P]\in F^2H^3(E,\Z)\ \ {\rm  and}\ \ [\widehat P]\in F^2H^3(\widehat E,\Z).
$$
\end{lem}

\begin{pf}
Let $C$ be any 1-dimensional CW-complex and  $f:C\to B$  any 
continuous map. 
We have to show that the class  $[P]$ is in the kernel of the induced map $f^*:H^3(E,\Z)\to H^3(C\times_B E,\Z)$.
As $C$ is 1-dimensional, $H^2(C,\Z)=0$ and there are no non-trivial torus bundles over it. Thus, if we pull back the topological triple $(\kappa,(P,E),(\widehat P,\widehat E))$
along $f$, it becomes 
$$
(\kappa_C,(P_C,C\times_B E\cong C\times\T^n),(\widehat P_C,C\times_B\widehat E\cong C\times\hat\T^n))$$
and the centre  part of the  corresponding diagram 
degenerates to the projections
$$
\xymatrix{
&&C\times\T^n\times\hat\T^n\ar[dr]\ar[dl]&&\\
&C\times \T^n&& C\times \hat\T^n.&
}
$$
We extend this diagram to
$$
\xymatrix{
&&&C\times \T^n\ar[dl]\ar[dd]\ar[dr] &\\
&&C\times\T^n\times\hat\T^n\ar[dr]\ar[dl]&&C\ar[ld]\\
&C\times \T^n&& C\times \hat\T^n&,
}
$$
with the obvious inclusions and projection, so the diagonal
composition is the identity.
When we apply the the cohomology functor $H^3(\ .\ ,\Z)$ to this diagram
the vertical arrow becomes zero as it factors over $H^2(C,\Z)=0$, but
as $f^*[P]=[P_C]\in H^3(C\times\T^n,\Z)$ is mapped  by the identity from the left lower to the right upper group and as $[P_C]$ and $[\widehat P_C]$ equal when
pulled back to $H^3(C\times \T^n\times\hat\T^n,\Z)$, we conclude it is 
zero, since its image in the right upper group coincides with the image of $f^*[\widehat P]=[\widehat P_C]$ under the vertical arrow. 

This proves that $f^*[P]=0\in H^3(C\times_B E,\Z)$ and the same argument 
shows that the corresponding statement is true for $[\widehat P]$.
\end{pf}

So we observe that our functor 
${\rm Top}:\{ {\rm base\ spaces}\}\to \{{\rm sets}\}$
is the same functor which is introduced in \cite[Def. 2.11]{BS}
under the name ${\tt Triple}_n$, and below we are going to 
use a central result  of \cite{BS} about this functor, 
namely that the functor ${\rm Top}={\tt Triple}_n$  is
representable by a space $R_n$ (see Lemma \ref{ThmOfBSOnTheRepresentability} below).
\\\\
We now want to discuss Theorem 
\ref{ThmTheMainTheoremOfThisThesisNumberThree}
in the case of $G=\R^n, N=\Z^n$.
In this case the question  which of the 
topological triples are almost strict has a trivial
answer and also the torsor structure of Theorem 
\ref{ThmTheMainTheoremOfThisThesisNumberThree}
becomes trivial.
The criterion we make use of is the following lemma.

\begin{lem}
Let $Z$ be a topological abelian group which is 
contractible as topological space and is equipped with
a continuous (right) $G/N$-action.
Then 
$$\check H^k(B,\underline Z,[g_{..}])=0\qquad k=1,2.
$$
\end{lem}

\begin{pf}
Similar to the proof of Lemma \ref{LemTheCrucialLemmaToConstructDelta},
the proof makes  use of the standard Zorn's lemma argument.
\\\\
k=1:
Let $\{U_i\}_{i\in I}$ be an open cover of $B$, and let
$
\varphi_{..}\in \check Z^1(U_\bullet,
\underline Z)$
be a twisted 1-cocycle. We shall construct 
functions $\chi_{i}:U_{i}\to Z$
such that $(\delta_g\chi)_{..}=\varphi_{..}$.
$B$ is paracompact, hence without restriction we 
can assume that even the closed cover $\{\overline U_i\}$ is
locally finite and all $\varphi_{ij}$ are defined on the whole  of 
$\overline U_{ij}$.
Let
\begin{eqnarray*}
K:=\Big\{ (J,\chi_{.}) &|& J\subset I, {\rm for\ all\ }j\in J:
\chi_{j}:\overline U_{j}\to Z
,\\
&&{\rm such \ that\ for\ all\ }j,k\in J\ {\rm  and}\  u\in \overline U_{jk}:\\
&&\chi_{j}(u)\ \chi_{k}(u)^{-1}\cdot g_{kj}(u)=\varphi_{kj}(u)
\Big\}
\end{eqnarray*}
$K$ is non-empty, since $(\{i\}, \{\varphi_{ii}\})\in L$, for each $i\in I$.
We define a partial order on $K$
such that every chain has an upper bound. We let
$(J,\chi_{.})\le(J',\chi'_{.})$ if and only if 
$J\subset J'\subset I$ and $\chi_{j}=\chi'_{j},$
for all $j\in J$. 
By Zorn's lemma, let $(J,\chi_{.})$ denote a maximal element of $K$.

Assume $J\not= I$, so there is some $a\in I\backslash J$.
Let $R:=\bigcup_{j\in J}(\overline U_j\cap \overline U_{a})\subset \overline U_{a}$.
For $u\in R$ we define
$\tilde \chi_{a}(u):=
\varphi_{ja}(u)
 \chi_{j}(u)\cdot g_{ja}(u)$ if 
$u\in \overline U_j$.
This definition is independent of $j\in J$.
We end up with a diagram  
$$
\xymatrix{
R\ar[r]^{\tilde \chi_{a}}\ar[d]_\cap&Z\\
\overline U_{a}\ar@{.>}[ru]_{ \chi_{a}}&.
}
$$
Since our  cover  is locally finite, $R$ is closed, but 
$Z$ is contractible, therefore
an extension $\chi_{a}$ exists \cite[Lem. 4]{DD}.
This contradicts the maximality of $(J, \chi_.)$,
so $J=I$. 
\\\\
k=2:
Let
$
\varphi_{...}\in \check Z^2(U_\bullet,
\underline Z
)$
be a twisted 2-cocycle. We shall construct 
functions $\chi_{ij}:U_{ji}\to Z
$
such that $(\delta_g\chi)_{...}=\varphi_{...}$.
Again by paracompactness of $B$
we can assume that even the closed cover $\{\overline U_i\}$ is
locally finite and all $\varphi_{ijk}$ are defined on the whole  of 
$\overline U_{ijk}$.
Let
\begin{eqnarray*}
L:=\Big\{ (J,\chi_{..}) &|& J\subset I, {\rm for\ all\ }i,j\in J:
\chi_{ij}:\overline U_{ij}\to Z
,\\
&&{\rm such \ that\ for\ all\ }i,j,k\in J\ {\rm  and}\  u\in \overline U_{ijk}:\\
&&\chi_{ji}(u)\ \chi_{ki}(u)^{-1}\
\chi_{kj}(u)\cdot g_{ji}(u)=\varphi_{kji}(u)
\Big\}
\end{eqnarray*}
The set $L$ is non-empty, since for each $i\in I$
$(\{i\}, \{\varphi_{iii}\})\in L$.
We define a partial order on $L$
such that every chain has an upper bound. We let
$(J,\chi_{..})\le(J',\chi'_{..})$ if and only if 
$J\subset J'\subset I$ and $\chi_{ij}=\chi'_{ij},$
for all $i,j\in J$. 
By Zorn's lemma, let $(J,\chi_{..})$ denote a maximal element of $L$.
Assume $a\in I\backslash J$. Then
let 
\begin{eqnarray*}
M_{a}&:=&\Big\{ (K,\psi_{.a})\ |\ K\subset J, {\rm for\ all\ }k\in K:
\psi_{ka}:\overline U_{ka}\to Z
,\\
&&{\rm such \ that\ for\ all\ }k,l\in K {\rm \ and \ for\ } u\in \overline U_{lka}\\
&&\psi_{ka}(u)\ \psi_{la}(u)^{-1}\
\chi_{kl}(u)\cdot g_{la}(u)=\varphi_{lka}(u)
\Big\}
\end{eqnarray*}
The set $M_a$ is non-empty as for each 
$j\in J$ we find $(\{j\},\{ 1 \})\in M_a$. This because we 
always have 
$1 \cdot 1 \cdot\chi_{jj}(u)\cdot g_{ja}(u)=\varphi_{jjj}(u)\cdot g_{ja}(u)
=\varphi_{jja}(u)$.
In the same manner as before,
let
$(K,\psi_{.a})\le(K',\psi'_{.a})$ if and only if 
$K\subset K'\subset I$ and $\psi_{ka}=\psi'_{ka},$
for all $k\in K$. $\le$ is a partial order on $M_a$ and every chain has 
an upper bound. Let $(K,\psi_{.a})$ be a maximal element.
Assume $b\in J\backslash K$.
Let $S:=\bigcup_{k\in K}(\overline U_k\cap \overline U_{b}\cap\overline U_{a})\subset 
\overline U_{ba}$.
For $u\in R$ we define
$\tilde\psi_{ba}(u):=\varphi_{kba}(u)\chi_{kb}(u)\cdot g_{ba}(u)^{-1}
\psi_{ka}(u)$, if 
$u\in \overline U_{kba}, k\in K$.
By a one line calculation we find that
this definition is independent of $k\in K$.
So we have a diagram
$$
\xymatrix{
S\ar[r]^{
\tilde\psi_{ba}}\ar[d]_\cap&
Z
\\
\overline U_{ba}\ar@{.>}[ru]_{\psi_{ba}}&.
}
$$
$S$ is closed, since our  cover  is locally finite. 
Thus, since $Z
$ is contractible, 
there is an extension $\psi_{ba}$ \cite[Lem. 4]{DD}.
This contradicts the maximality of $(K,\psi_{.a})$,
so $K=J$. 
We define $\psi_{aa}:=\varphi_{aaa}$
and then 
$\psi_{aj}(u):=\varphi_{aja}(u)\cdot g_{aj}(u)
\psi_{aa}(u)\cdot g_{aj}(u)^{-1}\psi_{ja}(u)\cdot g_{aj}(u)$, for $ j\in J$.
We let $J':=\{a\}\cup J$
and extend $\chi_{..}$ to $J'$ by
$$
\chi'_{ij}:=
\begin{cases}
\chi_{ij},&{\rm if\ } i,j\in J,\\
\psi_{ia},& {\rm if\ } j=a, i\in J,\\
\psi_{aa},& {\rm if\ } i=j=a,\\
\psi_{aj},& {\rm if\ } i=a, j\in J 
\end{cases}.
$$
It is straight forward to check that $(J',\chi'_{..})\in L$,
and as clearly $(J,\chi_{..})\le ( J',\chi'_{..})$
we have a  contradiction.
Hence $J=I$, and the lemma is proven.
\end{pf}

\begin{thm}\label{ThmTheSpecialCaseOfRZ}
Let $G=\R^n, N=\Z^n$. Then every topological triple 
is almost strict, i.e. in this case ${\rm Top}^{\rm as}(B)={\rm Top}(B)$.
Moreover, the group $Q_{_\R}$ defined in diagram (\ref{DiagOfExatSequences}) vanishes, hence the three natural transformations of 
Theorem \ref{ThmTheMainTheoremOfThisThesisNumberThree} reduce to one single transformation
$$
\delta:{\rm Top}\to{\rm Dyn^\dag}.
$$
\end{thm}

\begin{pf}
By the lemma above,  
it is sufficient to show that
$$
Z:=Z^1_{\rm cont}(\R^n,{\rm Map}(\T^n,\U(1)))
$$
is contractible.
Then it follows that $\check H^2(B,\underline Z,g_{..})=\{1\}$ and
every topological triple is almost strict.
Because $Q_{_\R}=\check H^2(B,\underline Z,g_{..})=\{1\}$,
the torsor structures in Theorem \ref{ThmTheMainTheoremOfThisThesisNumberThree}
are trivial and we obtain a natural map
$\delta(B):{\rm Top}(B)\to{\rm Dyn^\dag}(B)$.

We now show that $Z$ is contractible.
For a cocycle  $\alpha\in Z$
we have 
$
\alpha(0)(z)=\alpha(0+0)(z)=\alpha(0)(z)\alpha(0)(z)
$
 for all $z\in \T^n$, so 
$\alpha(0)(z)=1$ and each $\alpha(g):\T^n\to \U(1)$ is null homotopic as 
$\R^n$ is path connected.
Thus
$$Z=Z^1_{\rm cont}(\R^n,{\rm Map}(\T^n,\U(1)))
=Z^1_{\rm cont}(\R^n,{\rm Map}_0(\T^n,\U(1))),
$$ for 
${\rm Map}_0(\T^n,\U(1)):=\{f\simeq {\rm const.}\}$.
But for each $\alpha:\R^n\to {\rm Map}_0(\T^n,\U(1))$
there exists a lift $\overline\alpha$ in
$$
\xymatrix{
&&{\rm Map}(\T^n,\R^n)\ar[d]\\
\R^n\ar[rr]^\alpha\ar[rru]^{\overline\alpha}&&{\rm Map}_0(\T^n,\U(1)),
}
$$
and because of $d\alpha=1$, we have $d\overline\alpha=m\in 
Z^2_{\rm cont}(\R^n,{\rm Map}(\T^n,\Z))\cong\Z$. Then
$\alpha'(g)(z):=\overline\alpha(g)(z)-m$ defines a unique 
element  $\alpha'\in Z^1_{\rm cont}(\R^n,{\rm Map}_0(\T^n,\R))$ such
that $\alpha'(g)(z)\Z=\alpha(g)(z)$.
In fact, the mapping $\alpha\mapsto \alpha'$ is a 
homoeomorphism 
$$
Z^1_{\rm cont}(\R^n,{\rm Map}_0(\T^n,\U(1)))\to Z^1_{\rm cont}(\R^n,{\rm Map}(\T^n,\R)),
$$
and the latter space is easily seen to be contractible by 
$h(t,\alpha')(g)(z):= t\ \alpha'(g)(z)$, $t\in [0,1]$.
This proves the theorem.
\end{pf}

\begin{cor}
Let $G=\R^n, N=Z^n$. Then a pair 
has an extension to a toplogical T-duality
triple if and only if it has an extension to 
a dualisable dynamical triple.
\end{cor}
\begin{pf}
This follows from Corollary \ref{CorTheCorOfTheMainTheorem}
and Theorem \ref{ThmTheSpecialCaseOfRZ}.
\end{pf}

So far we have seen that in the case 
of $G=\R^n$, $N=\Z^n$
\begin{equation}\label{EqSomeInclusions}
{\rm Top^{im}}(B)\subset
{\rm Top^{s}}(B)\subset
{\rm Top^{as}}(B)\stackrel{!}{=}
{\rm Top}(B),
\end{equation}
and the composition 
of the natural transformations $\tau$ and $\delta$
is the identity transformation
$\delta\circ\tau={\rm id}: {\rm Dyn^{\dag}}\to{\rm Dyn^\dag}.
$
The opposite composition
$
J:=\tau\circ\delta:{\rm Top}\to{\rm Top}
$ 
is an idempotent $J\circ J=J$.
Unfortunately, inside our local
theory 
we won't be able to answer the 
question whether $J$ is  the identity 
on ${\rm Top}$, i.e. all inclusions 
in  (\ref{EqSomeInclusions}) are
equalities and $\delta$ and $\tau$
are both equivalences of functors 
inverse to each other.
Nevertheless we can give this
appealing answer when we restrict 
the functors ${\rm Top}$ and ${\rm Dyn^\dag}$ to
the subcategory ${\rm CW}\subset \{{\rm bases\ spaces}\}$
of CW-complexes 
(or, more general, of base spaces with the homotopy type of 
a CW-complex).
The point is that we can use a main result of \cite{BS} about 
the representability of the functor 
${\rm Top}|_{\rm CW}={\tt Triple}_n$.
However, the statement therein is not completely 
precise. We restate it for convenience.

\begin{lem}[{\cite[Thm. 7.24]{BS}}]\label{ThmOfBSOnTheRepresentability}
The functor
${\rm Top}|_{CW}:{\rm CW}\to\{ {\rm sets} \}$  is representable 
by a space $R_n\in {\rm CW}$, i.e.
there is an equivalence of functors
$$
\Psi:[\ .\ ,R_n] \cong {\rm Top}|_{\rm CW}(\ .\ ),
$$
where the squared brackets denote
the homotopy classes of continuous maps.

The definition of $\Psi$  is via pullback
of a certain topological triple 
$x_{\rm univ}\in {\rm Top}(R_n)$ (the universal triple) over $R_n$;
for $B\in {\rm CW}, [f]\in[B,R_n]$ it is
$$
\Psi(B)([f]):=f^*x_{\rm univ}.
$$
\end{lem}

Further, it is shown in \cite[Sec. 4]{BS} that the space $R_n$
has an homotopy action of the
so-called T-duality  
group $O(n,n,\Z)$ which
 is the group of $2n\times2n$-matrices 
that fix  the form 
$$
\Z^{2n}\ni(a_1,\dots,a_n,b_1,\dots,b_n)\mapsto \sum a_ib_i\in\Z.
$$
So each element of
$O(n,n,\Z)$ defines a homotopy class
of maps $R_n\to R_n$.
In particular, the element
$$
\Bigg(
\begin{matrix}
 {\bf 0}_n &{\bf 1}_n\\
{\bf 1}_n& {\bf 0}_n
\end{matrix}
\Bigg)\in O(n,n,\Z)
$$
defines  (the homotopy class of) a function $T:R_n\to R_n$.
$T$ is constructed such that $T\circ T={\rm id}_{R_n}$
and such that 
the pullback
$T^*:{\rm Top}(R_n)\to {\rm Top}(R_n)$
exchanges the underlying torus bundles\footnote{
This is meaningful, as
the case of the groups $G=\R^n,N=\Z^n$ is  self-dual, i.e.
$G\cong\widehat G, N\cong N^\perp$.
}.
To be precise, the construction of $T$
is such that if $[(\kappa, (P,E),(\widehat P,\widehat E))]\in{\rm Top}(R_n)$
is a topological triple, 
then $T^*$ maps this triple to a triple 
$[(\kappa', (P', E'),(\widehat P',\widehat E'))]$,
such that
$E'\cong\widehat E$ and $\widehat E'\cong E$.
Note that
$$
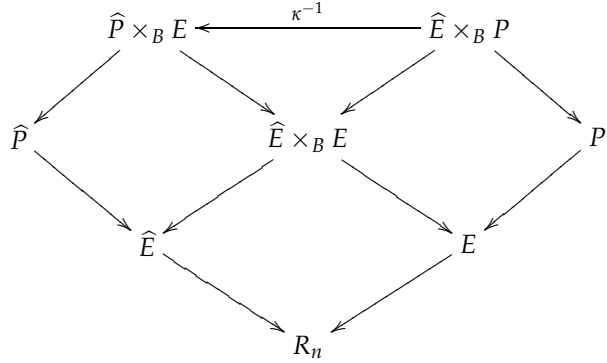
\xymatrix{
&\widehat P\times_B E\ar[dr]\ar[dl]&&\widehat E\times_B  P\ar[dl]\ar[dr]\ar[ll]_{\kappa^{-1}}&\\
\widehat P\ar[dr]&&\widehat E\times_B E\ar[dr]\ar[dl]&& P\ar[dl]\\
&\widehat E\ar[dr]&& E\ar[dl]&\\
&&R_n&&
}
$$ 
is not a topological triples as 
the inverse of $\kappa$ does not satisfy the Poincaré condition,
but also note that
$$
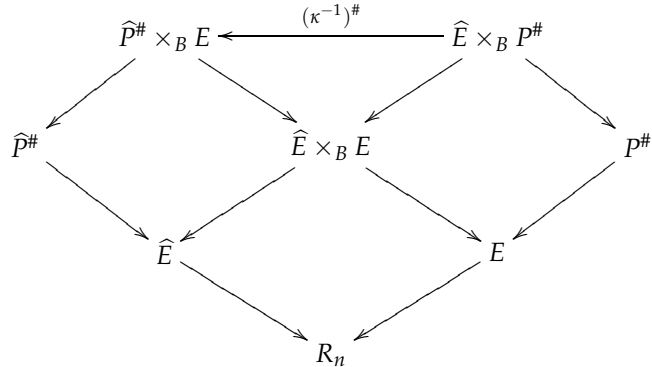
\xymatrix{
&\widehat P^\#\times_B E\ar[dr]\ar[dl]&&\widehat E\times_B  P^\# \ar[dl]\ar[dr]\ar[ll]_{(\kappa^{-1})^\#}&\\
{\widehat P}^{\#}\ar[dr]&&\widehat E\times_B E\ar[dr]\ar[dl]&& P^{\#}\ar[dl]\\
&\widehat E\ar[dr]&& E\ar[dl]&\\
&&R_n&&
}
$$ 
is a topological triple, wherein the superscript $\#$ denotes the 
complex conjugate\footnote{Complex conjugation may be defined by taking 
$\HH=l^2\N$ and then defining the complex conjugate bundles and isomorphism
by complex conjugation of the local transition functions and local isomorphisms.
} bundles and isomorphism.

We claim that
\begin{eqnarray}\label{EqDitsch}
T^*[(\kappa, (P,E),(\widehat P,\widehat E))]
=
[((\kappa^{-1})^\#, (\widehat P^\#,\widehat E),( P^\#,\widehat E))].
\end{eqnarray}
In fact,  by \cite[Prop. 7.4]{BS} 
the set of topological triples over $R_n$ with fixed torus bundles 
$E',\widehat E'$  is a torsor over $H^3(R_n,\Z)$ and
by \cite[Lem. 3.3]{BS} $H^3(R_n,\Z)=0$.
In other words, there only exists one triple over $R_n$ 
which has underlying torus bundles $E'\cong \widehat E$ and 
$\widehat E'\cong E$,
hence equation (\ref{EqDitsch}) is valid.

The properties of the map $T$ which is called universal T-duality 
enables us  to prove the statement indicated above:

\begin{thm}\label{Fischteichsegler}
Let $G=\R^n$ and $N=\Z^n$. Then
we have an equivalence of functors
$$
\delta|_{\rm CW}:{\rm Top}|_{\rm CW}\cong{\rm Dyn^\dag}|_{\rm CW}:\tau|_{\rm CW}.
$$
In particular
\begin{equation*}
{\rm Top^{im}}(B)=
{\rm Top^{s}}(B)=
{\rm Top^{as}}(B)=
{\rm Top}(B),
\end{equation*}
for all $B\in{\rm CW}$.
\end{thm}
\begin{pf}
It suffices to show that $J(B)={\rm id}$, for all $B\in{\rm CW}$.
$J$ is a natural transformation, so
$$
\xymatrix{
{\rm Top}(R_n)\ar[rr]^{J(R_n)}\ar[d]^{T^*}&&{\rm Top}(R_n)\ar[d]^{{T}^*}\\
{\rm Top}(R_n)\ar[rr]^{J(R_n)}&&{\rm Top}(R_n)\\
}
$$
is a commutative diagram.
But by construction $J(B)$ does not change the underlying pairs
of the corresponding topological triples, for any $B$, so 
by the commutativity of the diagram $J(R_n)$ at least does not
change the underlying pairs {\it and} underlying dual pairs of the corresponding triples.
That $J(R_n)$ is the identity, i.e. it does not change 
the equivalence class of the isomorphisms $\kappa$, follows 
from the trivial $H^3(R_n,\Z)$-torsor structure of 
the set of topological triples with fixed torus bundles.
Hence $J(R_n)={\rm id}$.

Now, let $B\in{\rm CW}$ and $x\in{\rm Top}(B)$ be any 
topological triple. By the universal property of $x_{\rm univ}$,
there is $f:B\to R_n$ such that $x=f^*x_{\rm univ}$ and, 
by naturality of $J$, 
$$
J(B)(x)=J(B)(f^*x_{\rm univ})= f^*(J(R_n)(x_{\rm univ}))=f^*(x_{\rm univ})=x.
$$
So $J(B)={\rm id}$, for all $B\in{\rm CW}$.
\end{pf}

We end this section with a remark on homotopic deckers.
Its content is based on the fact that the  functor 
${\rm Top}|_{\rm CW}\cong[\ .\ ,R_n]$ is homotopy invariant.

\begin{rem}
Let $(P,E)$ be a pair over $B\in {\rm CW}$, then 
$(P\times[0,1],E\times[0,1])$  is a pair over
$B\times[0,1]$. 
Let $(\rho,P\times[0,1],E\times[0,1])$ 
be a dualisable dynamical triple --
in other words $\rho$ is a homotopy of 
dualisable deckers $\rho_0$ and $\rho_1$ on $(P,E)$ --
then the classes $[(\rho_0,P,E)]$ and $[(\rho_1,P,E)]$
coincide.
\end{rem}
\begin{pf}
As the the inclusions $i_t:B\hookrightarrow B\times[0,1], t=0,1,$
are homotopic, we have an equality
\begin{eqnarray*}
\tau(B)([(\rho_0,P,E)])&=&
i_0^*\tau(B\times[0,1])([(\rho,P\times[0,1],E\times[0,1])])\\
&=&i_1^*\tau(B\times[0,1])([(\rho,P\times[0,1],E\times[0,1])])\\
&=&\tau(B)([(\rho_1,P,E)]),
\end{eqnarray*}
and the claim  follows from $\delta(B)\circ\tau(B)={\rm id}$.
\end{pf}

\subsection{The  Structure of the Associated $C^*$-Dynamical Systems}
\label{SecTheStructureOfTheAssociatedCStarDynSys}

Let $(\rho,P,E)$ be a dualisable dynamical triple and
let us denote by $F:=P\times_{\PU(\HH)}\K(\HH)$
the associated $C^*$-bundle. The decker $\rho$
induces a $G$-action on $F$ by $[x,K]\cdot g:= [\rho(x,g),K]$
$x\in P, K\in\K(\HH)$.
This action defines another action $\alpha^\rho$ of $G$ on the $C^*$-algebra of sections
$\Gamma(E,F)$ such that $(\Gamma(E,F),G,\alpha^\rho)$ becomes 
a $C^*$-dynamical system. This action is given by
\begin{eqnarray*}
({\alpha^\rho}_g s)(e):=s(e\cdot gN )\cdot(-g),
\quad e\in E, g\in G,s\in \Gamma(E,F).
\end{eqnarray*}
In the same manner we 
obtain a dual $C^*$-dynamical system
$(\Gamma(\widehat E,\widehat F),\widehat G,\alpha^{\hat\rho})$
for the associated \label{PageOfCStarThings}
$C^*$-bundle
 $\widehat F:=\widehat P\times_{\PU(L^2(G/N)\otimes\HH)}\K(L^2(G/N)\otimes\HH)$
of  the dual triple $(\hat\rho,\widehat E,\widehat P)$ of $(\rho, E,P)$.  
The essence of this section is that we can establish  
an isomorphism 
of $C^*$-dynamical systems from the crossed product\footnote{See 
the section on crossed products on page \pageref{SecTheSecOnCrossedProducts}
for notation.}
of the first to the dual $C^*$-dynamical system  (Thm. \ref{ThmTDualityGenaralCase})
$$
\big(G\times_{\alpha^\rho}\Gamma(E,F),\widehat G,\widehat{\alpha^\rho}\big)
\cong \big(\Gamma(\widehat E,\widehat F),\widehat G,\alpha^{\hat \rho}\big).
$$
We are going to calculate the 
crossed product under a series of 
isomorphisms which again  will be a local
calculation. We start with the  description of the simplest case,
namely $B$ being a point.
\\\\
In the situation of the trivial pair over the point $B=\{*\}$ the
sections can be identified with the continuous functions 
$$\Gamma({\rm triv.\ pair})\cong C(G/N,\K(\HH)),$$
since any section $s:G/N\to G/N\times\K(\HH)$ is uniquely given by a function
$f$ such that $s(z)=(z,f(z)).$
Then the action of $g\in G$ on such a function $f$ is 
obtained from
\begin{eqnarray}
({\alpha^\rho}_g s)(z)f
&=&s(z+ gN )\cdot (-g)\nonumber\\
&=&(z+gN,f(z+gN)\cdot(-g) \nonumber\\
&=&(z, \mu(-g,z+gN)(f(z+gN)))\nonumber\\
&=&(z, (\mu(g,z))^{-1}(f(z+gN)))\nonumber\\
&=:&(z,({\alpha^\mu}_g f)(z)) \label{EqActOnFct}
 \end{eqnarray}
for the 1-cocycle $\mu$ determined by $\rho$.

For the $C^*$-algebra of the dual trivial pair over the point we have
$$\Gamma({\rm dual\ triv.\ pair})\cong 
C(\widehat G/N^\perp,\K(L^2(N^\perp,\HH))).$$
It will be convenient at some point to deal with the 
Hilbert space $L^2(\widehat{G/N},\HH)$ rather than 
with $L^2(G/N)\otimes\HH$.
So we have to transform \label{PageMoerFoureirere}
the  cocycle $\hat\mu$ of the dual decker $\hat\rho$ from
equation (\ref{EqTheDualDeckerCocycle})
by Fourier transform $\mathcal F:L^2(G/N)\to L^2(\widehat{G/N})$,
\begin{eqnarray*}
(\mathcal F_*\hat\mu)(\chi,\hat z):=(\mathcal F\otimes\Eins_\HH)\circ
\hat\mu(\chi,\hat z) \circ(\mathcal F^{-1}\otimes\Eins_\HH),
\end{eqnarray*}
and as in eq. (\ref{EqActOnFct}) this 
 gives us an action 
 $$
\alpha^{\mathcal F_*\hat\mu}:\widehat G\to {\rm Aut}(C(\widehat G/N^\perp,\K(L^2(N^\perp,\HH))))
 $$
 by $({\alpha^{\mathcal F_*\hat\mu}}_\chi f)(\hat z):=(\mathcal F_*\hat\mu)(\chi)(\hat z)^{-1}(f(\hat z+\chi N^\perp))$,
 for $\chi\in \widehat G, \hat z\in\widehat G/N^\perp.$ 
 
The next lemma is a simple link between $\mathcal F_*\hat\mu$ and 
the left regular representation  $\lambda_{_{\widehat{G/N}}}$
on $L^2(\widehat{G/N})$.
$\lambda_{_{\widehat{G/N}}}\in {\rm Hom}(\widehat{G/N},\U(L^2(\widehat{G/N})))$
is a continuous homomorphism.

\begin{lem}\label{LemExtOfShiftOp}
 There exists a continuous extension 
$\Lambda\in {\rm Hom}(\widehat G,\U(L^2(\widehat{G/N})))$
such that
$$
\xymatrix{
\widehat{G/N}\ar[r]^{\!\!\!\!\!\!\!\!\!\!\!\!\!\!\!\!\!\!\lambda_{{\widehat{G/N}}}}\ar[d]^\cap&\U(L^2(\widehat{G/N}))\\
\widehat G\ar[ru]_{\!\!\!\!\! \Lambda}&
}
$$
commutes.
\end{lem}
\begin{pf}
Let $\beta\in  \widehat {G/N}$.
The Fourier transform $\mathcal F$ turns $\lambda_{_{\widehat{G/N}}}(\beta)$  
into the multiplication operator
$\langle -\beta,\_\rangle=\mathcal F^{-1} \circ \lambda^\perp(\beta) \circ \mathcal F 
\in \U(L^2(G/N))$.
Let $\sigma:G/N\to G$ be our chosen Borel section.
 We define an extension of
$\beta\mapsto \langle -\beta,\_\rangle$ by
$\widehat G\ni\chi\mapsto\langle \chi,-\sigma(\_)\rangle$.
This is a strongly continuous homomorphism, and it follows that 
$\Lambda(\chi):= \mathcal F\circ \langle \chi,-\sigma(\_)\rangle \circ \mathcal F^{-1}$ 
gives us the desired homomorphism $\Lambda:\widehat G\to \U(L^2(\widehat{G/N})).$
\end{pf}

With the notation of this lemma and equation (\ref{EqTheDualDeckerCocycle})
we can reformulate the defintion of $\mathcal F_*\hat\mu$, we have
$(\mathcal F_*\hat\mu)(\chi,\hat z)={\rm Ad}(\Lambda(\chi)\otimes\Eins_\HH)$
which is the adequate formula for the next theorem.
Therein $\alpha\mapsto\alpha^\perp$ denotes the natural
isomorphism $\widehat{G/N}\cong N^\perp\subset \widehat G$ which
is defined by $\langle\alpha^\perp,g\rangle:=\langle\alpha,gN\rangle$.

\begin{thm}\label{ThmTDualityOverThePoint}
Assume the considered  dynamical triple over $B=\{*\}$ 
is dualisable. Then there is an isomorphism of $C^*$-dynamical systems
$$\Big(G\times_{\alpha^\mu} C(G/N,\K(\HH)),\widehat G,\widehat{\alpha^\mu}\Big)
\stackrel{\cong}{\longrightarrow}
\Big(C(\widehat{G}/N^\perp,\K(L^2(\widehat{G/N},\HH))),\widehat G,\alpha^{\mathcal F_*\hat\mu}\Big).$$
\end{thm}

\begin{pf}
First note
that naturally
$$
G\times_{\alpha^\mu} C(G/N,\K(\HH))\cong
\overline{C_c(G\times G/N,\K(\HH))}^{\|.\|}\subset\LL(L^2(G\times G/N,\HH)),
$$
wherein $f\in C_c(G\times G/N,\K(\HH))$ acts
on $F\in L^2(G\times G/N,\HH)$ by
$$
(f\times F)(g,z)=\int_G (\mu(-g)(z))^{-1}(f(h,z-gN))\ F(g-h,z)\ dh.
$$
We assume the dynamical triple to be dualisable, so  we can lift $\mu$ to a unitary (Borel) cocycle
$\overline\mu:G\times G/N\to \U(\HH)$.
Note that for $K\in\K(\HH))$
$$\overline\mu(-g,z)\mu(-g,z)^{-1}(K)
=K\overline\mu(-g,z):\HH\to\HH.$$
We now define 
a unitary isomorphism
$$u:L^2(G\times G/N,\HH)\to L^2(\widehat G,L^2( \widehat{G/N},\HH))$$
by the composition 
$u:= {\rm shift}\circ ({\rm Fourier\ trans.})\circ {\rm mult}_{\overline\mu(^\_\_,\_)},$
explicitely 
\begin{eqnarray*}
(u F)(\chi)(\alpha)&=&\widehat{(\overline\mu(^\_ \_,\_) F)} (\chi+\alpha^\perp,-\alpha)\\
&=&\int_{G\times G/N} \langle \chi+\alpha^\perp,g\rangle\ \langle-\alpha,z\rangle\ 
\overline\mu(-g,z)\ F(g,z)\ d(g,z),
\end{eqnarray*}
for $\chi\in\widehat G, \alpha\in \widehat{G/N}$.
The next step is to calculate $u (f\times\_ ) u^{-1}.$
This is straightforward, but to keep the calculation
readable we first introduce some short hands:
$_{\overline\mu} F(g,z):=\overline\mu(-g,z)F(g,z),$
$f_{\overline\mu}(g,z):=f(g,z)\overline\mu(g,z)^{-1}$
and 
$\hat f^{{^{2}}}$
 the Fourier transform in the second, the $G/N$ variable only.
 Now
\begin{eqnarray}
&&u(f\times F)(\chi)(\alpha)\nonumber\\
&=&\int_{G\times G/N} \int_G\langle\chi+\alpha^\perp,g\rangle\ \langle-\alpha,z\rangle\ 
f(h,z-gN)\overline\mu(-g,z)\ F(g-h,z)\ dh\ d(g,z)\nonumber\\
&=&\int_{G\times G/N} \int_G\langle\chi+\alpha^\perp,g\rangle\ \langle-\alpha,z\rangle\
f_{\overline\mu}(h,z-gN)\ _{\overline\mu}F(g-h,z)\ dh\ d(g,z)\nonumber\\
&=&\int_{G\times G}\int_{\widehat{G/N}}\langle\chi+\alpha^\perp,g\rangle\
\widehat{f_{\overline\mu}}^{_{2}}(h,\beta)\ \langle\beta,gN\rangle\ 
\widehat{_{\overline\mu}F}^{_{2}}(g-h,(-\alpha)-\beta)\
d\beta\ d(h,g)\nonumber\\
&=&
\int_{\widehat{G/N}}\int_G \widehat{f_{\overline\mu}}^{_{2}}(h,\beta)\
\langle\chi+\alpha^\perp+\beta^\perp,h\rangle\
\widehat{_{\overline\mu}F}(\chi+\alpha^\perp+\beta^\perp,(-\alpha)-\beta)\ dh\ d\beta\nonumber\\
&=&\int_{\widehat{G/N}}\widehat{f_{\overline\mu}}(\chi+\alpha^\perp+\beta^\perp,\beta)\
\widehat{_{\overline\mu}F}(\chi+\alpha^\perp+\beta^\perp,(-\alpha)-\beta)\  d\beta\nonumber\\
&=&\int_{\widehat{G/N}}
\underbrace{\widehat{f_{\overline\mu}}(\chi+\gamma^\perp,\gamma-\alpha)}\
\underbrace{\widehat{_{\overline\mu}F}(\chi+\gamma^\perp,-\gamma)}\  
d\gamma,\label{EqTildeF}\\
&&\qquad\qquad\quad\ \ =: f^{\overline\mu}(\chi)(\alpha,\gamma)\qquad=u F(\chi)(\gamma)\nonumber
\end{eqnarray}
and we think of $ f^{\overline\mu}$ as a continuous family of Hilbert-Schmidt operators
$f^{\overline\mu}:\widehat G\to \K(L^2(\widehat{G/N},\HH))$, i.e. we
 do not decide in notation between the operator $f^{\overline\mu}(\chi)$ 
 and its integral kernel.
From the definition of the kernel  $f^{\overline\mu}$ we obtain 
$f^{\overline\mu}(\chi+\beta^\perp)(\alpha,\gamma)=f^{\overline\mu}(\chi)(\alpha+\beta,\gamma+\beta)$,
$\alpha,\gamma,\beta\in \widehat{G/N}$,
so
the operator $f^{\overline\mu}(\chi)$ satisfies the identity
$$f^{\overline\mu}(\chi+\beta^\perp)= (\lambda_{_{\widehat {G/N}}}(\beta)^{-1}\otimes\Eins)
f^{\overline\mu}(\chi) (\lambda_{_{\widehat {G/N}}}(\beta)\otimes\Eins)\in\K(L^2(\widehat{G/N},\HH)),$$
for the left regular representation $\lambda_{_{\widehat {G/N}}}$.
We now use the chosen extension  $\Lambda\in{\rm Hom}(\widehat G,\PU(L^2(\widehat{G/N},\HH)))$
from Lemma \ref{LemExtOfShiftOp} to
 define  
$T_{\overline\mu}:C_c(G\times G/N,\K(\HH)) \to C(\widehat G/N^\perp,\K(L^2(\widehat{G/N},\HH)))$
by 
\begin{eqnarray}\label{EqTheTOperator}
(T_{\overline\mu}f)(\chi N^\perp)&:=& 
(\Lambda(\chi)\otimes\Eins)\ f^{\overline\mu}(\chi)\ (\Lambda(\chi)^{-1}\otimes\Eins)\nonumber\\
&=&{\rm Ad}(\Lambda(\chi)\otimes\Eins))  f^{\overline\mu}(\chi).
\end{eqnarray}
It
is now a lengthy but straight forward calculation  to check that, firstly, 
$T_{\overline\mu}$ commutes with
the $^*$-operation, i.e.
$T_{\overline\mu}(f^\times)(\hat z)=(T_{\overline\mu}f)(\hat z)^*$, for
$\hat z \in \widehat G/N^\perp$ and
$f^\times(g,z)=\mu(g)(z)^{-1}(f(-g,z+gN)^*)$.
Secondly,  $T_{\overline\mu}$ preserves the product, i.e.
$T_{\overline\mu}(f_1\times f_2)(\hat z)=T_{\overline\mu}f_1(\hat z) T_{\overline\mu}f_2(\hat z)$.
Thirdly, $T_{\overline\mu}$ preserves the norm, i.e.
 $\| f\times\!\_\|=\sup_{\hat z\in\widehat G/N^\perp} \| T_{\overline\mu}f(\hat z)\|$.

Moreover $T_{\overline\mu}$ has dense image,
so it extends uniquely to a $C^*$-algebra isomorphism
$$T_{\overline\mu}: G\times_{\alpha^\mu} C(G/N,\K(\HH))\to C(\widehat G/N^\perp,\K(L^2(\widehat{G/N},\HH))).$$
 It remains to show that
 $T_{\overline\mu}(\widehat{\alpha^\mu}_\chi f)={\alpha^{\mathcal F_*\hat\mu}}_\chi(T_{\overline\mu}f)$:
 
 By definition $\widehat{\alpha^\mu}_\chi$ is just the multiplication with
 the character  $\chi$, so the Fourier transform gives us simply 
 a shift in the argument by $\chi$. We get
\begin{eqnarray}\label{EqEstablishDualDeckerOverPoint}
T_{\overline\mu}(\widehat{\alpha^\mu}_\chi f)(\chi'N^\perp)&=&{\rm Ad}(\Lambda(\chi')\otimes\Eins)  f^{\overline\mu}(\chi'+\chi)\nonumber\\
&=&{\rm Ad}(\Lambda(\chi)^{-1}\otimes\Eins) {\rm Ad}(\Lambda(\chi'+\chi)\otimes\Eins)
(f^{\overline\mu}(\chi'+\chi))\nonumber\\
&=&(\mathcal F_*\hat\mu)(\chi)(\chi'N^\perp)^{-1} T_{\overline\mu}f(\chi'N^\perp+\chi N^\perp)\nonumber\\
&=&{\alpha^{\mathcal F_*\hat\mu}}_\chi(T_{\overline\mu}f)(\chi'N^\perp).
\end{eqnarray}
\end{pf}


The discussion in the previous theorem will serve as a description
of the local situation of the general case.
Let $(\rho,P,E)$ be a dualisable  dynamical triple, and
let $(\hat \rho,\widehat P,\widehat E)$ be its dual.
We defined the associated $C^*$-bundles $F,\widehat F$ above.
Recall the definition of $\hat g_{ji}$ out of $\phi_{ji}$ (p. \pageref{PageThePageOfHatGji}).
A priori there need not exist a continuous lift $\varphi_{ji}$ 
in the diagram
$$
\xymatrix{
&&\widehat G\ar[d]\\
U_{ji}\ar[rru]^{\varphi_{ji}}\ar[rr]_{\hat g_{ji}}&&\widehat G/N^\perp,
}
$$
but by Lemma \ref{LemLiftForLocTrivFB} we assume without 
restriction that our atlas is sufficiently refined such that $\varphi_{ji}$ exists.
We define 
$\varphi'_{ji}(u)(gN):=\phi_{ji}(u)(g,0)\ \langle\varphi_{ji}(u),g\rangle^{-1}$.
Although the function $u\mapsto\varphi_{ji}'(u)\in L^\infty(G/N,\U(1))$ 
need {\sl not} to be continuous, the function
$u\mapsto{\rm Ad}(\varphi_{ji}'(u))\in \PP L^\infty(G/N,\U(1))$
is continuous by Lemma \ref{LemATechLemmaForMu}, and we 
have the identity 
\begin{eqnarray*}
\phi_{nm}(u)(g,hN)&=&\phi_{nm}(u)(g+h,0)\ \phi_{nm}(u)(h,0)^{-1}\\
&=&\varphi_{nm}(u)(g+h)\ \varphi'_{nm}(u)(gN+hN)\ \varphi_{nm}(u)(h)^{-1}\ \varphi'_{nm}(u)(hN)^{-1}\\
&=&\varphi_{nm}(u)(g)\ d(\varphi'_{nm}(u))(g,hN).
\end{eqnarray*}

\begin{thm}\label{ThmTDualityGenaralCase}
There 
is an isomorphism of $C^*$-dynamical systems
$$\Big(G\times_{\alpha^\rho}\Gamma(E,F),\widehat G,\widehat{\alpha^\rho}\Big)
\stackrel{\cong}{\longrightarrow}
\Big(\Gamma(\widehat E,\widehat F),\widehat G,\alpha^{\hat\rho}\Big).
$$
\end{thm}

\begin{pf}
We generalise the proof of Theorem \ref{ThmTDualityOverThePoint}.
The sections $s\in\Gamma(E,F)$ are in 
one to one correspondence with families of functions
$s_i\in C(U_i\times G/N,\K(\HH)),$ $i\in I,$ which
satisfy
\begin{equation}\label{EqLocStrOfSections}
s_i(u,z)=\zeta_{ji}(u)(z)^{-1}(s_j(u,g_{ji}(u)+z)),\quad u\in U_{ij},z\in G/N.
\end{equation}
It is 
$
G\times_{\alpha^\phi}\Gamma(E,F)=\overline{C_c(G,\Gamma(E,F))}^{\|.\|},
$
and for each chart $U_i$ of the pair $(P,E)$ we have induced restriction maps
\begin{eqnarray*}
&G\times_{\alpha^\rho}\Gamma(E,F)\to G\times_{\alpha^{\mu_i}}C(U_i\times G/N,\K(\HH))&\\
&\cup\qquad\qquad\qquad\qquad\cup\qquad\qquad&\\
&C_c(G,\Gamma(E,F))\to C(U_i, C_c(G\times G/N,\K(\HH))).&\\
&f\mapsto f_i\qquad\ \qquad\qquad &
\end{eqnarray*}
The triple $(\rho, P,E)$ is dualisable,
so we can lift (if necessary after a refinement of the atlas)
the cocycles $\mu_i$ to
unitary Borel cocycles
$\overline\mu_i:U_i\to Z^1_{\rm Bor}(G,L^\infty(G/N,\U(\HH))).$
For each such $\overline\mu_i$ we define
an operator
$$
T_i:C(U_i, C_c(G\times G/N,\K(\HH)))\to
C(U_i\times \widehat G/N^\perp,\K(L^2(\widehat{G/N},\HH)))
$$
as in equation (\ref{EqTheTOperator}) by
$$
T_i f_i(u,\chi N^\perp):= (T_{\overline\mu_i(u)}f_i(u))(\chi N^\perp)={\rm Ad}( \Lambda(\chi)) 
f_i(u)^{\overline\mu_i(u)}(\chi),
$$
wherein we used the notation of equation (\ref{EqTildeF}).
In view of equation (\ref{EqLocStrOfSections}) we are
interested in the relation between $T_if_i$ and $T_jf_j$ on
overlaps $U_{ij}.$
So let $u\in U_{ji}$, then  
\begin{eqnarray*}
&&{f_i(u)^{\overline\mu_i(u)}}(\chi)(\alpha,\gamma)\\
&\stackrel{(\ref{EqTildeF})}{=}&\int_{G\times G/N} f_i(u)(g,z)\ \overline\mu_i(u)(g,z)^{-1}\ 
\langle\chi+\gamma^\perp,g\rangle\ \langle\gamma-\alpha,z\rangle\ d(g,z)\\
&\stackrel{(\ref{EqLocStrOfSections})\textrm{ and def. of }\phi_{ji}}{=}&
\int_{G\times G/N} \overline\zeta_{ji}(u)(z)^*\ f_j(u)(g,g_{ji}(u)+z)\
\overline\mu_j(u)(g,g_{ji}(u)+z)^{-1}\\
&&\qquad\qquad\quad \overline\zeta_{ji}(u)(gN+z)\ \phi_{ji}(u)(g,z)\
\langle\chi+\gamma^\perp,g\rangle\ \langle\gamma-\alpha,z\rangle\ d(g,z)\\\\
&\stackrel{}{=}&
\int_{G\times G/N} \overline\zeta^\varphi_{ji}(u)(z)^*\
f_j(u)_{\overline\mu_j(u)}(g,g_{ji}(u)+z)\\
&&\qquad\qquad\quad \overline\zeta^\varphi_{ji}(u)(gN+z)\ 
\langle\varphi_{ji}(u)+\chi+\gamma^\perp,g\rangle\ 
\langle\gamma-\alpha,z\rangle\ d(g,z)\\
&=:&\clubsuit.
\end{eqnarray*}
The notation
$f_j(u)_{\overline\mu_j(u)}$ is as in the previous theorem,
and $\overline\zeta^\varphi_{ji}(u)(z):=\overline\zeta_{ji}(u)(z)\varphi_{ji}'(u)(z)$.

Since $G/N$ is compact and $\overline\zeta_{ji}(u)(z)\in \U(\HH)$
it follows that $\overline\zeta^\varphi_{ji}(u)(\_)v\in L^2(G/N,\HH)$,
for each $v\in\HH$. So it possesses a Fourier decomposition 
\begin{eqnarray*}
\overline\zeta^\varphi_{ji}(u)(gN+z)&=&
\int_{\widehat{G/N}} \widehat{\overline\zeta^\varphi_{ji}(u)}(\delta)\ 
\langle-\delta,z\rangle\ \langle-\delta^\perp,g\rangle\ d\delta
\end{eqnarray*}
and
$$
\overline\zeta^\varphi_{ji}(u)(z)^*=\int_{\widehat{G/N}} 
\widehat{\overline\zeta^\varphi_{ji}(u)}(\varepsilon)^*\ 
\langle\varepsilon,z\rangle\  d\varepsilon,
$$
for some $ \widehat{\overline\zeta^\varphi_{ji}(u)}(\_):\widehat{G/N}\to\LL(\HH)$.
The calculation continues inserting this
into the previous equation,
\begin{eqnarray*}
\clubsuit&=&\int_{\widehat{G/N}}\int_{\widehat{G/N}}
\widehat{\overline\zeta^\varphi_{ji}(u)}(\varepsilon)^*\ 
\widehat{f_j(u)_{\overline\mu_j(u)}}
(\varphi_{ji}(u)+\chi+\gamma^\perp-\delta^\perp,\gamma-\delta-\alpha+\varepsilon)\\
&&\quad\quad\quad\quad\quad\quad\quad\quad
\langle\gamma-\delta-\alpha+\varepsilon,-g_{ji}(u)\rangle\
\widehat{\overline\zeta^\varphi_{ij}(u)}(\delta)\ d\varepsilon\ d\delta\\\\
&=&\int_{\widehat{G/N}}\int_{\widehat{G/N}}
\widehat{\overline\zeta^\varphi_{ji}(u)}(\varepsilon)^*\ 
{f_j(u)}^{\overline\mu_j(u)}
(\varphi_{ji}(u)+\chi)(\alpha-\varepsilon,\gamma-\delta)\\
&&\quad\quad\quad\quad\quad\quad\quad\quad
\langle\gamma-\delta-\alpha+\varepsilon,-g_{ji}(u)\rangle\
\widehat{\overline\zeta^\varphi_{ij}(u)}(\delta)\ d\varepsilon\ d\delta\\\\
&=&\int_{\widehat{G/N}}\int_{\widehat{G/N}}
\widehat{\overline\zeta^\varphi_{ji}(u)}(\alpha-\varepsilon')^*\ 
\langle-\varepsilon',-g_{ji}(u)\rangle\ {f_j(u)}^{\overline\mu_j(u)}
(\varphi_{ji}(u)+\chi)(\varepsilon',\delta')\\
&&\quad\quad\quad\quad\quad\quad\quad\quad
\langle\delta',-g_{ji}(u)\rangle\
\widehat{\overline\zeta^\varphi_{ij}(u)}(\gamma-\delta')\ d\varepsilon'\ d\delta'.
\end{eqnarray*}
Now, bearing in mind that ${f_i(u)}^{\overline\mu_i(u)}(\chi)(\alpha,\gamma)$
is the kernel of the integral operator 
${f_i(u)}^{\overline\mu_i(u)}(\chi)\in\K(L^2(\widehat{G/N},\HH))$, 
the last expression becomes 
\begin{eqnarray*}
{f_i(u)}^{\overline\mu_i(u)}(\chi)
&=&
\eta_{ji}(u)^{*}\ {f_j(u)}^{\overline\mu_i(u)}(\varphi_{ji}(u)+\chi)\ \eta_{ji}(u)\\
&=&{\rm Ad}(\eta_{ji}(u)^*)({f_j(u)}^{\overline\mu_j(u)}(\varphi_{ji}(u)+\chi))
\in\K(L^2(\widehat{G/N},\HH)),
\end{eqnarray*}
wherein $\eta_{ji}(u)\in \U(L^2(N^\perp,\HH))$ is the 
composition of the two operators
\begin{eqnarray}
L^2(\widehat{G/N},\HH)&\to&L^2(\widehat{G/N},\HH)\label{EqTheConvOpera}\\
F&\mapsto&\widehat{\overline\zeta^\varphi_{ji}(u)}(-\_)*F:
\beta\mapsto\int_{\widehat{G/N}} 
\widehat{\overline\zeta^\varphi_{ji}(u)}(\beta-\alpha)\ F(\alpha^\perp)\ d\alpha\nonumber\\
{\rm and}\nonumber\\
\nonumber\\ 
L^2(\widehat{G/N},\HH)&\to&L^2(\widehat{G/N},\HH)\label{EqTheMultiOpera}\\
F&\mapsto& \Big(\beta\mapsto \langle\beta,-g_{ji}(u)\rangle\ F(\beta)\Big).\nonumber
\end{eqnarray}
Note that both of these operators are in fact unitary, so
$\eta_{ji}(u)$ is.
The calculation done so far can now 
give us the relation we are looking for
\begin{eqnarray}\label{EqTheDualSectLoc}
T_if_i(u,\chi N^\perp)&=&{\rm Ad}(\Lambda(\chi)\otimes\Eins){f_i(u)}^{\overline\mu_i(u)}(\chi)\nonumber\\
&=&{\rm Ad}(\Lambda(\chi)\otimes\Eins)\ {\rm Ad}(\eta_{ji}(u)^{-1})\  {f_j(u)}^{\overline\mu_j(u)}(\varphi_{ji}(u)+\chi)\nonumber\\
&=&\underbrace{{\rm Ad}\big((\Lambda(\chi)\otimes\Eins) \eta_{ji}(u)^{-1} 
(\Lambda(\varphi_{ji}(u)+\chi)^{-1}\otimes\Eins)\big)}\nonumber\\\
&&\qquad\qquad\qquad\qquad\quad\quad
=:(\mathcal F_*\hat\zeta_{ji})(u)(\chi N^\perp)^{-1}\nonumber\\
&&
{\rm Ad}(\Lambda(\varphi_{ji}(u)+\chi)\otimes\Eins)\ {f_j(u)}^{\overline\mu_j(u)}(\varphi_{ji}(u)+\chi)\nonumber\\
&=&(\mathcal F_*\hat\zeta_{ji})(u)(\chi N^\perp)^{-1}\ T_j f_j(u,
\underbrace{\varphi_{ji}(u)N^\perp}+\chi N^\perp).\\
&&\qquad\qquad\qquad\qquad\qquad\qquad\qquad =\hat g_{ji}(u) \in \widehat G/N^\perp\nonumber
\end{eqnarray}
$(\mathcal F_*\hat\zeta_{ji})(u)(\chi N^\perp)$ is in fact well defined for 
$\chi N^\perp\in\widehat G/N^\perp$, since
$\Lambda(\chi +\gamma^\perp)=\Lambda(\chi)\Lambda(\gamma^\perp)=
\Lambda(\chi) \lambda_{_{\widehat{G/N}}}(\gamma),$
and the left regular representation operator $\lambda_{_{\widehat{G/N}}}(\gamma)$ 
commutes with the convolution operator (\ref{EqTheConvOpera})
and commutes with the multiplication operator (\ref{EqTheMultiOpera})
up to $\langle\gamma,-g_{ji}(u)\rangle$, a $\U(1)$-valued multiple of the identity. Hence
the commutator of ${\rm Ad}(\Lambda(\gamma^\perp)\otimes\Eins)$ and 
${\rm Ad}(\eta_{ji}(u))$ vanishes in $\PU(L^2(\widehat{G/N},\HH))$.

In view of equation (\ref{EqLocStrOfSections}),
equation (\ref{EqTheDualSectLoc}) shows that 
the family $\{T_if_i\}_{i\in I}$ defines a section 
in a $\K(L^2(\widehat{G/N},\HH))$-bundle over $E$ with
transition functions $\mathcal F_*\hat\zeta_{ji}$.
Up to Fourier transform, 
this bundle is nothing but $\widehat F$ itself, for we find
\begin{eqnarray*}
&&
(\mathcal F_*\hat\zeta_{ji})(u)(\hat z)\\
&=&
{\rm Ad}\Big(
(\mathcal F\otimes\Eins)\circ
(\langle\varphi_{ji}(u)+\hat\sigma(\hat z),-\sigma(\_)\rangle\otimes\Eins)\circ
(\mathcal F^{-1}\otimes\Eins)\circ
\eta_{ji}(u)(\hat z)
\circ\\
&&(\mathcal F\otimes\Eins)\circ
(\langle \hat\sigma(\hat z),\sigma(\_)\rangle\otimes\Eins)\circ
(\mathcal F^{-1}\otimes\Eins)\Big)\\
&=&
{\rm Ad}\Big(
(\mathcal F\otimes\Eins)\circ
(\langle\varphi_{ji}(u)+\hat\sigma(\hat z),-\sigma(\_)\rangle\otimes\Eins)\\
&&(\lambda_{_{G/N}}(-g_{ji}(u))\otimes\Eins)\ \overline\zeta_{ji}(u)(-\_)\
(\varphi_{ji}'(u)(-\_)\otimes\Eins)\\
&&(\langle \hat\sigma(\hat z),\sigma(\_)\rangle\otimes\Eins)\circ
(\mathcal F^{-1}\otimes\Eins)\Big)\\
&=&{\rm Ad}\Big(
(\mathcal F\otimes\Eins)\circ
(\langle\varphi_{ji}(u)+\hat\sigma(\hat z),\sigma(\_+g_{ji}(u)-\sigma(\_)\rangle\otimes\Eins)\\
&&(\lambda_{_{G/N}}(-g_{ji}(u))\otimes\Eins)\ \overline\zeta_{ji}(u)(-\_)\
(\langle\varphi_{ji}(u),-\sigma(\_)\rangle\ \varphi_{ji}'(u)(-\_)\otimes\Eins)\\
&&\circ
(\mathcal F^{-1}\otimes\Eins)\Big)\\
&=&{\rm Ad}\Big(
(\mathcal F\otimes\Eins)\circ
(\langle\hat\sigma(\hat g_{ji}(u)+\hat z),-\sigma(\_)\rangle\otimes\Eins)\\
&&(\lambda_{_{G/N}}(-g_{ji}(u))\otimes\Eins)\ \overline\zeta_{ji}(u)(-\_)\
(\phi_{ji}(u)(-\sigma(\_),0)\otimes\Eins)
\circ(\mathcal F^{-1}\otimes\Eins)\Big)\\
&=&
(\mathcal F\otimes\Eins)\circ
 {\hat\zeta}_{ji}(u)(\hat z)\
\circ(\mathcal F^{-1}\otimes\Eins).
\end{eqnarray*}
So the the family $\{T_if_i\}$ defines a section
$Tf\in\Gamma(\widehat E,\widehat F)$, and
we have constructed a map
$
T:C_c(G,\Gamma(E,F))\to \Gamma(\widehat E,\widehat F)
$
which extends to an isomorphism of $C^*$-algebras
$$
T:G\times_{\alpha^\rho}\Gamma(E,F))\to \Gamma(\widehat E,\widehat F)
$$
Then the relation $T({\widehat{\alpha^\rho}}_\chi f)={\alpha^{\hat\rho}}_\chi (Tf)$
is established by local calculation in the same manner as over
the point in equation (\ref{EqEstablishDualDeckerOverPoint}).
Thus $T$ is in fact an isomorphism of $C^*$-dynamical systems.
\end{pf}

\newpage
\begin{appendix}

\section{Some Notation and Basic Lemmata}

\subsection{Groups}\label{SuperSectonGoops}
Let  $G$ be a Hausdorff locally compact abelian group and 
$N$ some  discrete, cocompact subgroup, i.e. $G/N$
is compact. 
\begin{lem}\label{LemBorelSectionsForG/N}
\begin{enumerate}
\item[$(i)$]The quotient map $G\to G/N$ has local sections.
\item[$(ii)$]The quotient map $G\to G/N$ has a Borel section.
\end{enumerate}
\end{lem}
\begin{pf}
$(i)$  
$N\subset G$ is discrete, i.e. there exists an open
neighbourhood $U$ of $0\in G$ such that $U\cap N=\{0\}.$
Let $+:G\times G\to G$ be the addition. $+$ is continuous 
so $+^{-1}(U)$ is an open neighbourhood of $(0,0)\in G\times G$.
So there is an open neighbourhood $V\subset G$ of $0\in G$ such
that $V\times V\subset +^{-1}(U)$. Let $W:=V\cap(-V)$.
Then $W$ is an open neighbourhood of $0\in G$, and for all
$x\in W$ and $n\in N\backslash\{0\}$ the sum  $x+n\notin W$, 
for $x\in W$ implies $-x\in W$ and in case $x+n\in W$ 
we would find $(x+n)+(-x)=n\in U$
-- a contradiction.
Therefore $W$ maps injectively to $G/N$, and as the quotient map
is open it defines a homoeomorphism from $W$
to its image $W/N$. 
This defines a local section from  $W/N$ to $G$,
and using addition in $G/N$ we can move $W/N$ 
all over $G/N$ to get a local section in the neighbourhood  
of each point in $G/N$.
\\\\
$(ii)$ 
This follows from compactness of $G/N$.
For let $s_m:U_m\to G$, $m=1,\dots,n$ be a family of local 
sections such that $\bigcup_{m=1}^{n} U_m=G/N$, then
$$
\sigma(z):=\begin{cases}
s_1(z),& {\rm if\ } z\in U_1,\\
s_2(z),& {\rm if\ } z\in U_2\backslash U_1,\\
\vdots&\vdots\\
s_n(z),& {\rm if\ } z\in U_n\backslash \bigcup_{m=1}^{n-1}U_m
\end{cases}
$$
defines a Borel section.
\end{pf}
The dual group of $G$ is \label{PageOfDualOfG}
$\widehat G:={\rm Hom}(G,\U(1))$.
With compact-open topology it becomes again a Hausdorff,
locally compact group, and if $G$ is second countable, then 
also $\widehat G$ is.
For parings of a group and its dual
we will use bracket notation \label{PageOfSkp}
$\langle\chi,g\rangle,\langle\alpha,n\rangle,\dots\in\U(1)$, for
$g\in G,\chi\in\widehat G,n\in N,\alpha\in\widehat N$.

We recall part of the classical duality theorems \cite{Ru}.
Pontrjagin Duality states that the canonical
map $G$ to $\widehat{\widehat G}$ is an isomorphism 
of topological groups.  Moreover, if  \label{PageOfDualOfN}
$N^{\perp}:=\{\chi\in\widehat G\ |\ \chi|_N=1\}$
is the annihilator of $N$, then
there is a canonical isomorphism 
$\widehat{G/N}\ni\alpha\mapsto\big( g\mapsto\langle\alpha,gN\rangle\big) \in N^{\perp}$, 
and by the same means ${\widehat G/N^\perp}\cong \widehat N$.
Further,  the dual group of a discrete group is compact and 
vice versa, so $N^\perp\subset\widehat G$ is a discrete 
cocompact subgroup, thus the situation is completely 
symmetric under exchange of $N,G$ by $N^\perp,\widehat G$. 
Let us denote the integration of a (compactly supported, 
continuous) function  $f:G\to\C$ against 
the Haar measure of $G$ simply  
by $\int_G f(g)\ dg$. 
For the Fourier transform $\hat f$ of $f$ 
 we use the convention 
$\hat f (\chi):=\int_G \langle \chi,g\rangle f(g)\ dg$, $\chi\in \widehat G$.
It extends to an isomorphism 
$L^2(G)\stackrel{\hat\ }{\to} L^2(\widehat G)$.
%

\subsection{Group and $\check{\rm C}$ech  Cohomology}
\label{SecCechandGroupCoho}
For a  topological $G$-module 
$M$ let us denote by 
$C^k_{\rm cont}(G,M)$ (resp.
$C^k_{\rm Bor}(G,M)$)
the continuous (resp. Borel)
maps $G^k\to M, k=0,1,2,\dots$
($C^0_?(G,M):=M$). They are topological spaces 
with the compact-open topology.
The differential $d:C^k_{\rm ?}(G,M)
\to C^{k+1}_{\rm ?}(G,M)$,  given
by
\begin{eqnarray*}
d f(g_1,\dots,g_{k+1}) &:=&
(-1)^{k+1}f(g_1,\dots,g_k)\\
&+&\sum_{i=1}^k (-1)^{i} f(g_1,\dots, g_i+g_{i+1},\dots,g_{k+1})\\
&+& g_1\cdot f(g_2,\dots,g_{k+1})
\end{eqnarray*}
makes $C^*_?(G,M)$ a cochain complex
with cohomology 
groups
$H^k_{\rm ?}(G,M):= 
Z^k_?(G,M)/B^k_?(G,M)$, for
the cocycles 
$Z^k_?(G,M):=\ker( d_k)$
and the boundaries \label{shllkashbf}
$B^k_?(G,M):={\rm im }(d_{k-1})$.
\\\\
Let $U_\bullet=\{U_i| i\in I\}$ be an open covering of
a space $B$. By $U_{i_0\dots i_n}$ we denote the
intersection $U_{i_0}\cap\dots\cap U_{i_n}$.\label{PageOfUij}
Let $\mathcal F$ be any abelian sheaf and let
$\check C^k(U_\bullet,\mathcal F):=\prod \mathcal F(U_{i_0\dots i_k})$.
The boundary operator 
$\delta:\check C^k(U_\bullet,\mathcal F)\to
\check C^{k+1}(U_\bullet,\mathcal F)$ is given by
$$
(\delta \varphi)_{i_0\dots i_{k+1}}:=
 \varphi_{i_1\dots i_{k+1}}|_{U_{i_0\dots i_{k+1}}}-
 \varphi_{i_0i_2\dots i_{k+1}}|_{U_{i_0\dots i_{k+1}}}+\dots
+(-1)^{k+1}  \varphi_{i_0\dots i_{k}}|_{U_{i_0\dots i_{k+1}}},
$$
and we use the standard notation for the cohomology groups 
$\check H^k(U_\bullet,\mathcal F)$.

We also use the notation $\underline A$ to denote
the locally constant sheaf of continuous functions
to $A$, for any abelian topological group $A$.

\subsection{The Unitary and the Projective Unitary Group} 
\label{SecTheUnitaryAndProjectiveUnitaryGroup}

Let $\HH$ be some 
infinite dimensional, separable
Hilbert space with unitary group $\U(\HH)$ which we
equip with the strong (or equivalently weak) operator topology.
We denote by \label{ThePageOfPUAdAndU}
$${\rm Ad}:\U(\HH)\to\PU(\HH):=\U(\HH)/ \U(1)
$$
the quotient map from the unitary  onto the projective 
unitary group, and we endow the latter with the quotient
topology. 
\\\\
Let $U\subset \U(\HH)$ be any subset, e.g. $U=\U(1)\cdot \Eins$
and let $M$ be some measure space, e.g. $G$ or $G/N$,
then we denote by
\begin{equation}\label{EqTheDefiOfLinfty}
L^\infty(M,U)\subset \U(L^2(M)\otimes\HH)
\end{equation}
the set of unitary operators which are 
given by (equivalence classes of) measurable functions 
$f:M\to U$ which act as multiplication operators.
In particular $L^\infty(M,U)$ has the subspace topology
of $\U(L^2(M)\otimes \HH)$ which is usually referred as
the weak topology on $L^\infty(M,U)$.
The image of $L^\infty(M,U)$ under 
${\rm Ad}:\U(L^2(M)\otimes \HH)\to \PU(L^2(M)\otimes\HH)$
is denoted by $\PP L^\infty(M,U)$.

Let us equip the Borel functions ${\rm Bor}(G,U)$
with the compact-open topology.

\begin{lem}\label{LemTheInclusionIsConti}
The natural map
$$
{\rm Bor}(G,U)\to L^\infty(G,U)
$$
 is continuous.
\end{lem}
\begin{pf}
A standard $\varepsilon/3$-argument.
Let $f_\alpha\to f$ be a converging net in
${\rm Bor}(G,U),$
i.e. for each compact $K\subset G$ and 
any $w\in \HH$ 
we have 
$$
\|f_\alpha-f\|_{_{K,w}}:=\sup_{g\in K}\|f_\alpha (g)w-f(g)w\|_{_{\HH}}\to 0.
$$
We have to show that for all $v\in L^2(G,\HH)$ 
$$
\|f_\alpha v- fv\|_{_{L^2(G,\HH)}}\to 0.
$$
Recall that in case of the Haar measure
the compactly supported functions
$C_c(G)$ are dense in $L^2(G)$.
So for any $v\in L^2(G,\HH)$ and $\varepsilon >0$
there exist 
compactly supported functions $h_i\in C_c(G),$
vectors  $w_j\in\HH$ and numbers $a_{ij}\in\C$
such that 
$\|v-\sum_{i,j=1}^N a_{ij} h_i\otimes w_j\|_{_{L^2(G,\HH)}} <\varepsilon/3.$

Choose $K:=\bigcup_{i=1}^N {\rm supp}\ h_i \subset G$ and
$C:=1+\sum_{i,j=1}^N |a_{ij}|^2\int_K |h_i(g)|^2\ dg >0.$
Then there exists an $\alpha_0$ such that
$\|f_\alpha-f\|_{K,w_j}<\frac{\varepsilon}{3 N\sqrt{ C}}$, for
all $\alpha>\alpha_0$ and all $j=1,\dots,N$.

We can estimate now
\begin{eqnarray*}
&&\|f_\alpha v- f v\|_{_{L^2(G,\HH)}}
\le\ 
 \|f_\alpha (v - \sum^N_{i,j=1} a_{ij} h_i\otimes w_j)\|_{_{L^2(G,\HH)}} 
\\&&\qquad
 +
 \|(f_\alpha - f)(\sum^N_{i,j=1} a_{ij} h_i\otimes w_j) \|_{_{L^2(G,\HH)}}
+
\|f(\sum^N_{i,j=1} a_{ij} h_i\otimes w_j - v)\|_{_{L^2(G,\HH)}}\\
&\le&\|f_\alpha\|_{_{\rm Op}}\cdot \varepsilon/3
 +
 \Big( \int_K \|\sum_{i,j=1}^N a_{ij}  h_i(g) (f_\alpha(g)w_j- f(g) w_j)\|_{_{\HH}}^2
\ dg\Big)^{1/2}
 +
 \|f\|_{_{\rm Op}}\cdot \varepsilon/3\\
&\le\ & \varepsilon/3
+\Big( N^2\sum_{i,j=1}^N  |a_{ij}|^2\int_K   |h_i(g)|^2 \| (f_\alpha(g)w_j- f(g) w_j)\|_{_{\HH}}^2
\ dg\Big)^{1/2}
 + \varepsilon/3\\
\\
&\le& \varepsilon/3
+\Big( N^2\sum_{i,j=1}^N  |a_{ij}|^2\int_K   |h_i(g)|^2 \frac{\varepsilon^2}{9 N^2 C}
\ dg\Big)^{1/2}
 + \varepsilon/3\\
 \\
 &\le &\varepsilon/3+\sqrt{\varepsilon^2/9}+\varepsilon/3 = \varepsilon,
\end{eqnarray*}
for all $\alpha>\alpha_0$.
\end{pf}

For $\PU(\HH)$-principal bundles we have the
following well-known classification theorem.
(See e.g.  \cite[Thm. 10.8.4]{Di} for the first  and
of \cite{PR} for the second statement; although therein
it is not stated as below, we can carry over the proofs.) 
To state the theorem we introduce some notation.
Let us denote by ${\rm Iso}(E)$ the set
of isomorphism classes of $\PU(\HH)$-principal bundles 
over $E$, and if $P\to E$ is a $\PU(\HH)$-principal bundle, we denote
by ${\rm Aut}_0(P,E)$ the group of bundle automorphisms
(over the identity of $E$). 
There is the subgroup 
${\rm Null}(P,E)\subset{\rm Aut}_0(P,E)$ which consists of all null-homotopic 
bundle automorphisms.

\begin{thm}\label{ThmClassifacationOfPUBundlesAndAutos}
Let $E$ be a paracompact Hausdorff space and 
$P\to E$ be any fixed $\PU(\HH)$-principal bundle.
\begin{enumerate}
\item 
$\PU(\HH)$-bundles over $E$ are classified by 
the second $\check{\rm C}$ech cohomology
with values in the locally constant sheaf of 
continuous functions to $\U(1)$,
$${\rm Iso}(E)\cong \check H^2(E,\underline{U(1)}).
$$
\item
There is a short exact sequence
$$
0\to {\rm Null}(P,E)\to {\rm Aut}_0(P,E)\to \check H^1(E,\underline{\U(1)})\to 0.
$$
\end{enumerate}
\end{thm}

In particular, if $P\to E$ is the trivial bundle $P=E\times \PU(\HH)$
we identify ${\rm Aut}_0(P,E)$
with the continuous functions $C(E,\PU(\HH))$.
If we keep this in mind,
a corollary of the classification theorem of bundle automorphisms is 
the following statement which is sometimes 
called  homotopy commutativity of the projective unitary group.

\begin{cor}\label{CorTheHomotopyCommutativityOfPU}
Let $E$ be a paracompact Hausdorff space and
$f,g:E\to \PU(\HH)$ two continuous functions. Then there
extists a continuous function $V:E\to \U(\HH)$ such that
$$
f(x)\ g(x)=g(x)\ f(x)\ {\rm Ad}(V(x)),\qquad x\in E.
$$
\end{cor}

\begin{pf}
$x\mapsto f(x)g(x)$ and 
$x\mapsto g(x)f(x)$ define the same
$\check{\rm C}$ech class.
\end{pf}

We do not give a proof  of Theorem \ref{ThmClassifacationOfPUBundlesAndAutos},
but we remark that it depends heavily on the the fact that the unitary group $\U(\HH)$ is contractible.
For the strong topology on $\U(\HH)$ this is not difficult 
to prove and may be found in \cite[10.8]{Di}.
We line out the proof for  convenience.
Assume $\HH=L^2([0,1])$, then let
$\varphi_t:L^2(0,t)\cong L^2(0,1)$ be the  
isometric isomorphism defined by
$\varphi_t(f)(x):=\sqrt t f(tx)$, for $t>0$.
We define $H:[0,1]\times\U(L^2(0,1))\to\U(L^2(0,1))$ 
by $H(0,U):=\Eins$ and
$$
(H(t,U)(f))(x):=
\begin{cases}
\big(\varphi_t^{-1}\circ U\circ\varphi_t(f|_{(0,t)})\big)(x),& {\rm if\ } 0<x<t,\\
f(x),& {\rm if\ } t<x<1 ,
\end{cases}
\quad{\rm for\ } t>0,
$$
then $H$ is a homotopy connecting the identity on $\U(\HH)$ and
the constant function with value $\Eins\in\U(\HH)$, thus $H$ is a 
contraction.

\begin{rem}\label{RemPUEquivariantContractionOfU}
For each fixed time slice $t\in [0,1]$ $H(t,\ .\ )$ is
a group homomorphism. So $\U(\HH)$ is contractible 
as group 
\end{rem}

\begin{rem}\label{RemAbelianU}
Theorem \ref{ThmClassifacationOfPUBundlesAndAutos}
also  holds when we
replace the unitary group $\U(\HH)$ by any contractible, abelian  (sub)group
$\UA(\HH)$ such that  $\U(1)\subset\UA(\HH)$.
I.e. the second $\check{\rm C}$ech cohomology $\check H^2(E,\underline{\U(1)})$
also classifies  $\PP\UA(\HH)$-bundles over paracompact Hausdorff spaces $E$
and the first  $\check{\rm C}$ech cohomology $\check H^1(E,\underline{\U(1)})$ also classifies 
$\PP\UA(\HH)$-bundle automorphisms, for $\PP\UA(\HH):=\UA(\HH)/\U(1)$.
\end{rem}
\label{PageOfUA}

The next lemma shows that there exists
a appropriate commutative 
version $\U(1)\to \UA(\HH)\to\PP\UA(\HH)$ of
 $\U(1)\to \U(\HH)\to\PU(\HH)$.
\begin{lem}
There exists a contractible, commutative 
subgroup $\UA(\HH)\subset \U(\HH)$ 
such that $\U(1)\cdot\Eins\subset\UA(\HH)$.
\end{lem}
\begin{pf}
We let $\HH=L^2(0,1)$ and $\UA(\HH):= L^\infty([0,1],\U(1))$.
It is a  commutative subgroup and $\U(1)\cdot\Eins\subset \UA(\HH)$.
The contraction $H$ from above restricts to
$h:=H|_{[0,1]\times \UA(\HH)}:[0,1]\times \UA(\HH)\to \UA(\HH)$
and is given by  
$$
h(t,g)(s):=\begin{cases}
g(\frac{1}{t}s),& {\rm if\ } s<t,\\
1,& {\rm if\ }s>t
\end{cases}.
$$
\end{pf}

\subsection{Crossed Products}
\label{SecTheSecOnCrossedProducts}
Let $(A,G,\alpha)$ be a $C^*$-dynamical system,
i.e. $A$ is some $C^*$-algebra equipped with 
a (strongly) continuous action $\alpha:G\to {\rm Aut}(A).$
Let $(\pi,\HH)$ be some faithful representation of
$A$. 
We embed $C_c(G,A)$, the 
compactly supported continuous functions $G\to A$,
 into $\LL(L^2(G,\HH))$:
define 
$$f\times\_ :L^2(G,\HH)\to L^2(G,\HH)$$
by
$$(f\times F)(g):=\int_G \pi(\alpha_{-g}(f(h))) F(g-h)\ dh,$$
for $f\in C_c(G,A), F\in L^2(G,\HH).$
The adjoint operator of $f\times\_$ 
 is easily calculated and is 
$(f\times\_ )^*={f^\times\times\_}$, wherein
$f^\times(g):=\alpha_{g}(f(-g))^*$.
Clearly $f^\times$ has compact support,
and $f^\times$ is continuous (all $\alpha_g$ have
norm $1$).
So   $C_c(G,A)$
 is closed under the $^\times$-operation.
Furthermore one has\footnote{We use the same 
symbol $\times$ for different maps.}
$f_1\times (f_2 \times F)=(f_1\times f_2)\times F$,
wherein 
$$
(f_1\times f_2)(g):=\int_G f_1(h) \alpha_h(f_2(g-h))\ dh
$$
again defines an element of $C_c(G,A)$.
Thus $C_c(G,A)\hookrightarrow\LL(L^2(G,\HH))$
is a $^*$-subalgebra.
The {\bf crossed product} of $G$ and $A$ is then
defined as the norm completion 
$$G\times_\alpha A:=\bigg(\overline{C_c(G,A)}^{\|.\|},\times,^\times\bigg).
$$
(This is in fact well-defined, since the operator norm of 
$f\times\_$ is independent of the faithful representation 
$(\pi,\HH)$.)
\\\\
For $\chi\in\widehat G, f\in C_c(G,A)$ we set 
$\hat\alpha_\chi( f)(g):=\langle\chi,g\rangle f(g)$ which 
extends to a strongly continuous action
$\hat \alpha:\widehat G\to{\rm Aut}(G\times_\alpha A)$,
so $(G\times_\alpha A, \widehat G, \hat \alpha)$
again defines a $C^*$-dynamical system.
Going once more through the process of building the crossed product
 gives a $C^*$-dynamical system
 $(\widehat G\times_{\hat\alpha}(G\times_\alpha A), G, \hat{\hat{\alpha}})$,
and
a  key statement in the analysis of crossed products
is the following Takai Duality Theorem (see e.g. \cite{Pe}).
\begin{thm}\label{ThmTakaiduality}
There is an isomorphism of 
$C^*$-dynamical systems
$$\big(\widehat G\times_{\hat\alpha}(G\times_\alpha A\big), G, \hat{\hat{\alpha}})
\cong \big(A\otimes_{C^*} \K(L^2(G)), G, \alpha\otimes {\rm Ad}\circ\varrho\big),
$$
for the right regular representation  $\varrho:G\to \U(L^2(G))$,
i.e. $(\varrho_gf)(h):=f(h+g)$ and  
${\rm Ad}(\varrho_g)(K)=\varrho_g K\varrho_g^{-1},$ for $K\in \K(L^2(G))$.
\end{thm}

\subsection{Some Topology}

We will sometimes use the word space as 
abbreviation for topological space. Let $X,Y$ 
be topological spaces. By  ${\rm Bor}(X,Y)$ 
we denote the set of Borel functions from $X$
to $Y$. We endow this space with the compact-open 
topology, i.e. we define a basis of the topology
by all sets of the form $U_{K,V}:=\{ f:X\to Y| f(K)\subset V\}$ for
compact $K\subset X$ and open $V\subset Y$.
The subspace of continuous maps will be denoted
by  ${\rm Map}(X,Y)\subset {\rm Bor}(X,Y)$. \label{PageOfMapAndBor}
We use the notation  $C(X,Y)$ for the continuous functions 
  if we do not want to specify the topology on  it, or in case 
$Y$ is a normed space, then we put  the 
supremums norm on $C(X,Y)$.
Recall the exponential law \cite{Sch} for ${\rm Map}(.,..)$.
\begin{lem}\label{LemExpLawForMap}
Let $X,Y,Z$ be topological spaces.
Assume $X$ and $Y$ to be Hausdorff  and 
$Y$ locally compact. Then we have a homoeomorphism
\begin{eqnarray*}
{\rm Map}(X\times Y,Z)&\cong&{\rm Map}(X,{\rm Map}(Y,Z))\\
f&\mapsto& \big(x\mapsto f(x,\_)\big).
\end{eqnarray*}
\end{lem}

Let $E\to B$ be a surjective fibration and
assume $Y$ is locally compact and
Hausdorff. By use of the 
exponential law it is then immediate
that ${\rm Map}(Y,E)\to {\rm Map}(Y,B)$
still has the homotopy lifting property 
with respect to all Hausdorff spaces;
we denote this property by $T_2HLP$.
For contractible $E$ the map
${\rm Map}(Y,E)\to {\rm Map_0}(Y,B):=
\{f\simeq {\rm const.}\}$
becomes  a surjection having the $T_2HLP$.

We now restrict ourselves to the particular
 case of $E=\U(\HH), B=\PU(\HH).$
Since $\U(\HH)$  is Polish  $\cite{Ke}$ the quotient
map $\U(\HH)\to\PU(\HH)$ admits a Borel
section $s:\PU(\HH)\to\U(\HH)$ \cite[\rm Thm. 12.7]{Ke}.

Let $f\in {\rm Map}(Y,\PU(\HH))$ be a
map not homotopic to a constant map.
 The Borel function 
$s_*f=s\circ f:Y\to\U(\HH)$ gives rise to a bijection of sets
\begin{eqnarray*}
{\rm Map}(Y,\U(\HH))&\to& {\rm Map}(Y,\U(\HH))\cdot s_*f.\\
g&\mapsto& g\cdot s_*f
\end{eqnarray*}
We turn this map into a homoeomorphism
by defining the topology of the right hand side to
be the image of the topology of the left hand side.
Then
$$
{\rm Map}(Y,\U(\HH))\cdot s_*f\to {\rm Map_0}(Y,\PU(\HH))\cdot f.\\
$$
as well as the assembled 
map
$$\coprod_{f\in\pi} \big( {\rm Map}(Y,\U(\HH))\cdot s_*f\big)\to 
\coprod_{f\in\pi}\big( {\rm Map_0}(Y,\PU(\HH))\cdot f\big)
$$
becomes  a surjection satisfying
the $T_2HLP$.
Here $\pi:=\{f: Y\to \PU(\HH)\}$ is a set of 
representatives of the first
homotopy group $\pi_0({\rm Map}(Y,\PU(\HH)))$.
If $Y$ is compact ${\rm Map_0}(Y,\PU(\HH))\subset
{\rm Map}(Y,\PU(\HH))$ is open, i.e.
$${\rm Map}(Y,\PU(\HH))=\coprod_{f\in\pi} \big({\rm Map_0}(Y,\PU(\HH))\cdot f\big).$$
 
To consider non compact $Y$
the following easy statement is 
helpful.
\begin{lem}
Let $Z$ be any (pointed) space, and let
$Z_z$ denote the path-connected 
component of $z\in Z$. Then
$$\coprod_{[z]\in\pi_0(Z)}Z_z\to \bigcup_{[z]\in\pi_0(Z)}Z_z=Z$$
is a fibration.
\end{lem}
\begin{pf}
Give any test space $X$ and a diagram 
$$
\xymatrix{
X\ar[r]^{h_0}\ar[d]^{i_0}&\coprod Z_z\ar[d]\\
X\times I\ar[r]^{h}&Z
}$$
we  define the (unique)  homotopy $\tilde h:X\times I\to\coprod Z_z$ that 
fits into the above diagram simply by $\tilde h(x,t):=h(x,t)$. One has to
check that $\tilde h$ is continuous. For this it 
is sufficient to show that $\tilde h^{-1}(Z_z)$ is open. 
Since $h$ maps path-connected components into
path-connected components we have $h^{-1}(Z_z)=h_0^{-1}(Z_z)\times I$.
Since
$h_0^{-1}(Z_z) \subset X$ is open,
it follows that $\tilde h^{-1}(Z_z)=h^{-1}(Z_z)=h_0^{-1}(Z_z) \times I\subset
X\times I$ is open.
\end{pf}

This lemma illustrates that
\begin{equation}\label{EqFwjbfuklargf}
\coprod_{f\in\pi} \big( {\rm Map}(Y,\U(\HH))\cdot s_*f\big)\to {\rm Map}(Y,\PU(\HH))
\end{equation}
is  surjective and satisfies  the $T_2HLP$ for all locally compact Hausdorff
spaces $Y$, e.g. $Y=G\times G/N$. This implies the next

\begin{cor}\label{CorExistenceOfBorelCocycle}
Let $U$ be a contractible Hausdorff space, and let
$\mu:U\to {\rm Map}(G\times G/N,\PU(\HH))$ be continuous.
Then there exists a continuous lift
$\overline\mu: U\to {\rm Bor}(G\times G/N,\U(\HH))$ such
that ${\rm Ad}_*\circ\overline\mu=\mu$, i.e. 
${\rm Ad}\circ\overline\mu(u)=\mu(u)$ for all $u\in U$.
\end{cor}

\begin{pf}
(\ref{EqFwjbfuklargf}) has the $T_2HLP$ and $U$ is contractible and 
Hausdorff.
\end{pf}

For the remainder of this paragraph we
stick to compact $Y$, e.g. $Y=G/N$. 
We shall  consider the map
${\rm Ad}_*:{\rm Bor}(Y,\U(\HH))\to {\rm Bor}(Y,\PU(\HH))$.
\begin{lem}
${\rm Ad}_*$ is continuous and open.
\end{lem}
\begin{pf}
Note first that ${\rm Ad}:\U(\HH)\to \PU(\HH)$ is continuous and open.
\\
Continuity:
Let $U_{K,V}:=\{g | g(K)\subset V, K\subset Y\ {\rm compact}, 
V\subset\PU(\HH)\ {\rm open}\}$. Then
${\rm Ad}_*^{-1}(U_{K,V})=\{f | f(K)\subset {\rm Ad}^{-1}(V)\}$  is
open in ${\rm Bor}(Y,\U(\HH))$.
\\
Openness:
Let $U_{K,W}:=\{f | f(K)\subset W, K\subset Y\ {\rm compact}, 
W\subset\U(\HH)\ {\rm open}\}$. Then the inclusion
${\rm Ad}_*(U_{K,W})\subset \{g |  g(K)\subset {\rm Ad}(W)\}$  is
obvious; we show equality.
To do so it is sufficient to construct a Borel section $s$ of
$\U(\HH)\to\PU(\HH)$ such that $s({\rm Ad}(W))\subset W$.
$\PU(\HH)$ is separable. Take a countable dense set $\{ x_i\}_i\in {\N}$
and a local trivialisation $U_0\subset \PU(\HH)$ of
$\U(\HH)\to\PU(\HH)$.
Let  $U_i:= U_0 x_i \subset\PU(\HH)$ and 
$V_i:={\rm Ad}(W)\cap U_i$,
so $\bigcup_i V_i={\rm Ad}(W)$.
It suffices to construct  Borel
sections $s_i:V_i\to W$ and puzzling them together
by
$$
s(z):=\begin{cases}
s_1(z),& {\rm if\ } z\in V_1,\\
s_2(z),& {\rm if\ } z\in V_2\backslash V_1,\\
\vdots&\vdots\\
s_n(z),& {\rm if\ } z\in V_n\backslash \bigcup_{m=1}^{n-1}V_m\\
\vdots&\vdots
\end{cases},
$$
then $s$ is Borel, since we put together a countable family.
The sections $s_i:V_i\to W$ may be obtained by similar manners:
Let $h_i:{\rm Ad}^{-1}(V_i)\to V_i\times\U(1)$ be a trivialisation, and
let $W_i:= h_i({\rm Ad}^{-1}(V_i)\cap W)$.
Then $W_i\subset V_i\times \U(1)$ is open.
Thus for each 
$l\in L:=\{l\in \N | x_l\in V_i\}$ there is $\varphi_l\in \U(1)$ and
an open neighbourhood $V_i^l\ni x_l$ such
that $V_i^l\times\{\varphi_l\}\subset W_i$.
We define
$$
s_i(z):=\begin{cases}
h_i^{-1}(z,\varphi_{l_1}),& {\rm if\ } z\in V_i^{l_1},\\
h_i^{-1}(z,\varphi_{l_2}),& {\rm if\ } z\in V_i^{l_1}\backslash V_i^{l_2},\\
\vdots&\vdots\\
h_i^{-1}(z,\varphi_{l_n}),& {\rm if\ } z\in V_i^{l_n}\backslash \bigcup_{m=1}^{n-1}V_i^{l_m}\\
\vdots&\vdots
\end{cases},
$$
for a counting $l_1,l_2,\dots$ of $L$.
This completes the proof.
\end{pf}

\begin{lem}\label{LemMapIsBaseOfPFB}
 If $Y$ is compact then
 ${\rm Ad}_*:{\rm Bor}(Y,\U(\HH))\to {\rm Bor}(Y,\PU(\HH))$ 
 is a locally trivial
 principal fibre bundle with structure group ${\rm Bor}(Y,\U(1))$.
\end{lem}

\begin{pf}
 Let $t:V\to \U(\HH)$ be a local section 
 of $\U(\HH)\to\PU(\HH)$ for some open
 $V\subset\PU(\HH).$
 Then $t_*: U_V
 \ni f\mapsto t\circ f\in {\rm Bor}(Y,\U(\HH))$
 is a local section on $U_V:=\{f:Y\to \PU(\HH) | f(Y)\subset V\}$
 which is open, since $Y$ is compact.
 We can cover the whole of ${\rm Bor}(Y,\PU(\HH))$ by
 translates of $U_V$ under the action of
 ${\rm Bor}(Y,\PU(\HH))$ on itself.
 The lemma will be proven if we can show that
 ${\rm Bor}(Y,\PU(\HH))$ and 
 ${\rm Bor}(Y,\U(\HH))/{\rm Bor}(Y,\U(1))$ are homoeomorphic.
 Take a Borel section $\sigma:\PU(\HH)\to \U(\HH)$. Then 
 \begin{eqnarray*}
 \sigma_*:{\rm Bor}(Y,\PU(\HH))&\to&{\rm Bor}(Y,\U(\HH))/{\rm Bor}(Y,\U(1))\\
 f&\mapsto&[\sigma\circ f]
 \end{eqnarray*}
 is easily seen to be a bijection such that
 $$\xymatrix{
 {\rm Bor}(Y,\U(\HH))\ar[d]_{{\rm Ad}_*}\ar[rd]&\\
 {\rm Bor}(Y,\PU(\HH))\ar[r]^{\sigma_*\quad\quad\quad}&
 {\rm Bor}(Y,\U(\HH))/{\rm Bor}(Y,\U(1))
 }
 $$
 commutes. It follows that  $\sigma_*$ is a
 homoeomorphism, for  the quotient map
 and ${\rm Ad}_*$ are both continuous  and open.
\end{pf}

We will use the above lemma in combination with the next.

\begin{lem}\label{LemLiftForLocTrivFB}
 Let $P\to M$ be a locally trivial fibre bundle, and 
 $B$ a paracompact  space. Then for each covering 
 $\{U_i | i\in I\}$ of $B$ and maps $\zeta_{ij}:U_i\cap U_j\to M$,
 there exists a refinement
 $\{V_{ix} | ix\in I\times B\}$ $(V_{ix}\subset U_i)$ and continuous 
 $\overline\zeta_{ix,jy}:V_{ix}\cap V_{jy}\to P$ such that the diagram
$$
\xymatrix{
&&P\ar[d]\\
V_{ix}\cap V_{jy}
\ar[urr]^{\overline\zeta_{ix,jy}}
\ar[rr]_{\ \ \ \ \ \ \ \ \ \zeta_{ij}|_{ V_{ix} \cap V_{jy} } }
&&M
}
$$
 commutes.
\end{lem}
\begin{pf}
 Without restriction we can assume that the 
 covering $\{U_i\}_{i\in I}$ is locally finite, so for each $x\in X$
$I_x:=\{ k \in I | \{x\}\cap U_k\not= \emptyset\}$
 is finite. 
 If $x\notin U_i$, we define $V_{ix}:=\emptyset.$
  If $x\in U_i$, then $i\in I_x$.
 For each $\zeta_{ik}(x), k\in I_x,$ choose a local trivialisation
 $M_{i,k,x}\ni \zeta_{ik}(x)$  of $P\to M$.
Then 
 $V_{ix}:=\bigcap_{k \in I_x} \zeta_{ik}^{-1}(M_{i,k,x})\subset U_i$
 is open.
By construction the image of 
$\zeta_{ij}|_{V_{ix}\cap V_{jy}}$
is contained in $\bigcap_{k\in I_x}M_{i,k,x}$,
and we can compose $\zeta_{ij}|_{V_{ix}\cap V_{jy}}$
with any local section, say $M_{i,i,x}\to P$, to define $\overline\zeta_{ix,jy}$.
\end{pf}

\end{appendix}

\newpage

\begin{theindex}
\addcontentsline{toc}{section}{$\bullet\ \ \!$   Index}

\begin{supertabular} {lr}
Symbol & Page\\ 
\\
$\langle\ .\ ,\ . \ \rangle$&\pageref{PageOfSkp}\\
$[\ .\ ,\ . \ ]$&\pageref{PropTheClassifyingSpaceOfPairs},
\pageref{ThmOfBSOnTheRepresentability}\\
\\
$\A $ &\pageref{SubsecThelocalStrOfPairsAndTheisClassifying} \\
$\A_1 $ &\pageref{PageOfA1EA1BA1} \\
$A_{ji} $ &\pageref{PageOfAandBeta} \\
Ad &\pageref{ThePageOfPUAdAndU}\\
$\alpha_{ji}$&\pageref{EqTheDefiOfAlpha}\\
$\alpha^\rho$&\pageref{PageOfCStarThings}\\
$\alpha^\mu$&\pageref{EqActOnFct}\\
${\rm Aut}, {\rm Aut}_0, {\rm Aut}_1$&\pageref{PageofAut}\\
$B$&\pageref{DefiOfBaseSpace}\\
$\B\A_1 $ &\pageref{PageOfA1EA1BA1} \\
$\beta_{ji} $ &\pageref{PageOfAandBeta} \\
$\delta_g $ &\pageref{SecPairsAndTwistedCech} \\
${\rm Bor}$&\pageref{PageOfMapAndBor}\\
$c_i$&\pageref{DefiOfExteriorEquivalence}\\
$c^\tau$&\pageref{PageOfCTau}\\
$\chi,\chi_{ji}$&\pageref{EqTheDefiOfMi}\\
Dyn&\pageref{PageOfDyn}\\
$\widehat{\rm Dyn}$&\pageref{PageOfDualFunctors}\\
${\rm  Dyn}^\dag$&\pageref{PageOfDualisableFunctor}\\
$\delta$&\pageref{ThmTheSpecialCaseOfRZ}\\
$\delta^{\rm s},\delta^{\rm as}, \delta^{\rm im}$&
\pageref{ThmTheMainTheoremOfThisThesisNumberThree}\\
$E$&\pageref{DefiOfBaseSpace}\\
$\widehat E$&\pageref{PageOfAllTheDuals}, \pageref{PageDualOfEanfFdd}\\
$E_{\rm dyn}$&\pageref{PageOfPDynStbilll}\\
$\E\A_1 $ &\pageref{PageOfA1EA1BA1} \\
$F$&\pageref{DefiOfBaseSpace}\\
$\mathcal F$&\pageref{PageFouriertrrrr}, \pageref{PageMoerFoureirere}\\
$\hat{\mathcal F}$&\pageref{PageOfInvFourier}\\
$f^{\overline\mu}$&\pageref{EqTheTOperator}\\
$G, G/N$&\pageref{DefiOfBaseSpace}\\
$\widehat G, \widehat G/N^\perp$&\pageref{PageOfDualOfG}\\
$\Gamma,\Gamma_0$&\pageref{DefiOfBaseSpace}, \pageref{PageOfCStarThings}\\
$\gamma_{ji}$&\pageref{PageOfUA}\\
$g_{ji}$ & \pageref{PageOfTransiFunc}\\
$\hat g_{ji}$ & \pageref{PageThePageOfHatGji}\\
$\hat g^{\rm dyn}_{ji}$&\pageref{PageOfPhiDyn}\\
$\HH$&\pageref{DefiOfBaseSpace}\\
$I$&\pageref{DiagTheLocalStructureOfAPair}\\
$\kappa$&\pageref{DiagAnAlmostTDualityDiagram}\\
$\kappa_i$&\pageref{PageOfKappai}\\
$\kappa_i^a,\kappa^b_i$&\pageref{PageOfPrimedTransis}\\
$\kappa^\sigma,\overline\kappa^\sigma $ & \pageref{LemTheDefiOfThePoincareClass}\\
$\overline{\hat\kappa}^{\hat\sigma} $ & \pageref{EqTheDualAnalyticalExpression}\\
$\kappa^{\rm top}$&\pageref{ThmTheTopologicalisationMapTau}\\
$\kappa_i^{\rm top}$&\pageref{EqTheLocalDefinitionOfKappaTop}\\
$\K(\HH) $ & \pageref{PageOfCompacts}\\
$L^\infty $ & \pageref{EqTheDefiOfLinfty}\\
$l_i,l_i^0,\lambda_i$&\pageref{PageOfLambdaili}\\
$\lambda_{_{G}}$&\pageref{PageTheExampleOfStableEquivalence}\\
$\lambda_{_{G/N}}$&\pageref{EqTheDualPUCocycle}\\
$m_i$&\pageref{EqInClaim1Forbbbbr}\\
${\rm Map}$&\pageref{PageOfMapAndBor}\\
$\mu_i$&\pageref{EqLift}\\
$\overline\mu_{i}$&\pageref{EqThisIsTheDefiOfPhiOhYear}\\
$\hat\mu_i$&\pageref{EqTheDualDeckerCocycle}\\
$\mu_i^{\rm dyn}$&\pageref{pageOfDefiOfDualisableDecker}\\
$N$&\pageref{DefiOfBaseSpace}\\
$N^\perp$&\pageref{PageOfDualOfN}\\
$\omega,\omega_{i}$&\pageref{EqThisIsTheDefiOfPhiOhYear}\\
$P, (P,E)$&\pageref{DefiOfBaseSpace}\\
$\mathcal P$& \pageref{PageOfTopImAndP}\\
$\widehat P$&\pageref{PageOfAllTheDuals}, \pageref{PageDualOfEanfFdd}\\
$P_{\rm top}$&\pageref{PageOfStableP}\\
$P_{\rm dyn}$&\pageref{PageOfPDynStbilll}\\
$P_\HH$&\pageref{EqStablilisedPUBundleP}\\
${\rm Par}$&\pageref{EqStablilisedPUBundleP}\\
$\widehat{\rm Par}$&\pageref{PageOfDualFunctors}\\
$\PU(\HH)$&\pageref{ThePageOfPUAdAndU}\\
$\pi$&\pageref{DefiOfPoincareClass}\\
$\varphi, \varphi_{kji}$&\pageref{EqThisIsVarphiTheObstructionCocycleAgainstDeckers}\\
$\varphi^*P$&\pageref{DefiOfBaseSpace}\\
$\phi,\phi_{ji}$&\pageref{EqThisIsTheDefiOfPhiOhYear}\\
$\hat\phi,\hat\phi_{ji}$&\pageref{PageOfDualPhiTotalCocyc}\\
$\phi^{\rm dyn}_{ji}$&\pageref{PageOfPhiDyn}\\
$\psi,\psi_{kji}$&\pageref{PageOfPsi},\pageref{EqThisIsTheDefiOfPhiOhYear}\\
$q, Q_{_G},Q_{_{G/N}}$&\pageref{DiagOfExatSequences}\\
$\rho$&\pageref{SubsecDynamicalTriples}\\
$\rho^\tau$&\pageref{EqLift}\\
$\hat \rho$&\pageref{PageOfAllTheDuals}\\
$\rho^{\rm dyn}$&\pageref{pageOfDefiOfDualisableDecker}\\
$\sigma$&\pageref{LemTheDefiOfThePoincareClass}\\
$\hat\sigma$&\pageref{EqTheDualAnalyticalExpression}\\
$\T^n,\hat\T^n(B)$&\pageref{PageOfTorus}\\
$T_{\overline\mu}$&\pageref{EqTheTOperator}\\
Top& \pageref{PageOfTopFunc}\\
${\rm Top^s},{\rm Top^{as}}$& \pageref{PageOfAllTheTops}\\
${\rm Top^{im}}$& \pageref{PageOfTopImAndP}\\
$\tau,\tau(B)$&\pageref{ThmTheTopologicalisationMapTau}\\
$\nou$ & \pageref{PageOfNoU}\\
$\U(\HH)$&\pageref{ThePageOfPUAdAndU}\\
$\UA(\HH)$&\pageref{RemAbelianU}\\
$U_\bullet, U_i,U_{ji},\dots$&\pageref{DiagTheLocalStructureOfAPair}, \pageref{PageOfUij}\\
$w_{ji}$&\pageref{PageOfUA}\\
$Z^k_{\rm  cont},Z^k_{\rm Bor}, \check Z^k$&\pageref{shllkashbf}\\
$\zeta_{ji}$ & \pageref{PageOfTransiFunc}\\
$\overline\zeta_{ji}$ & \pageref{PageOfPsi}\\
$\hat\zeta_{ji}$ & \pageref{EqTheDualPUCocycle}\\
$\zeta'_{ji},\hat\zeta'_{ji}$&\pageref{PageOfPrimedTransis}\\
$\zeta_{ji}^{\rm dyn}$&\pageref{pageOfDefiOfDualisableDecker}\\
$\widehat{\zeta^{\rm dyn}_{ji}}$&\pageref{PageOfPhiDyn}\\
\end{supertabular}

\end{theindex}

\end{document}